\input amstex
\documentstyle{amsppt} 
\refstyle{C}
\TagsOnRight
\NoRunningHeads


\def\eqnum#1{\eqno (#1)}

\def\fnote#1{\footnote}

\def\ord{\text {\rm ord}}

\def\wt{\text {\rm wt}}
\def\XOR{\operatorname{\text {\rm XOR}}}
\def\AND{\operatorname{\text {\rm AND}}}
\def\OR{\operatorname{\text {\rm OR}}}
\def\NEG{\operatorname{\text {\rm NEG}}}


\topmatter
\title Uniformly distributed sequences of 
$p$--adic integers, II \endtitle
\author Vladimir Anashin\endauthor
\address
Faculty for the Information Security. 
Russian State University for the Humanities.
Miusskaya Square, 6, Moscow 125267, Russia.
\endaddress
\email vladimir\@anashin.msk.su, anashin\@rsuh.ru\endemail


\abstract
The paper describes ergodic (with respect to the Haar measure) functions 
in the class of all functions, which are defined on (and take values in) 
the ring $\Bbb Z_p$ of $p$-adic integers, and which satisfy (at least, locally) Lipschitz 
condition with coefficient 1. Equiprobable (in particular, measure-preserving) 
functions of this class are described also. In 
some cases (and especially for $p=2$) the descriptions are given by explicit formulae. 
Some of the results may be viewed as descriptions of ergodic isometric dynamical 
systems on the $p$-adic unit disk. The study is motivated by the problem of pseudorandom 
number generation  for computer simulation and cryptography. From this view the 
paper describes nonlinear congruential pseudorandom generators modulo $m$ which 
produce stricly periodic uniformly distributed sequences modulo $m$ with maximal 
possible period length (i.e., exactly $m$). Both the state change function
and the output function 
of these generators could be, e.g., meromorphic on 
$\Bbb Z_p$
functions (in particular, polynomials with rational, but not necessarily integer 
coefficients), or compositions of arithmetical operations (like addition, 
multiplication, exponentiation, raising to integer powers, including negative ones) 
with standard computer operations, such as bitwise logical operations 
(e.g.,$\XOR$, $\OR$, $\AND$, $\NEG$, etc.). The linear complexity of the produced 
sequences 
is also studied.    
\endabstract

\endtopmatter
\document

\head {1.} Introduction.\endhead
\par
A number of applications in computer simulation, numerical analysis (especially
Quasi Monte Carlo) and cryptography demand regular methods to generate
successively 
a uniformly distributed sequence. The corresponding literature
is so vast that we could not even mention here the most important monographs
in the area. We refer only \cite 2, where a reader could found a rather substantial
survey of relevant methods as well as a comprehensive bibliography. The
major part of 
these methods are certain recursive procedures, which 
may be viewed also as automata. The latter are commonly
referred as pseudorandom (or quasirandom) generators.

The typical one is the so-called linear congruential generator, which has
been developed more than half a century ago. It produces a sequence 
$\{x_n: n=0,1,2,\ldots\}$ over a set $\{0,1,\ldots,m-1\}$ (the latter is
commonly treated
as the residue class ring $\Bbb Z/m$ of the ring $\Bbb Z$ of rational integers 
modulo natural $m>1$), which 
is a first order recurrence sequence, defined by $x_{n+1}\equiv a+bx_n\pmod m$
with integer rationals $a,b$. The sequence is uniformly distributed
iff it is purely periodic with period length $m$. The latter condition
implies that {\it each} element of $\Bbb Z/m$ occurs at the period {\it exactly once}; and vice versa. The
necessary and suffcient conditions $a$ and $b$ must (for a given
$m$) satisfy  to provide the  maximal period length (i.e., $m$) of the produced 
sequence, are well known --see \cite {2, section 3.2.1.2, theorem
A}.    

The undoubtful advantage of linear congruential generators is the simplicity 
(especially for $m=2^k$) of their program implementations. One of the key reasons
of their disadvantages (e.g., lack of statistical quality of the produced sequences,
for certain
applications) is 
their linearity. For instance, as the state change function $f(x)=a+bx$ of the
generator is a polynomial of degree 1, the produced sequence
has linear complexity 2 over the ring $\Bbb Z/m$, i.e., it is a linear recurrence
sequence of order 2 over $\Bbb Z/m$ (defined by  $x_{n+2}\equiv
(1+b)x_{n+1}-bx_n\pmod m$). Hence, {\it for each $m$}  
the points $(\frac{x_{n+2}}{m},\frac{x_{n+1}}{m},\frac{x_n}{m})$ fall into
the parallel planes $c+X-(1+b)Y+bZ$ ($c\in\Bbb Z$), which intersect the
unit cube of Euclidean space. The well known result due to George Marsaglia \cite 7 states
that similar effect also holds in higher dimensions $>3$: all the
points fall into the relatively small number of parallel hyperplanes
(rather than fill this cube more or less
uniformly), 
and the reason is again that $\deg f=1$. 

During the past decades these considerations stimulated the developement
of various alternatives to linear congruential generators. The significant
part of these are {\it nonlinear} congruential generators with state change
function 
$f$ being either a polynomial over $\Bbb Z$ of degree $>1$,
(in particular, quadratic \cite 2, or of higher degree \cite{15}), or
some non-polynomial transformations, which gave rise to exponential generators
(with $f(x)=a^{g(x)}$),
or to so-called inversive generators, involving raising to negative powers 
(for the survey
of different generators
we again refer to \cite 2). 
Very often some authors seem to be more concerned with the 
linear complexity of the produced 
sequence, then with its uniform distribution, admitting non-maximal period length,
i.e., they admit state change functions $f$, for which
the sequence never reaches the period length $m$, and hence,
in a strict sence, is not uniformly distributed in $\Bbb Z/m$. In such
cases the authors have not only to estimate possible period lengths, but also
to choose the initial state (the {\it seed}) $x_0$ of the generator 
according to certain (sometimes,
sophisticated) procedures, which are to assure logging on the sufficiently
long period, rather then to choose the seed at random. 

Increasing the degree of a polynomial (as well as the use in the 
composition other arithmetical operations like exponentiation or taking
an inverse) is also has to be paid for by certain rise of complexity
of program implementation. The most disappointing here is that the statistical
quality and complexity of the program implementation of the generators
often occur to be in inverse dependence: the better the quality the slower the
performance; fast generators sometimes demonstrate lack of quality.

So it is still important to find new classes of functions 
$f\colon\Bbb Z/m\rightarrow\Bbb Z/m$, for which the corresponding generators,
\roster
\item firstly, achieve maximal possible (i.e., exactly $m$) period length  of
the recurrence sequence defined by the relation $x_{n+1}\equiv f(x_n)\pmod m$,
hence producing uniformly distributed sequence in $\Bbb Z/m$ (we refer
such transformations $f$ as {\it transitive modulo} $m$); 
\item secondly, guarantee
the suffuciently large linear complexity of the produced sequence over $\Bbb Z/m$,
i.e., absence of `short' (in some definite sence) linear dependencies
of the form $\sum_{i=0}^{r-1} c_i x_{n+i}\equiv 0\pmod m\ (n=0,1,2,\ldots)$
among the members of the sequence;
\item and, thirdly, basically are `easy-to-implement', namely,  are
`flexible', i.e., have certain (critical to the
performance) parameters, varying which it is possible to gain speed without loosing quality.
\endroster
The paper presents wide classes of transformations $f$ which to
some extend satisfy these conditions. 

At the first turn we obtain transitivity modulo $m$ conditions 
for functions, which could be implemented as compositions of arithmetical
operations (addition and multiplication of integers), as well as of standard computer
ones, like bitwise logical operations, shifts, masking, etc. These compositions
might involve as well exponentiation and taking a multiplicative inverse,
hence, raising to negative powers (see 4.9, 4.11,
2.5) and/or  $\OR, \XOR, \AND$, etc., see 2.5, 2.8.

In particular, we describe wide classes of transitive modulo $m$ functions
which could be expressed as integer-valued polynomials with rational (and
not
necessarily integer) coefficients (see 4.7), as well as by analytic functions (4.11, 4.9, 2.5)
or meromorphic (in particular, raional) functions (4.9, 4.11, 4.12). 
These conditions are easy-to-verify,
and with the use of them the various explicit formulae for transitive modulo $m$
transformations could be (and are) obtained -- see e.g. 2.3, 2.4, also
2.5--2.8 (as well 4.11, 4.12) together with 2.1, and other examples here
and there in the paper.  

To illustrate, we start with some of these examples: theorem 2.7 together with lemma 4.11 imply that
each transformation $f$ of the form
$$f(x)=1+x+2(g(x+1)-g(x))$$
is transitive modulo $m=2^k$ {\it for all} $k=1,2,\ldots$ and for {\it arbitrary} composition
$g$ of 
\roster
\item arithmetical operations --- an addition $(y,z)\mapsto
y+z$, a multiplication $(y,z)\mapsto yz$, an exponentiation $(y,z)\mapsto (1+2y)^z$
(in particular, taking an inverse $y\mapsto (1+2y)^{-1}$),
\item bitwise logical operations --- such as conjunction $(y,z)\mapsto y\AND z$, 
disjunction
$(y,z)\mapsto y\OR z$,  exclusive `or' $(y,z)\mapsto y\XOR z$, negation
$z\mapsto \NEG z$, etc.,
\item machine operations (which could be derived from the bitwise logical
ones) --- an $s$-step shift towards most significant bits $z\mapsto 2^sz$, masking
$z\mapsto z\AND M$, $M$ being a mask, `reduction modulo
$2^s$', i.e., a truncation of the
most significant bits $z\mapsto z\bmod 2^s=z\AND (2^s-1)$, and some others.
\endroster
We assume here that all the operands are non-negative integer rationals
which are represented as  base 2 expansions; so, for instance,  
$2=1\XOR 3 = 2\AND 7\equiv\NEG 13\pmod 8$, $3^{-1}\equiv
11\equiv -5\pmod {16}$, $3^{-\frac{1}{3}}\equiv
3^{11}\equiv 3^{-5}\equiv 11\pmod{16}$, etc.
Up to this agreement the functions $g$ and $f$ are correctly defined on
$\Bbb Z/m$, the efficiency of their program implementation depends only
on the number of `fast' and `slow' operations in the composition $g$ and hence
one may vary it in wide range to achieve the desired performance.

We emphasize, in the  example just mentioned the transitivity modulo $m=2^k$
{\it does not depend} neither on $k$ nor on actual form of the composition
$g$ --- both for 
$g(x)=x\XOR(2x+1)$ and for 
$$g(x)=\Biggl(1+2\frac{x\AND x^{2}+x^3\OR x^4}{3 + 4(5+6x^5)^{x^6\XOR x^7}}\Biggr)^{7+\frac{8x^8}{9+10x^9}}$$   
the sequence $\{x_n\}$ defined by the recurrence relation
$x_{n+1}\equiv 1+x_n+2(g(x_n+1)-g(x_n))\pmod {2^k}$ is uniformly distributed
in $\Bbb Z/2^k$  for each $k=1,2,3\ldots$. Actually, this sequence
is strictly periodic with period length $2^k$, and  {\it each} element
of 
$\{0,1,\ldots,2^k-1\}$ occurs at the period {\it exactly once}. 

Similar assertions also hold for arbitrary composite $m$: e.g.,
4.11 and 4.12 imply that the transformation 
$$f(x)=1+x+\pi(m)^2u(x)(1+\pi(m)v(x))^{w(x)}$$
with $\pi(m)$ being a product of all prime factors of $m$,
is transitive
modulo $m$ for arbitrary polynomials $u(x), v(x), w(x) \in \Bbb Z[x]$
over $\Bbb Z$.  A variety of results of such kind may be obtained in much more general situation
for integer-valued polynomials with rational (not necessarily integer) coefficients
by applying the techniques of section 4. 

Note that this example also demonstrates how by minor changes of the recurrence
relation one may achieve the transitivity of both inversive generator (for
which $f(x)=a+bx^{-1}$
or $f(x)=(a+bx)^{-1}$) and exponential
generator (with $f(x)=a^x$): for $w(x)=const=-1$ the introduced generator
is of inversive type, for $v(x)=const\ne 0$ it is of exponential type. 

As for linear dependencies 
$\sum_{i=0}^{r-1} c_i x_{n+i}\equiv 0\pmod m\ (n=0,1,2,\ldots)$ of fixed
length
$r$
in produced sequences $\{x_n\equiv f(x_{n-1})\pmod m: n=1,2,\ldots\}$,
one may say that from this view among all congruential generators linear ones
are rather exceptions
than the law. For instance, if $f\colon\Bbb Z\rightarrow\Bbb Z$ is
represented by a transitive
modulo some prime power $m=p^k$ ($k\ge 3$)  polynomial
of degree $\ge 2$ with integer rational coefficients, no such dependencies
with {\it $r$ and $c_i$ not depending on $k$} do
exist. Moreover, in this case the minimal order of linear recurrence sequence
over $\Bbb
Z/p^k$, which represents
the produced sequence, tends to infinity together with $k$  
(in fact, much more general result holds --- see 5.1--5.4 for exact statements).

The  paper also studies equiprobable modulo $m$ functions, i.e.,
mappings $F$ of the $s$th Cartesian power $\Bbb (Z/m)^{(s)}$ onto the
$t$th Cartesian
power $\Bbb (Z/m)^{(t)}$ of the ring $\Bbb Z/m$, $(s\ge
t)$, such that all preimages of all elements are of the same cardinality.
In particular, for
$s=t$ equiprobable modulo $m$ functions are bijections of the corresponding rings
and throughout the paper are referred as bijective modulo $m$ functions. A very particular case of the
equiprobable modulo $m$ functions studied here are so-called permutation
polynomials modulo $m$, the latter being polynomials over $\Bbb Z$ which induce bijections of the ring
$\Bbb Z/m$ onto itself. The results of the paper concerning equiprobability
modulo $m$ generalize known
\cite 8 results on permutation polynomials to much  wider classes
of functions. The study  was motivated by the observation that application
of equiprobable
modulo $m$ functions as output functions to uniformly distributed
in $\Bbb Z/M$
periodic sequences with period length $M$ leads to new uniformly distributed
in  $\Bbb Z/N$ (with $N|M$) sequences of the same period length $M$. In
other words, each element of $\Bbb Z/N$ occurs at the period of such 
sequence the same number of times (but not necessarily once). 
Hashing with equiprobable modulo $m$ functions
the sequences, generated by already introduced methods, 
seems to be useful to design secure stream
ciphers. Yet this will be an issue of the forthcoming paper and is out of the scope 
of the present one.

Note that proofs  of our basic assertions imply $p$-adic techniques. The
problems stated above are firstly restated in these terms. Actually the
paper studies ergodic with respect to the Haar measure (as well as preserving
this measure or equiprobable with respect to it) functions, which are defined on
(and
which take values in) the space $\Bbb Z_p$ of all $p$-adic integers, and which
are non-expanding functions, i.e., satisfy Lipschitz condition
with coefficient 1. From this view the results of the paper could be of
interest for non-Archimedean dynamical systems theory: a number of statements
could be easily interpreted as descriptions of ergodic dynamical
systems with discrete time and with $\Bbb Z_p$ as a phase space.

The paper continues the study started in \cite{11}: here we prove
some results announced in \cite{11, 12, 14, 17} 
and establish new ones. Moving towards exact statements, for reader's
convenience we recall some facts from the
$p$-adic analysis and the theory of uniform distribution of sequences, following
\cite 6, \cite 3 and \cite 2; we recall some necessary results, definitions
and notations from \cite{11} as well.
\par
Here and after let  $p$ be a prime number. Consider a canonic representation 
$z=z_0+z_1p+z_2p^2+\ldots$ of $p$-adic integer $z\ne 0$, where 
$z_j\in\{0,1,\ldots,p-1\} \  (j=0,1,2,\ldots)$; we denote $\ord_p\,z=\min\{j:z_j\ne 0\}$
the exponent of a maximal power of $p$ which is a factor of $z$. By 
definition,  $\|z\|_p=p^{-\ord_p\,z}$
is $p$-adic norm of $z$, $\|0\|_p=0$. The valuation $\| \ \|_p$ could be expanded
to the whole field  $\Bbb Q_p$ of $p$-adic numbers (which is a quotient
field of the ring $\Bbb Z_p$ of $p$-adic integers) in a standard way; so
this valuation induces on $\Bbb Q_p$ a distance $d_p(u,v)=\|u-v\|_p$, with $\Bbb Q_p$ being
a completion of the space $\Bbb Q$ of all rationals with respect to this
distance. Note that often they use another terminology, where 
a distance is called a {\it metric},
a $p$-adic norm is called a {\it $p$-adic  value}, whereas the term `$p$-adic
valuation' is reserved
for $\ord_p$. However, throughout the paper we mainly use the terminology of \cite
3, with the only exception, speaking of `$p$-adic norms' instead of 
`$p$-adic values'.  

The ring 
$\Bbb Z_p=\{u\in \Bbb Q_p : \|u\|_p\le    1\}$ is compact in the space 
$\Bbb Q_p$, being a closure of the set $\Bbb N_0=\{0,1,2,\ldots \}$.
Hence, $\Bbb Z_p$ is a separable compact metric space. The set of all cosets
with respect to all ideals of the ring $\Bbb Z_p$ forms a base of the corresponding
 topology. Each coset $a+p^k\Bbb Z_p$ ($a\in \Bbb Z_p$, $k=0,1,2,\ldots$) 
is an open (and simultaneously closed) ball of radius $p^{-k}$. 

There exists
a natural measure $\mu$ on $\Bbb Z_p$:
putting $\mu(a+p^k\Bbb Z_p)=p^{-k}$, we then expand $\mu$ to the correponding
$\sigma$-ring generated by all compact subsets of $\Bbb Z_p$ (these compact
subsets are exactly all closed subsets of  $\Bbb Z_p$). So we define uniquelly
a measure on $\Bbb Z_p$, which is non-negative  $\sigma$-additive
regular normalized Borel and Haar measure in this case. Thus, $\mu$ is a natural probability measure
on $\Bbb Z_p$. The probability measure on $n$-dimensional space 
$\Bbb Z_p^{(n)}$ could be defined in a similar way as a corresponding
normalized Haar measure.

Now let $f\colon {\Bbb Z}_{p}\rightarrow {\Bbb Z}_{p}$ be a function, which
preserves all congruences of the ring
$\Bbb Z_p$, i.e.,
$a\theta b$ implies $f(a)\theta f(b)$ for each congruence  $\theta$ and
all $a,b\in\Bbb Z_p$. As each nontrivial congruence of the ring 
$\Bbb Z_p$ is an equivalence modulo certain ideal
$p^k\Bbb Z_p$ (we denote this congruence  via
$\cdot\equiv\cdot\pmod{p^k}$), it can be easily shown that the function
$f$ preserves all the congruences of the ring $\Bbb Z_p$ iff it satifies
Lipschitz condition with coefficient $1$:\quad $\| f(x)-f(y)\| _{p}\le
\| x-y\| _{p}$. Function preserving all congruences of a universal
algebra is called {\it compatible}; we will use this term instead of the
term `conservative' of \cite{11}, since the latter in numerous papers
on algebraic systems has another meaning, see  \cite {8, p. 45}.
\par
The class of all compatible functions on $\Bbb Z_p$ is rather wide: it
contains all functions represented by polynomials with rational integer
or $p$-adic integer coefficients, integer-valued analytic on $\Bbb Z_p$ functions,
as well as integer-valued and meromorphic (in particular, rational) on 
$\Bbb Z_p$ functions with denominators equivalent to 0 modulo $p$ at no
point of $\Bbb Z_p$. Some other examples will be introduced further in
the paper.

Recall that a function, which is defined on some field $F$, and which takes
values there, is called
{\it integer-valued} iff all its values are integers of
$F$ providing arguments take integer values in $F$. 
Further we study integer-valued functions on the field $\Bbb Q_p$; hence
they map $\Bbb Z_p$ into $\Bbb Z_p$. In particular, we study integer-valued
functions on $\Bbb Q$. A polynomial over a field $F$ is called integer-valued
iff it induces an integer-valued function on $F$. Note that integer-valued
function $f$ on $F$ defines on the ring $Z$ of all integers of
 $F$ a function $f|_Z:Z\rightarrow Z$, which is not necessarily compatible
 on $Z$, i.e., does not necessarily preserve all congruences of  $Z$; yet,
 if $f|_Z$ is compatible as a function on
$Z$, then (in cases which do not lead to misunderstanding) we also call
$f$ compatible. Moreover, if a compatible integer-valued function $f$ could be 
defined by a polynomial over $F$, we
call compatible the corresponding polynomial too.      
\par
Note that the notion of compatible integer-valued function could be naturally
expanded to the multivariate case --- a valuation (and hence, a
distance)
on the space $\Bbb Z_p$ induces a (pseudo)-valuation (hence, a distance) on the  
$n$-dimensional space $\Bbb Z_p^{(n)}$ in a standard manner:  for $ \bold u=(u_1,\ldots,
u_n)\in\Bbb Z_p^{(n)}$ we assume $\|\bold u\|_p=\max\{\|u_i\|_p:
i=1,2,\ldots,n\}$. So, the function
$$F=(f_1,\ldots,f_m)\colon {\Bbb Z}_{p}^{(n)}\rightarrow {\Bbb Z}_{p}^{(m)}$$
is compatible iff it satisfies Lipschitz condition with coefficient 1. In
particular, all compatible on $\Bbb Z_p$ functions are continuous as functions
of $p$-adic  variables.

This obvious conclusion is important for applications. Each machine word,
i.e., a word of some finite length in the alphabet $\{0,1\}$, could be
treated as a 2-base expansion of a non-negative integer rational. Then
all the above mentioned bitwise logical operations 
and machine operations could be naturally continued to the set $\Bbb Z_2$
of all 2-adic integers in their canonic representations. Moreover, the
above mentioned arithmetical operations could be continued to
$\Bbb Z_2$ either. It could be
easily demonstrated that all these operations (to be more precise, their uniquelly
defined continuations
to $\Bbb Z_2$) and all their compositions are compatible (hence, continuous)
integer-valued functions on $\Bbb Z_2$: for exponentiation
$(y,z)\mapsto (1+2y)^z$, and, in particular, for the inversion
$y\mapsto (1+2y)^{-1}$ 
see 4.11, for the rest the assertion follows immediately from the corresponding
definitions. We note here that an $m$-step shift towards less
significant bits (i.e., the operation $\lfloor\frac{\cdot}{2^m}\rfloor$
of `integer division', a division succeeded by a truncation of the fractional part
of the quotient) is {\it not} compatible, yet continuous, integer-valued
function on $\Bbb Z_2$ (hence the results of the paper remain valid for
compositions including the latter operation either, providing the whole
composition is compatible). 

These considerations give an opportunity to apply, while studing compositions
of the above mentioned operations,
certain methods of non-Archimedean ($p$-adic) analysis. Certainly, these techniques
could be applied only to problems which are stated in appropriate terms
(measures, distances, limits, derivatives, etc.).
\par
It turnes out that some properties of functions, which traditionally
have been treated as discrete mathematics issues, could be restated in these
terms. We have already introduced one of such properties, namely, compatibility. 
It worth a brief notice in this connection that so-called `determinate functions 
on superwords' of automata theory (which
are functions defined on infinite sequences of $\{0,1\}$) after natural
identification of superwords with elements of
$\Bbb Z_2$ could be considered as compatible functions on  $\Bbb Z_2$.
\par
There exist other properties which could be restated in such manner. We
consider a property of a compatible function 
$F=(f_1,\ldots ,f_m)\colon {{\Bbb Z}_{p}}^{(n)}\rightarrow {{\Bbb Z}_{p}}^{(m)}$
to be {\it equiprobable modulo} $p^k$. The latter by definition means
that the function $F$ induces on the $n$th Cartesian power  $(\Bbb Z/p^k)^{(n)}$ of $\Bbb Z/p^k$ 
an equiprobable function
$\bar F=(\bar f_1,\ldots ,\bar f_m)\colon ({\Bbb Z}/{p^k})^{(n)}\rightarrow ({\Bbb Z}/{p^k})^{(m)}$, 
i.e., each point of ${(\Bbb Z/p^k)}^{(m)}$ has 
the same number of $F$-preimages in ${(\Bbb Z/p^k)}^{(n)}$. In particular, for $m=n$  equiprobable
 modulo $p^k$ functions are exactly {\it bijective modulo} $p^k$ functions. We
 consider also an important (especially for pseudorandom generation) property
 of a bijective modulo $p^k$ function  $F$ to be
{\it transitivite modulo} $p^k$, which means that $F$ induces on ${(\Bbb Z/p^k)}^{(n)}$ 
a single cycle permutation.
Note that while defining notions of equiprobability, bijectivity or transitivity
of a function  $F$ modulo $p^k$, we have assumed the compatibility of $F$. 

A
value of induced function $\bar f_i(x)$ in the ring $\Bbb Z/p^k$
is, by definition, $f_i(x)\bmod p^k $, the least non-negative residue   modulo
$p^k$ of $f_i(x)$, i.e., 
$f_i(x)\bmod p^k =\alpha
\in \{0,1,\ldots ,p^k-1\}$, with $\| f_i(x)-\alpha \|_p \le p^{-k}$.
In view of compatibility of the function $f_i$, the value of the function
 $\bar f_i(x)$
does not depend on choice of the representative $x$ in a coset of
the ring $\Bbb Z_p^{(n)}$
with respect to the ideal
$(p^k \Bbb Z_p)^{(n)}$; hence, the function $F$ correctly defines on 
$(\Bbb Z/p^k)^{(n)}$ a function $F\bmod{p^k}=(f_1(x)\bmod p^k,\ldots,f_m(x)\bmod p^k)$, 
which takes values in $(\Bbb Z/p^k)^{(m)}$. Throughout the paper the latter function
is  denoted via $F\bmod{p^k}$,
or via $\bar F$, when it does not lead to misunderstanding.
\par
Now recall some definitions of the theory of measurable functions (cf.
\cite{1}).
Let $S$ and $T$ be spaces with nonnegative normalized measures 
$\mu $ and $\tau $, respectively, and let $f\colon S\rightarrow T$ be a
measurable function, i.e., each full $f$-preimage $f^{-1}(U)$ of $U\subseteq
T$ is  $\mu $-measurable for each $\tau $-measurable
$U$. 

We say that the function $f$ is  $(\mu,\tau)$-{\it proportional}, iff 
for each pair of $\tau$-measurable subsets $U,V\subseteq T$ 
the equality $\tau (U)=\tau (V)$ implies the equality
$\mu (f^{-1}(U))=\mu (f^{-1}(V))$. 
In case both $\mu ,\tau $ are probability measures (e.g., are properly normalized
Haar measures), then 
$f$ is called $(\mu,\tau)$-{\it equiprobable} (or {\it equiprobable with
respect to} $\mu$ and $\tau$) iff  
$\mu (f^{-1}(U))=\tau (U)$ for each $\tau$-measurable $U\subseteq T$. For
$S=T$ and $\mu =\tau $ we say that $f$ {\it preserves measure} $\mu$, iff 
$\mu (f^{-1}(U))=\mu (U)$ holds for each $\mu$-measurable 
$U$.  Finally, if $f$ preserves measure $\mu$, and for $\mu$-measurable
subset 
$U$ the equality 
$f^{-1}(U)=U$ implies that  either $\mu (U)=0$, or $\mu (U)=1$,
we say that $f$ is $\mu$-{\it ergodic} (or {\it ergodic with respect to} $\mu$).

Note that in metric theory instead of terms `measure-preserving function'
or `equiprobable function' they  often use terms `metric endomorhism' 
and `metric homomorphism', and in dynamical systems theory they sometimes speak about `metric
transitivity' instead of ergodicity.
Since throughout the paper we deal with the only measure, the properly
normalized Haar measure, we omit mentioning this measure, so preserving
the Haar measure, equiprobable (accordingly, ergodic) with respect to the
Haar measure functions are referred as  {\it measure-preserving}, {\it equiprobable} (or,
accordingly, {\it ergodic}).

The following theorem holds:
\proclaim
{1.1 Theorem}{
A compatible function
$F\colon {{\Bbb Z}_{p}}^{(n)}\rightarrow {{\Bbb Z}_{p}}^{(m)}$
is equiprobable {\rom(}respectively, measure-preserving or 
ergodic{\rom)} iff
it is equiprbable  {\rom(}respectively, is bijective or transitive{\rom)} modulo
$p^k$
for all $k=1,2,\ldots$. A compatible and measure-preserving function
$F$ is bijective {\rom(}consequently, is a metric automorphism{\rom)};
moreover, it is an isometry of the space $\Bbb Z_p^{(n)}$.}
\endproclaim

Note that further throughout the paper while proving ergodicity (equiprobability)
of a compatible function with respect to the Haar measure we actually prove
its transitivity (equiprobability) each modulo $p^k$,
$k=1,2,\ldots$, i.e., directly establish the properties we are interested
in view of the problems mentioned above. That is why we omit the proof
of this theorem 1.1: it is not related directly to the aims of this paper.
Nevertheless throughout the paper we use the relevant terminology (e.g.,
we commonly speak of `ergodicity' instead of  `transitivity modulo
$p^k$ for all $k=1,2,\ldots$', etc.)
  
In connection with theorem 1.1 it is worth noticing, however, that the results
of the paper related to description of measure-preserving or ergodic
functions may be treated as description of non-Archimedean (i.e.,
ultrametric) dynamical systems 
$(\Bbb Z_p^{(n)}, F)$ with phase space $\Bbb Z_p^{(n)}$, discrete time,
and with nonexpanding $F$ (i.e. for each pair of points $\bold a,\bold
b$  a distance between their $F$-images
$F(\bold a)$ and $F(\bold b)$ does not exceed a
distance between these points).
In this sence theorem 2.2, for instance, might be condidered as a complete
description (in terms of explicit formulae) of ergodic dynamical systems of the above mentioned kind
when
$p=2$ and $n=1$; together with theorem 3.11 it gives full description
of twice integer-valued (i.e., having everywhere integer-valued derivative)
ergodic dynamical systems. These themes, however, are not covered by this
paper and will be considered in forthcoming one.

Returning to the leading theme of the paper we note that for a wide class
of compatible functions, which are in some (properly defined in section 3)
sence generalizations of uniformly differentiable on
$\Bbb Z_p$ functions,  the bijectivity modulo 
$p^k$ of a function for a {\it certain} $k$ is equivalent to the 
property of being measure-preserving; the latter is equivalent to its bijectivity
modulo 
$p^k$ {\it for all} 
$k=1,2,3,\ldots$ . The property of being transitive modulo $p^k$ 
for a {\it certain} $k$ turned out to be equivalent to the ergodicity of
a function; the latter implies that the function is transitive modulo
$p^k$ {\it for all} 
$k=1,2,3,\ldots$. Finally, the equiprobability of a function modulo $p^k$ 
for a {\it certain} $k$ implies its equiprobability with respect to the
Haar measure; the latter property is equivalent to equiprobability modulo
$p^k$ {\it for all} $k=1,2,3\ldots$. The results of this kind are proved
in section 3.

These results demonstrate the same remarkable effect originally enlighted
by Hensel lemma: the Hensel lift, that is, a situation when a behavior
of a function modulo
$p^{k_0}$ for a certain $k_0$ controls its behavior modulo $p^k$ for all
$k=k_0+1,k_0+2,\ldots$ and on the whole space $\Bbb Z_p$. This effect have
been already
observed while studying transitivity of some transformations. For instance,
the necessary and sufficient conditions for the polynomial
$f(x)=a+bx$ with integer rational $a,b$ (see e.g., \cite {2; 3.2.1.2, theorem A})
could be restated as follows: a polynomial $a+bx$ is transitive modulo $p^k$ 
for some (that is, for all) $k\ge 2$ iff it is transitive  modulo $p$
for odd $p$ or, respectively, modulo $p^2$ for 
$p=2$. The general criterion for the transitivity modulo $p^k$ of the polynomial 
$f$ of arbitrary degree over integer rationals \cite {15} demonstrates this effect
either: for $p\ne 2,3$ a polynomial $f$ is transitive modulo $p^k,\  k\ge3$,
iff it is transitive modulo $p^2$; respectively, for $p=2$ or $p=3$ ---
iff it is transitive modulo $p^3$. Note by the way that the latter result
holds for a much wider class of functions, even not necessarily analytic (see
4.9--4.10). 

The results of section 3 show that Hensel lift of such properties as bijectivity
or transitivity modulo $p^k$ is basically due to the specific character of $p$-adic distance
and holds for various rather wide classes of functions. The values
of  $k_0$ from which the lift starts are estimated in section 4.
\par
The results of this kind are useful if for a given $f$ one has to establish
whether it shares some property (e.g., transitivity or bijectivity)
modulo  $p^k$ for a definite rather large $k$, for which direct verification
is not accessible. However, if one needs to construct out of prescribed operations
a certain function, which is to be transitive or bijective modulo
 $p^k$, then explicit formulae are more convenient. Such formulae for
 bijective modulo $2^k$ polynomials over $\Bbb Z$ were obtained in
\cite{13}, for transitive modulo $2^k$ polynomials over $\Bbb Z$
 --- in \cite{15}. Explicit formulae for ergodic or measure-preserving
 compatible functions (in particular, for compatible integer-valued polynomials
 over $\Bbb Q$), which are defined on (and take values in) $\Bbb Z_2$ were
 obtained in \cite{11}. 
 The current paper presents explicit formulae for compatible ergodic (or
 measure-preserving) functions on $\Bbb Z_p$ for odd $p$ --- see the next section.
\par
\head {2.} Explicit formulae \endhead 

Recall (see \cite 3) that each function $f\colon{\Bbb N}_{0}\rightarrow {\Bbb Z}_{p}$
(or, respectively, $f\colon{\Bbb N}_{0}\rightarrow {\Bbb Z}$)
admits one and only one representation in the form of   
so-called {\it interpolation series}
$$
f(x)=\sum^{\infty }_{i=0}a_{i}{{x}\choose{i}},
\eqnum{\diamondsuit}$$
\noindent where $
\left(\matrix x\\ i\endmatrix\right)
=\dfrac{x(x-1)\cdots  (x-i+1)}{i!}$ 
for $i=1,2,\ldots$, and
$\left(\matrix x\\ 0\endmatrix\right)=1$;  $a_{i}\in {\Bbb Z}_{p}$ 
(respectively, $a_{i}\in {\Bbb Z}$), $i=0,1,2,\ldots $ .
\par 
If $f$ is uniformly continuous on ${\Bbb N}_{0}$ with respect to
$p$-adic distance, it can be uniquely continued to the uniformly continuous
function on ${\Bbb Z}_{p}$. Hence the interpolation series for 
$f$ converges uniformly on 
${\Bbb Z}_{p}$. The following is true: the series
$
f(x)=\sum^{\infty }_{i=0}a_{i}
{x\choose i}
, \quad$ 
($a_{i}\in {\Bbb Q}_{p}$, $i=0,1,2,\ldots \ $)
converges uniformly on ${\Bbb Z}_{p}$ iff 
$ \lim\limits^p_{i \to \infty }a_{i}=0,$
where $\lim\limits^p$ is a limit with respect to $p$-adic distance; 
hence the uniformly convergent series defines a uniformly continuos function
on ${\Bbb Z}_{p}$. The latter function is integer-valued iff
 $a_{i}\in {\Bbb Z}_{p}$ for all $i=0,1,2,\ldots \ $.
\par
Further throughout this section we assume that the function $f\colon{\Bbb Z}_{p}\rightarrow {\Bbb Z}_{p}$ 
is uniformly continuous on ${\Bbb Z}_{p}$, and that it is represented by
series $(\diamondsuit)$.
The following three criteria hold (see \cite{11}):
\par
\proclaim
{ 2.1 Theorem} { {\rom {(See 4.3 of \cite{11}; cf. \cite 5)}} A function
$f\colon{\Bbb Z}_{p}\rightarrow {\Bbb Z}_{p}$ 
is compatible iff
$$
a_{i}\equiv 0\pmod {p^{\lfloor\log_{p}i\rfloor}}
$$
\noindent для for all $i=p,\ p+1,\ p+2,\ldots  $ . {\rom { (Here and after
for a real $\alpha $ we denote
$\lfloor\alpha \rfloor$ an integral part of $\alpha $, i.e., the nearest
to $\alpha $ integer rational not exceeding $\alpha $.)}}
\endproclaim
\par
\proclaim
{2.2 Theorem} {\rom{(See 4.5 of \cite{11})}} {A function $f\colon{\Bbb Z}_{2}\rightarrow {\Bbb Z}_{2}$ 
is compatible and measure-preserving iff it could be represented as 
$$f(x)=c_0+x+\sum^{\infty }_{i=1}c_{i}\,2^{\lfloor \log_2 i \rfloor +1}{x \choose i},$$
\noindent where $c_0, c_1, c_2 \ldots \in {\Bbb Z}_2$.}
\endproclaim
\par
\proclaim
{ 2.3 Theorem} {\rom{(See 4.7 of \cite{11})}} { A function $f\colon{\Bbb Z}_{2}\rightarrow {\Bbb Z}_{2}$
is compatible and ergodic iff it could be represented as 
$$
f(x)=1+x+\sum^{\infty }_{i=0}c_{i}\,2^{\left\lfloor \log_{2}(i+1)\right\rfloor+1}{ {x}\choose {i}},
$$
\noindent where $c_0, c_1, c_2 \ldots \in {\Bbb Z}_2$.}
\endproclaim
\par
For an arbitrary prime $p$ the necessity of condtions of theorems  2.2 and
2.3 does not hold, yet the sufficientness remains true. Namely, in this section
we prove the following: 
\par
\proclaim {2.4 Theorem} {Let $p$ be an odd prime. A function $f\colon{\Bbb Z}_{p}\rightarrow {\Bbb Z}_{p}$,
which is represented in the form $(\diamondsuit)$,  is  
compatible and measure-preserving if the following
congruences hold simultaneously:
$$\displaylines{
a_{1}\not\equiv 0\pmod{p};
\cr
a_{i}\equiv 0\pmod{p^{\left\lfloor{\log_{p}i}\right\rfloor+1}},\ \  (i=2,3,\ldots  ).
\cr}
$$
\par 
The function $f$ if compatible and ergodic if the following congruences
hold simultaneously:
$$\displaylines{
a_{0}\not\equiv 0\pmod{p};
\cr
a_{1}\equiv 1\pmod p;
\cr
a_{i}\equiv 0\pmod{p^{\left\lfloor{\log_{p}(i+1)}\right\rfloor+1}},\ \  (i=2,3,\ldots  ).
\cr}$$}
\endproclaim
For the proof of the theorem we will need two additional results which are of interest
by their own.
\par
\proclaim {2.5 Lemma} {Let $p$ be an arbitrary prime, let $v\colon{\Bbb Z}_{p}\rightarrow {\Bbb Z}_{p}$ 
be a compatible function, and let $c,d$ be $p$-adic integers,
with $c\not\equiv 0\pmod p$. Then the function $g(x)=d+cx+pv(x)$ 
preserves measure, and the function 
$h(x)=c+x+p\Delta v(x)$ is ergodic. {\rom {(Here and after $\Delta$ is
a difference operator:
$\Delta v(x)= v(x+1)-v(x)$. Note that both $g$ and
$h$ are obviously compatible since they are compositions of compatible
functions.)}}}
\endproclaim
\par
\demo {Proof of the lemma 2.5} Firstly by induction on  $l$ we show that 
$g$ is bijective modulo $p^{l}$ for all  $l=1,2,3,\ldots  $ . The assumption
is obviously true for $l=1$.
\par
Assume it is true for $l=1,2,\ldots  ,k-1$. Prove that it holds for
$l=k$ either.
Let $g(a)\equiv g(b)\pmod{p^{k}}$ for some $p$-adic integers $a,b$. 
Then $a\equiv b\pmod{p^{k-1}}$ by the induction hypothesis. Hence  
$pv(a)\equiv pv(b)\pmod{p^{k}}$ since $v$ is compatible. 
Further, the congruence $g(a)\equiv g(b)\pmod {p^{k}}$ implies that 
$ca+pv(a)\equiv cb+pv(b)\pmod{p^{k}}$, and consequently, $ca\equiv cb\pmod{p^{k}}$. 
Since $c\not\equiv 0\pmod p$,
the latter congruence implies that $a\equiv b\pmod{p^{k}}$, proving the
first assertion of the lemma.
\par
To prove the rest part of the statement we note firstly that the assertion
just proven implies that $h$ preserves measure.
To prove the transitivity of $h$ modulo $p^{k}$ for all $k=1,2,3,\ldots  $ 
we apply induction on $k$ once again.
\par
It is obvious that $h$ is transitive modulo $p$. Assume that $h$ is transitive
modulo $p^{k-1}$.
Then, since $h$ induces a permutation on ${\Bbb Z}/p^{k}$
and since it is a compatible function, we conclude that the length of each
cycle of this permutation must be a multiple of $p^{k-1}$. So to prove
this permutation is single cycle it is sufficient to prove that the function 
$$h^{p^{k-1}}(x)=\underbrace {h(h\ldots  (h}_{p^{k-1}\;\text{ раз}}(x))\ldots)$$ 
induces a single cycle permuation on the ideal $(p^{k-1})$, generated by
the element $p^{k-1}$
of the ring ${\Bbb Z}/p^{k}$.  In other words, it is sufficient to demonstrate
that the function ${\frac {1} {p^{k-1}}}h^{p^{k-1}}(p^{k-1}x)$ 
is transitive modulo $p$.
\par
Applying obvious direct calculations, we successively obtain that 
$$h^{1}(x)=c+x+pv(x+1)-pv(x),$$
$$\ldots \qquad \ldots \qquad \ldots  $$
$$\displaylines{h^{j}(x)=h(h^{j-1}(x))=cj+h^{j-1}(x)+pv(h^{j-1}(x)+1)-pv(h^{j-1}(x))\hfill\cr
\hfill{}=cj+x+p\sum^{j-1}_{i=0}v(h^{i}(x)+1)-p\sum^{j-1}_{i=0}v(h^{i}(x)),\cr}$$
\par
\noindent and henceforth. We recall that $h^{0}(x)=x$ by definition. So,
$$
h^{p^{k-1}}(x)=cp^{k-1}+x+p\sum^{p^{k-1}-1}_{i=0}v(h^{i}(x)+1)-p\sum^{p^{k-1}-1}_{i=0}v(h^{i}(x)).
\eqnum{1}$$
\par
Since $h$ is transitive modulo $p^{k-1}$ and compatible, we get now that 
$$
\sum^{p^{k-1}-1}_{i=0}v(h^{i}(x)+1)\equiv \sum^{p^{k-1}-1}_{i=0}v(h^{i}(x))\equiv \sum^{p^{k-1}-1}_{z=0}v(z)\pmod{p^{k-1}},
$$
\noindent and (1) implies then 
$h^{p^{k-1}}(x)\equiv cp^{k-1}+x\pmod{p^{k}}$. But  $c\not\equiv 0\pmod p$, 
so we conclude that the function $cp^{k-1}+x$ induces on the ideal $(p^{k-1})$ a
single cycle permutation, thus proving the lemma.\qed
\enddemo
\par
\proclaim {2.6 Corollary} Under the assumptions of lemma 2.5, let $p$ be
an
odd prime, and let
$r\equiv 1(\bmod p)$. Then the function 
$c+rx+p\Delta v(x)$ is compatible and ergodic.
\endproclaim
\par
\demo
{Proof of the collorary 2.6} We have that $r=1+ps$ for a suitable $s\in {\Bbb Z}_{p}$.
Now, since $p$ is odd, the function $s{ {x}\choose {2}}$ is compatible; consequently,
the function $v_{1}(x)=s{ {x}\choose {2}}+v(x)$ is compatible either. Yet
$\Delta v_{1}(x)=sx+\Delta v(x)$, and it is sufficient now to apply 
lemma 2.5 to finish the proof of the corollary.\qed
\enddemo
\demo
{Proof of the theorem 2.4} Note that according to 2.1 a compatible function
$v(x)$ could be represented as 
$$
v(x)=a+ \sum^{\infty }_{i=1}b_{i}p^{\left\lfloor{\log_{p}i}\right\rfloor}{ {x}\choose {i}},
$$
\par
\noindent where $a,b_{1},b_{2},\ldots  \in {\Bbb Z}_{p}$. As 
$\left\lfloor{\log_{p}i}\right\rfloor=\left\lfloor{\log_{p}(i+1)}\right\rfloor$ 
for all $i=1,2,\ldots  $ with the exception of  $i=p^{t}-1$, $(t=1,2,3,\ldots)$, and
as
$$
\Delta v(x)=\sum^{\infty }_{i=1}b_{i}p^{\left\lfloor{\log_{p}i}\right\rfloor}{ {x}\choose {i-1}},
\eqno{(1)}
$$
\noindent we finish the proof of the theorem, applying  2.5 and 2.6.\qed
\enddemo
\par

For $p=2$ the results just proven imply one more useful
criterion of ergodicity of a function (or being measure-preserving).
\proclaim
{2.7 Theorem}{A function $f\colon\Bbb Z_2\rightarrow\Bbb Z_2$ is compatible
and preserves measure  {\rom (}respectively, is compatible and ergodic{\rom )} iff it can be represented in the form
$f(x)=c+x+2v(x)$ {\rom (}respectively, in the form
$f(x)=1+x+2\Delta v(x)${\rom )}, where $c\in\Bbb Z_2$ and $v(x)$ is a
compatible function.}
\endproclaim
\demo
{Proof} Follows easily from 2.1--2.3 and 2.5 in combination with  (1)
of the proof of the theorem  2.4. \qed
\enddemo

Both 2.5--2.6 and theorem 2.7 could be applied to consruct 
measure-preserving or ergodic functions as compositios of given compatible
functions. For instance, putting
$v(x)=(x^2)\XOR(x+32\AND x)$ (this function is compatible as a composition
of compatible functions) we conclude that the function
$$
7+x+2((x^2+2x+1)\XOR(x+1+32\AND(x+1)))-2(x^2\XOR(x+32\AND x))
$$
is ergodic. This conclusion is not very easy to verify by direct application
of theorems  2.2 or 2.3.
 
By the way, for $p=2$ the statement of theorem 2.7 could be slightly modified
to make it a little bit more convenient for the construction of ergodic
functions out of addition and bitwise logical operations (like 
bitwise exclusive `or', $\XOR$, bitwise `and', $\AND$, or bitwise negation $\NEG$).
Namely, it could be easily seen that in the ring $\Bbb Z_2$ there holds
an 
identity 
$z+{{\NEG}}(z) = -1$. 
Hence, $\Delta v(x)=v(x+1)-v(x)=v(x+1)+{{\NEG}}(v(x))+1$, and we obtain
the following
\proclaim
{2.8 Proposition} A function $f\colon\Bbb Z_2\rightarrow\Bbb Z_2$ is compatible
and ergodic iff it can be represented in one {\rom(}hence, all{\rom)} of the following forms
$f(x)=1+x+2(v(x+1)+{{\NEG}}v(x))=2+x+2v(x+1)+{{\NEG}}(2v(x))=3+x+2v(x+1)+2{{\NEG}}v(x)$,
where $v\colon\Bbb Z_2\rightarrow\Bbb Z_2$ is an arbitrary compatible function.
\endproclaim
\par
Since multiplication by $2$ is just a 1-digit shift of 2-base expansion of a number
towards senior bits, the proposition
2.8 could be applied to construct pseudorandom number generators out of
the `fast' computer commands, like addition, bitwise logical oherations
and shifts towards senior bits, by implementing the function $v$ as a composition
of them.
\par
It worth noticing also that all the functions described in 2.4 -- 2.8
are `affine modulo $p$', i.e., induce on $\Bbb Z/p$ a transformation of
the form $x\mapsto a+bx$.

\head{3.} Hensel lift.\endhead
\par
This section studies conditions when a function of an important class of
uniformly differentiable modulo $p^k$ functions (which are properly defined below), 
is
equiprobable, measure-preserving or ergodic. As a rule, the results of the section
demonstrate the effect of Hensel lift, mentioned in the introduction:
speaking loosely, if a function $F$ has some property modulo $p^{k_0}$
then it has this property modulo $p^n$ for all
$n\ge k_0$. Besides, it worth noticing here that the results of this section, 
contrasting those
of the previous one, provide some tools to
construct measure-preserving or ergodic functions which are not necessarily
affine modulo 
$p$. In fact, a certain techniques
based on the ideas of this section 
could be developed;  these techniques enables one `to lift' an arbitrary transitive
transformation of the ring $\Bbb Z/p^{k_0}$ to the function on $\Bbb Z_p$, 
which is transitive modulo $p^k$ for all $k=k_0, k_0+1,
k_0+2,\ldots$. This is the main reason we introduce a notion of asymptotically
compatible function below. 
However, the techniques themselves are not discussed here being  
left
to the forthcoming paper.
\par
Firstly, recall some generalizations of our basic notions (see 5.1 of \cite{11}).
\definition 
{ 3.1 Definition}   Let $F=(f_{1},\ldots  ,f_{m})\colon{\Bbb Z}^{(n)}_{p}\rightarrow {\Bbb Z}^{(m)}_{p}$ be
a function, not necessarily compatible. The function $F$ is called ({\it 
asymptotically}) {\it equiprobable},
iff for all $k=1,2,\ldots$
(respectively, for all sufficiently large $k\in\Bbb N$) it is {\it equiprobable
modulo $p^{k}$}, that is,
the restriction $F\bmod{p^{k}}=(f_{1}\bmod {p^k},\ldots  ,f_{m}\bmod p^k)$ of
the function $F$ to the set 
$\{0,1,\ldots  ,p^{k}-1\}^{(n)}$
is an equiprobable function. (Note that in cases which do not lead to misunderstanding
we identify the set  
$\{0,1,\ldots  ,p^{k}-1\}^{(n)}$
with the set of all elements of the ring $(\Bbb Z/p^k)^{(n)}$). By the
analogy, we say that $F$ is {\it asymptotically measure-preserving} (respectively,
that $F$ is {\it asymptotically
ergodic}), iff $F\bmod{p^{k}}$ is a bijective (respectively, transitive) 
transformation of the ring $(\Bbb Z/p^k)^{(n)}$
for all sufficiently large $k$. Lastly, we say that   $F$ is {\it asymptotically
compatible} iff there exists positive integer rational
 $N$ such that for all $\bold a,\bold b\in {\Bbb Z}^{(n)}_{p}$
and all $k\ge N$ a congruence $\bold a\equiv \bold b\pmod{p^{k}}$ implies a congruence
 $F(\bold a)\equiv F(\bold b)\pmod{p^{k}}$. 
\enddefinition
\par 
By definition, for $ \bold a=(a_{1},\ldots  ,a_{n})$ and
 $\bold b=(b_{1},\ldots  ,b_{n})$ of
${\Bbb Q}^{(n)}_{p}$ 
the congruence $\bold  a\equiv  \bold b\pmod{p^{s}}$ means that
$\|a_{i}-b_{i}\| _{p}\le p^{-s}$ (or, the same, that $a_{i}=b_{i}+c_{i}p^{s}$ 
for suitable $c_{i}\in {\Bbb Z}_{p}$, $i=1,2,\ldots  ,s$); that is $\|\bold
a-\bold b\|_p\le p^{-s}$.
In other words, a function is asymptotically compatible iff for some
$N\in {\Bbb N}_{0}$ it satisfy Lipschitz condition with coefficient 1 for
each pair of points which are at least as close one to another as $p^{-N}$. 
Since $\Bbb Z_p^{(n)}$ is compact, $F$ is asymptotically
compatible iff it satisfy Lipschitz condition with coefficient 1 locally.
\par 
Now for reader's convenience we recall some 
facts of \cite{11}. A function $F=(f_{1},\ldots  ,f_{m})\colon{\Bbb Z}^{(n)}_{p}\rightarrow {\Bbb Z}^{(m)}_{p}$
is called {\it differentiable modulo $p^k$} at the point 
$ \bold u=(u_{1},\ldots  ,u_{n})\in {\Bbb Z}^{(n)}_{p}$, iff there exist a positive
integer rational
$N$ and $n\times m$ matrix $F^{\prime}_{k}(\bold u)$ over ${\Bbb Q}_{p}$
(called {\it a Jacobi matrix modulo} $p^{k}$ of the function $F$ at the
point
$\bold u$) such that for each positive integer rational 
$K\ge N$ and each $ \bold h=(h_{1},\ldots  ,h_{n})\in {\Bbb Z}^{(n)}_{p}$  
the inequality 
$\|\bold h\| _{p}\le     p^{-K}$ implies a congruence  
$$F( \bold u+\bold h)\equiv F(\bold u)+ \bold hF^{\prime}_{p}(\bold u)\pmod{p^{k+K}}.\eqnum{\heartsuit } $$
 In case $m=1$ a
Jacobi matrix modulo $p^k$ is called a {\it differential modulo $p^k$}. In
case $m=n$ a determinant of Jacobi matrix modulo $p^k$ is called a {\it Jacobian
modulo $p^k$}. The elements of Jacobi matrix modulo $p^k$
are called {\it partial derivatives modulo} $p^k$ of the function $F$ at
the point $\bold u$. 
A partial derivative (respectively, a differential) modulo $p^k$ are
sometimes  denoted as 
$\frac{\partial_k f_i (\bold u)}{\partial_k x_j}$ (respectively, as
$d_{k}F(\bold u)=\sum^n_{i=1} \frac {\partial_k F(\bold u)}{\partial_k x_i}d_{k}x_{i}$).
\par
The definition immediately implies that partial derivatives 
modulo $p^k$ of the function $F$ are defined up to the $p$-adic integer
summand which $p$-adic norm does not exceed $p^{-k}$. In cases when all partial derivatives
modulo $p^k$ at all points of  
$\Bbb Z_p^{(n)}$ are
$p$-adic integers, we say that the function 
$F$ has {\it integer-valued derivative modulo} $p^k$; 
in these cases we can associate to each partial derivative modulo $p^k$
a unique element of the ring $\Bbb Z/p^k$, 
and a Jacobi matrix modulo $p^k$ 
at each point $\bold u\in \Bbb Z_p^{(n)}$ 
thus can be considered as a matrix over a ring $\Bbb Z/p^k$. 
\par
Under the latter agreement the `rules of differentiation
modulo $p^k$' have the same (up to congruence modulo $p^k$ instead of equality)
form as for usual differentiation.
For instance, if both functions 
$G\colon{\Bbb Z}^{(s)}_{p}\rightarrow {\Bbb Z}^{(n)}_{p}$ and
$F\colon{\Bbb Z}^{(n)}_{p}\rightarrow {\Bbb Z}^{(m)}_{p}$ 
are differentiable modulo 
$p^{k}$ at the points, respectively, $\bold v=(v_{1},\ldots  ,v_{s})$
and $\bold u=G(\bold v)$, and their partial derivatives modulo $p^{k}$ at
these points are $p$-adic integers, then a composition 
$F\circ G\colon{\Bbb Z}^{(s)}_{p}\rightarrow {\Bbb Z}^{(m)}_{p}$ 
of these functions is uniformly differentiable modulo $p^{k}$ at the point
$\bold v$, all its partial derivatives 
modulo $p^{k}$ at this point are $p$-adic integers, and 
$(F\circ G)^\prime_k (\bold v)\equiv G^\prime_k (\bold v) F^\prime_k (\bold u)\pmod
{p^k}$.
\par
By the analogy with classical case we define for the function $F$ a notion
of {\it uniform differintiability modulo $p^k$ on $\Bbb Z_p^{(n)}$}; the
least  
$K\in\Bbb N$ such that $(\heartsuit)$ holds simultaneously for all 
$\bold u \in \Bbb Z_p^{(n)}$, whereas $\| h_{i}\| _{p}\le     p^{-K}$, $(i=1,2,\ldots  ,n)$,
is denoted via $N_k(F)$. The latter number plays an important role in
further coniderations.
\par
We recall that accordingly to  2.12 of \cite{11} all  partial derivatives
 modulo $p^k$ of the uniformly differentiable modulo $p^k$ function $F$
 are periodic functions with period  
$p^{N_k(F)}$. This in particular implies that each partial derivative modulo
$p^k$ can be considered as a function defined on $\Bbb Z/p^{N_k(F)}$. 
Moreover, if 
$F=(f_{1},\ldots , f_{m})\colon{\Bbb N}^{(n)}_{0}\rightarrow {\Bbb N}^{(m)}_{0}$  
could be continued to a function on the $\Bbb Z_p^{(n)}$, which is uniformly differentiable modulo $p^k$ on the
whole space $\Bbb Z_p^{(n)}$, this continuation could be done simultaneously 
with all its
(partial) derivatives modulo $p^k$.
\par
Here and after in this section let  $F=(f_{1},\ldots , f_{m})\colon{\Bbb Z}^{(n)}_{p}\rightarrow {\Bbb Z}^{(m)}_{p}$  
and $f\colon{\Bbb Z}^{(n)}_{p}\rightarrow {\Bbb Z}_{p}$ be functions, 
which are uniformly differentiable on $\Bbb Z_p^{(n)}$  modulo $p$. This is relatively
weak restriction since all uniformly differentiable on $\Bbb Z_p^{(n)}$ functions,
as well as functions, which are uniformly differentiable on $\Bbb Z_p^{(n)}$
modulo $p^k$ for some $k\ge
1$, are uniformly differentiable on $\Bbb Z_p^{(n)}$ modulo $p$. 
\par
The examples of functions which are not uniformly differentiable on $\Bbb Z_p^{(n)}$,
yet are  uniformly differentiable on $\Bbb Z_p^{(n)}$ modulo $p$, are 
the function $f(x,y)=x\XOR y$ for $p=2$
and its corresponding analogs for $p\ne 2$; all partial derivatives modulo
$p$ of
these functions are congruent to 1 modulo $p$ at all points (see \cite{11}). 
Note by the way, that previously introduced  function 
$\bmod{\,p^n}\colon \Bbb Z_p\rightarrow\Bbb Z/p^n$, the `reduction modulo $p^n$',
is uniformly differentiable on  $\Bbb Z_p$ (its derivative is $0$ at all
points);
the function $f(x,y)=x\AND y$ is differentiable modulo  $2$ at no point
of $\Bbb Z_2^{(2)}$, yet it is uniformly differentiable with respect to 
$x$ for each  $y\in \Bbb Z$: its derivative is 0 for $y\ge 0$, and it is
1 in the opposite case.

It turnes out that properties of being asymptotically 
compatible or asymptotically measure-preserving impose certain restrictions
on $p$-adic norms of derivatives modulo $p$ of a given function.
\proclaim
{ 3.2 Proposition} { If the function $f\colon \Bbb Z_p \rightarrow \Bbb Z_p$ 
asymptotically preserves measure, then 
$\| f^{\prime}_1(u)\| _{p}\nobreak \geqslant \nobreak 1$ at all points
  $u\in {\Bbb Z}_{p}$.}
\endproclaim
\demo 
{ Proof} Since a derivative modulo $p^k$ of the function $f$ 
is periodic with period $p^{N_k(f)}$, it is sufficient to prove the proposition
assuming
$u\in {\Bbb N}_{0}$. Definition of differentiability modulo 
$p^k$
implies that for $K\ge N_{1}(f)$ and for  $u\in {\Bbb N}_{0}$ the congruence
$$
f(u+h)\equiv f(u)+hf_1^{\prime}(u)\pmod{p^{K+1}}
\eqnum{1}
$$
\noindent holds as soon as $\| h\| _{p}\le   p^{-K}$. Assuming 
$\| f_1^{\prime}(u)\| _{p}<1$ for some $u\in {\Bbb N}_{0}$, the condition
$f_1^{\prime}(u)\equiv 0\pmod p$
and congruence (1) imply that  $f(u+p^{K})\equiv f(u)\pmod{p^{K+1}}$. 
The latter congrunce means that for all
$K\ge N_{1}(f)$, such that $u+p^{K}\le    p^{K+1}-1$, the function 
$f$ is not bijective modulo $p^{K+1}$. A contradiction.\qed
\enddemo
\proclaim 
{ 3.3 Corollary} { If under the assumptions of 3.2 the function $f$ 
is uniformly differentiable, then
$\| f^\prime (u)\| _{p}\ge 1$ for all $u\in {\Bbb Z}_{p}$.}
\endproclaim
\demo 
{Proof} Definition of a derivative modulo $p$ immediately implies
that 
$$f_1^{\prime} (u)\equiv f^\prime (u)\pmod p$$ 
\noindent for all $u\in {\Bbb Z}_{p}$. Thus $f^\prime (u)=f_1^{\prime}(u)+ps(u)$ 
for a suitable function $s\colon{\Bbb Z}_{p}\rightarrow {\Bbb Z}_{p}$. 
Yet if $\| f_1^{\prime}(u)\| _{p}\ge 1$, then the latter equality 
obviously implies that $\| f^\prime (u)\| _{p}\ge 1$ 
by the properties of $p$-adic distance. 
Now the conclusion follows from 3.2.\qed
\enddemo
\par
The inverse of 3.2 is not true: an obvious counterexample gives the function
$\frac {x^2 - x}{2}$  on ${\Bbb Z}_{2}$. It vanishes both at  0 and at
1, but the $2$-adic norm of its derivative is 2 everywhere on ${\Bbb Z}_{2}$. Nevertheless, functions of this
kind are locally injective. Namely, the following is true:
\proclaim
 { 3.4 Proposition} {If the function$f\colon{\Bbb Z}_{p}\rightarrow {\Bbb Z}_{p}$ 
is uniformly differentiable modulo $p$, and if 
$\| f_1^{\prime}(u)\| _{p}\ge 1$, then a space ${\Bbb Z}_{p}$ can be represented
as a disjoint union of a finite number of open 
{\rom (}and simultaneously closed{\rom )} balls $U$, for which the following holds: if
$a,b\in U$, $k\ge N_{1}(f)$ and $a\not\equiv b\pmod{p^{k}}$, then
$f(a)\not\equiv f(b)\pmod{p^{k}}$.}
\endproclaim
\demo 
{Proof} Consider a union
$${\Bbb Z}_{p}=\bigcup\limits_{a=0}^{p^N - 1}(a+p^{N}{\Bbb Z}_{p}),$$ 
where $N=N_{1}(f)$. Each set
$U=a+p^{N}{\Bbb Z}_{p}$ is an open (and at the same time closed) ball of
radius $p^{-N}$ 
(see \cite 3). 
Let $u,v\in U$, and let
$u\neq v$. Then $v=u+h$, where $\| h\| _{p}=p^{-K}$ for a suitable
positive integer rational $K\ge N$.
The definition of differentiability modulo $p$ implies that
$$
f(u+h)\equiv f(u)+hf_1^{\prime }(u)\pmod{p^{K+1}}.
\eqnum{1}$$
Thus, if  $f(u)\equiv f(v)\pmod{p^{K}}$, then (1) implies that
$\|f_1^\prime (u)\| _{p}=p^{-1}<1$. 
A contradiction.\qed
\enddemo
\par
The proposition 3.4 implies that if the $p$-adic norm of a uniformly differentiable
modulo $p$ function is not less then 1 everywhere 
on $\Bbb Z_p$, then this function might `glue together modulo $p^k$' for sufficiently
large $k$ only points which lie in distinct balls from the statement of 3.4. From here it
follows
\proclaim
{3.5 Proposition} { Let a function $f\colon{\Bbb Z}_{p}\rightarrow {\Bbb Z}_{p}$ be
uniformly differentiable modulo $p$ on $\Bbb Z_p$. Then $f$ asymptotically
preserves measure iff the following condition hold simultaneously:
\roster
\item $\| f_1^{\prime }(u)\| _{p}\ge 1$ at all points $u\in {\Bbb Z}_{p}$;
\item  $f(a)\not\equiv f(b)\pmod{ p^{n}}$ for all 
$n,a,b\in {\Bbb N}_{0}$ such that
$\| a-b\| _{p}\ge p^{-N_{1}(f)}$ и $0\le a,b\le p^{n}-1.$\qed\endroster}
\endproclaim
\par
A. A. Nechaev (private communication) noticed that the function $f(x)=\frac {x^2 + x}{2}$ 
on ${\Bbb Z}_{2}$ asymptotically preserves measure (this also follows from
3.5). Thus, if a compatible function $g\colon{\Bbb Z}_{2}\rightarrow {\Bbb Z}_{2}$
asymptotically preserves measure (all these functions are characterized
in 2.2), then a composition 
$h(x)=g(f(x))$ is uniformly differentiable modulo $p=2$ and asymptotically
measure-preserving function, and
$\| g_1^\prime (u)\| _{2}=2$ 
at all points
$u\in {\Bbb Z}_{2}$. There are no other functions 
$f\colon{\Bbb Z}_{p}\rightarrow {\Bbb Z}_{p}$, which are uniformly differentiable
modulo $p$, asymptoticaly preserve measure, and which derivatives modulo
$p$ have norms not less then 1 everywhere on $\Bbb Z_p$, \cite {10}.
The proof of the latter statement involves not only $p$-adic tools,
but algebraic geometry techniques as well.
\par 
The latter notice illustrates the fact that the second condition of the
criterion 3.5 is rather difficult to verify since one has to calculate
values of a function at infinite number of points. However, the problem
might be simplified by imposing certain restrictions on the function under
study. Namely, we will assume additionally that $f$ maps each ball of radius
$p^{-M}$ (with $M\ge N_{1}(f)$) into a ball of radius 
$p^{-M}$ (consequently, $f$ is asymptotically compatible). 
This restriction is equvalent to the property of derivative modulo $p$ to
be
integer-valued
everywhere on $\Bbb Z_p$.
\proclaim
{ 3.6 Proposition} { If for some $M\ge N_{1}(f)$ a uniformly differentiable
modulo $p$ function $f$ maps each ball of radius $p^{-M}$ into a ball of
radius $p^{-M}$, then
$f_1^\prime (a)\in {\Bbb Z}_{p}$ for all $a\in {\Bbb Z}_{p}$. Vice versa,
each uniformly differentiable modulo $p$ function, which has an integer-valued
derivative modulo $p$ everywhere on $\Bbb Z_p$, maps each ball of radius
$p^{-M}$
into a ball of radius $p^{-M}$ for all $M\ge N_{1}(f)$.}
\endproclaim
\demo
{Proof} If $M\ge N_{1}(f)$ and $\| h\| _{p}\le   p^{-M}$, then the definition
of uniform differentiability modulo $p^k$ (see 2.4 of \cite{11}) implies
that
$$
f(u+h)\equiv f(u)+hf_1^\prime (u)\pmod{p^{M+1}}
\eqnum{1}$$
\noindent for all $u\in {\Bbb Z}_{p}$. On the other hand, the inclusion
$f(a+p^{M}{\Bbb Z}_{p})\subseteq f(a)+p^{M}{\Bbb Z}_{p}$
implies that
$$
\| f(u+h)-f(u)\| _{p}\le     p^{-M}
\eqnum{2}$$
\noindent for all $h$ with $\| h\| _{p}\le   p^{-M}$. Comparing (1) and
(2) we see that
$\| f_1^\prime (u)\| _{p}\le   1$. The inverse statement is equivalent
to the asymptotic compatibility of $f$ (see 2.10 of \cite{11}).\qed
\enddemo
\par
Henceforth in the section we additionally assume that 
$f$ and $F$ have integer-valued derivatives modulo $p$. In particular,
this implies that both $f$ and $F$ are asymptotically compatible 
(see 2.10 and 2.11 of \cite{11}).
Now we state necessary and sufficient conditions 
the function $F$ must satisfy to be measure-preserving, and sufficient conditions
for equiprobability of $F$. 
\proclaim
{ 3.7 Theorem} {Let a function $F=(f_{1},\ldots  ,f_{m})\colon{\Bbb Z}^{(n)}_{p}\rightarrow {\Bbb Z}^{(m)}_{p}$ 
be uniformly differentiable modulo $p$ and let all its partial derivatives
modulo $p$ be integer-valued on $\Bbb Z_p$. Then $F$ is asymptotically
equiprobable if it is equiprobable modulo $p^{k}$ for some
$k\ge N_{1}(F)$ and the rank of its Jacobi matrix $F_1^\prime (u)$ modulo
$p$ is exactly $m$ at all points  $\bold u=(u_{1},\ldots  ,u_{n})\in ({\Bbb Z}/p^{k})^{(n)}$.}
\endproclaim
\demo 
{Proof} For $\xi \in ({\Bbb Z}/p^{s})^{(m)}$ denote  $F^{-1}_{s}(\xi )=\{\gamma \in ({\Bbb Z}/p^{s})^{(n)}:F(\gamma )\equiv \xi \pmod{p^{s}}\}$. 
Let
$s\ge k\ge N_{1}(F)$. Since $F$ is asymptotically compatible, and hence
$F$ is a sum of a compatible function and a periodic function with period
 $p^{N_1(F)}$ (see 2.10
of \cite{11}), we conclude that if $\eta \in F^{-1}_{s+1}(\xi )$,  then  $\bar{\eta }\in F^{-1}_{s}(\bar{\xi })$.  
Here, in accordance with our agreement in the introduction,  $\bar{\alpha }=(\bar{\alpha }_{1},\ldots  ,\bar{\alpha }_{m})\in ({\Bbb Z}/p^{s})^{(m)}$ 
stands for $\alpha\bmod{p^s}=(\alpha_1\bmod{p^s},\ldots,\alpha_m\bmod{p^s})$, where $\alpha =(\alpha _{1},\ldots  ,\alpha _{m})\in ({\Bbb Z}/p^{s+1})^{(m)}$. 
Put $\lambda =\bar{\eta }+p^{s}\sigma\in(\Bbb Z/p^{s+1})^{(n)}$, where $\sigma \in ({\Bbb Z}/p)^{(n)}$. 
In view of the uniform differentiability of 
the function $F$ modulo $p$ (see $(\heartsuit )$), we have
$$
F(\lambda )\equiv F(\eta )+p^{s}\sigma F_1^{\prime}(\bar{\eta })\pmod{p^{s+1}}.
\eqnum{1}$$
\noindent  Since  $F(\bar{\eta })\equiv \bar{\xi }+p^{k}\beta \pmod{p^{s+1}}$
and $ \xi =\bar{\xi }+p^{s}\gamma $ 
for suitable $\beta ,\gamma \in ({\Bbb Z}/p)^{(m)}$, then (1)
implies that $\lambda \in F^{-1}_{s+1}(\xi )$ iff $\bar{\lambda }\in F^{-1}_{s}(\xi )$ 
(i.e., $\bar{\eta }\in F^{-1}_{s}(\xi ))$ and $\alpha $
satisfies the following linear system over a field ${\Bbb Z}/p$:
$$
\beta +\alpha F_1^{\prime}(\bar{\eta })=\gamma .
\eqnum{2}$$
\noindent Thus, if columns of the matrix $F_1^\prime (\bar{\eta })$ are
linearly independent over ${\Bbb Z}/p$, then linear system (2) has exactly
$p^{n-m}$ distinct solutions for arbitrary $\beta ,\gamma \in ({\Bbb Z}/p)^{(m)}$. 
From here it follows that
$$
\vert F^{-1}_{s+1}(\xi )\vert =\vert F^{-1}_{s}(\xi )\vert p^{n-m}.
\eqnum{3}$$
\noindent Hence, if $F$ is equiprobable modulo $p^{k}$ (i.e., if $\vert F^{-1}_{s}(\bar{\xi })\vert $ does
not depend on $\bar{\xi })$
and a rank of a matrix $F_1^\prime (\bar{\eta })$ is $m$, then (3) implies
that $F$ 
is equiprobable modulo $p^{s+1}$.\qed
\enddemo
\proclaim
{ 3.8 Corollaries} { $1^{\circ} $ Under the assumptions of theorem 3.7 let $m=1$.  Then
$F$ if asymptotically equiprobable if $F$ is equiprobable
modulo $p^{k}$ for some $k\ge N_{1}(F)$, and differential $d_{1}F$  modulo
 $p$  of the function $F$
vanishes at no point of $({\Bbb Z}/p^{k})^{(n)}$.
\par
$2^{\circ }$ Let $f(x_{1},\ldots  ,x_{n})$ be a polynomial 
with integer $p$-adic coefficients and in variables
 $x_{1},\ldots  ,x_{n}$. A polynomial 
$f$ is equiprobable if it is equiprobable modulo
$p$ and all its partial derivatives vanishes simultaneously modulo $p$ at no point
of $({\Bbb Z}/p)^{(n)}$ {\rm (i.e., are simultaneously congruent modulo $p$
nowhere)}.}
\endproclaim
\demo 
{Proof} The assertion $1^{\circ }$ trivially follows from  3.7. In turn,  $2^{\circ }$ immediately
follows from $1^{\circ }$, since for all $f\in {\Bbb Z}[x_{1},\ldots  ,x_{n}]$ 
holds $N_{1}(f)\le 1$. We have only to prove the latter inequality.
\par
By Taylor formula,
$$
f(x_{1}+h_{1},\ldots  ,x_{n}+h_{n})=f(x_{1},\ldots  ,x_{n})+\sum\limits_{i=1}^{n}h_{i}\frac{\partial f}{\partial x_{i}}+Q
\eqnum{1}$$
\noindent where $Q\in {\Bbb Z}[x_{1},\ldots  ,x_{n},h_{1},\ldots  ,h_{n}]$, 
and each monomial in a canonic representation of
the polynomial $Q$
is of degree not less then 2 with respect to variables $h_{1},\ldots  ,h_{n}$. Since
$\| (h_{1},\ldots  ,h_{n})\| _{p}=p^{-s}$, where $s\ge 1$, for all values
of
$x_{1},\ldots  ,x_{n}$ we have $Q\equiv 0\pmod{p^{2s}}$. In view of (1)
this
proves the inequality.\qed
\enddemo
\par
For $m=n$ the above stated sufficient conditions of asymptotical equiprobability 
occur to be necessary as well.
\proclaim
{ 3.9  Theorem}  { Uniformly differentiable modulo $p$ function
$$F=(f_{1},\ldots  ,f_{m})\colon{\Bbb Z}^{(n)}_{p}\rightarrow {\Bbb Z}^{(n)}_{p}$$
with integer-valued derivatives modulo $p$, asymptotically preserves
measure if and only if it is bijective modulo 
$p^{N_{1}(F)}$ and its Jacobian modulo $p$ vanishes at no point of 
$({\Bbb Z}/p^{N_1(F)})^{(n)}$ {\rom (}equivalent condition: iff $F$ is bijective
modulo $p^{N_1{(F)+1}}$\rom ).}
\endproclaim
\demo 
{Proof} If $F$ is bijective modulo  $p^{N_{1}(F)}$ and its Jacobian modulo
$p$ vanishes nowhere, then in view of  3.7  $F$ is asymptotically equiprobable,
hence, asymptotically preserves measure, since $m=n$.
\par
Vise versa, let $F$ asymptotically preserves measure, i.e., let $F$ be
bijective modulo $p^{k}$ for all $k\ge N$, where $N$ is some positive integer
rational.
Now take $k\ge \max\{N,N_{1}(F)\}$, then the definition of uniform differentiability
modulo $p$ implies that
$$
F(u+p^{k}\alpha )\equiv F(u)+p^{k}\alpha F_1^{\prime}(u)\pmod{p^{k+1}}
\eqnum{1}$$
\noindent for all $u,\alpha \in {\Bbb Z}_{p}$. Here $F_1^\prime (u)$ is
an $n\times n$ matrix over a field ${\Bbb Z}/p$. If $\det  F_1^\prime (u)\equiv 0\pmod{p}$
for some $u\in {\Bbb Z}^{(n)}_{p}$ (or, the same, for 
some $u\in \{0,1,\ldots  ,p^{N_{1}(F)}-1\}^{(n)}$ in view of the periodicity
of partial derivatives modulo $p$), 
then there exists $\alpha \in \{0,1,\ldots  ,p-1\}^{(n)}, \alpha \not\equiv (0,\ldots  ,0)\pmod{p}$, 
such that
$\alpha F_1^\prime(u)\equiv  (0,\ldots,0)\pmod{p}$. But then (1) implies that 
$F(u+p^{k}\alpha )\equiv F(u)\pmod{p^{k+1}}$. The latter contradicts the
bijectivity modulo $p^{k+1}$ of the function
$F$, since for $u\in \{0,1,\ldots  ,p^{N_{1}(F)}-1\}^{(n)}$ we have
$u,u+p^{k}\alpha \in \{0,1,\ldots  ,p^{k+1}-1\}^{(n)}$ and $u+p^{k}\alpha \neq u$.
\par
Now we prove the criterion in the equivalent form. Let $F$ be bijective 
modulo $p^{N_{1}(F)}$. Then assuming $k=N_{1}(F)$ in the above argument,
we conclude that $\det F_1^\prime(u)\not\equiv 0\pmod{p}$ for all $u\in {\Bbb Z}^{(n)}_{p}$. 
According to 3.7, this implies that $F$ asymtotically preserves measure.
\par
Let $F$ asymptotically preserves measure, and let it be not 
bijective modulo $p^{k}$ for some $k\ge N_{1}(F)$. We prove that in this
case $F$ is not bijective modulo $p^{k+1}$. 

Choose  $u,v\in \{0,1,\ldots  ,p^{k}-1\}^{(n)}$ 
such that $u\neq v$ и
$F(u)\equiv F(v)\pmod{p^{k}}$. Then either $F(u)\equiv F(v)\pmod{p^{k+1}}$ (i.e., $F$ 
is not bijective modulo $p^{k+1})$, or $F(u)\not\equiv F(v)\pmod{p^{k+1}}$. 
Yet in the latter case we have
$F(u)\equiv F(v)+p^{k}\alpha \pmod{p^{k+1}}$ for some $\alpha \in \{0,1,\ldots  ,p-1\}^{(n)}$, $\alpha \not\equiv (0,\ldots  ,0)\pmod{p}$.
Consider $u_{1}=u+p^{k}\beta $, where $\beta \in \{0,1,\ldots  ,p-1\}^{(n)}$
with $\beta \not\equiv (0,\ldots  ,0)\pmod{p}$ and
$\beta F_1^\prime(u)+\alpha \equiv (0,\ldots  ,0)\pmod{p}$. Such $\beta $ exists,
since $F$ asymptotically preserves measure and, consequently,  
$\det F_1^\prime(u)\not\equiv 0\pmod{p}$, as
it have been proven already.
Now the definition of uniform differentiability modulo $p$ implies that 
$$
F(u+p^{k}\beta )\equiv F(u)+p^{k}\beta F_1^{\prime}(u)\equiv F(v)+p^{k}\alpha +p^{k}\beta F_1^{\prime}(u)\equiv F(v)\pmod{p^{k+1}},
\eqnum{2}$$
\noindent where $u+p^{k}\beta \in \{0,1,\ldots  ,p^{k+1}-1\}^{(n)}$ and $u+p^{k}\alpha \neq v$ (since $u\neq v$). 
Thus (2) in combination with our assumption imply that $F$ is not bijective
modulo 
$p^{k+1}$. Applying this argument sufficient number of times, we conclude
that $F$ is not bijective 
modulo $p^{s}$ for all $s\ge k$. But at the same time $F$ asymptotically
preserves measure. A contradiction.\qed
\enddemo
\proclaim 
{ 3.10 Corollaries} { $1^{\circ }$ If $n=1$ within assumptions of the
theorem 3.9, then $F$ asymptotically preserves measure iff it is bijective
modulo  $p^{N_{1}(F)}$ and its derivative modulo $p$ vanishes at no point
of $\{0,1,\ldots  ,p^{N_{1}(F)}-1\}$.
\par
$2^{\circ }$ {\rm (cf. \cite {8, Ch. 4, sections 4--5})} Let
$F=(f_{1},\ldots ,f_{m})\colon{\Bbb Z}^{(n)}_{p}\rightarrow {\Bbb Z}^{(n)}_{p}$, 
where 
$f_{i}(x_{1},\ldots  ,x_{n})\in {\Bbb Z}_{p}[x_{1},\ldots  ,x_{n}]$, $i=1,2,\ldots  ,n$.
Then $F$ preserves measure iff  $F$ is bijective
modulo $p$ and
$\det  F^\prime (u)\not\equiv 0\pmod{p}$ for all $u\in \{0,1,\ldots  ,p-1\}^{(n)}$ {\rom (}equivalent
statement: iff
$F$ is bijective modulo $p^{2}${\rom )}.
\par
$3^{\circ }$ Let $A=\langle {\Bbb Z}_{p};\Omega \rangle $ be a universal
algebra of finite signature
$\Omega $, and let all operations of $\Omega $ are uniformly differentiable
modulo  $p$ and have integer-valued
derivatives modulo $p$. Then 
a polynomial over $A$ defines an asymptotically measure-preserving function
iff it is bijective modulo $p^{k(A)}$, where $k(A)=\max\{N_{1}(\omega ):\omega \in \Omega \}+1$.}
\endproclaim
\demo
{Proof} Assertion $1^{\circ }$ trivially follows from 3.9. Assertion $2^{\circ }$ 
holds in view of  3.9, since $N_{1}(F)\le 1$ (see proof of the corollary
3.8). A compositin $F\circ G$ 
of functions $F$ and $G$, which are both uniformly differentiable modulo
 $p$ and have 
integer-valued derivatives modulo $p$, is uniformly differentiable modulo
$p$ function, which has an integer-valued derivative modulo $p$, and 
$N_{1}(F\circ G)\le \max\{N_{1}(F),N_{1}(G)\}$. The latter proves
assertion $3^{\circ }$.
\qed
\enddemo
Comparing statements 3.7 and 3.9 one may ask a natural question
whether sufficient conditions of 3.7 are necessary. The answer is
negative: the results of \cite 9 make it possible to construct 
the following counterexample. 
\par
Consider a function
$f(x,y)=2x+y^{3}$ on ${\Bbb Z}_{2}$.
As $f$ is a polynomial over ${\Bbb Z}$, then it is uniformly differentiable,
has integer-valued derivatives, and $df=2dx+3y^{2}dy$. So,
$df\equiv 0\pmod{2}$ if $y\equiv 0\pmod{2}$. Nevertheless, $f$ induces an equiprobable
function $({\Bbb Z}/2^{n})^{(2)}\rightarrow {\Bbb Z}/2^{n}$ for every $n=1,2,\ldots  $. Here is a proof.
\par
For $n=1$ we have that $f(x,y)\equiv y\pmod{ 2}$ is an an equiprobable function on ${\Bbb Z}/2$.
Let $n>1$. We will show that for every $z\in {\Bbb Z}/2^{n}$ there exist exactly $2^{n}$ pairs
$(x,y)$, such that $f(x,y)\equiv z\pmod{2^{n}}$ and $(x,y)\in \{0,1,\ldots  ,2^{n}-1\}^{(2)}$.
\par
In fact, if $z=1+2r$ for some $r\in \{0,1,\ldots  ,2^{n-1}-1\}$, then it follows that
$y=1+2k$ for some $k\in \{0,1,\ldots  ,2^{n-1}-1\}$. So $2x+(1+2k)^{3}\equiv 1+2r\pmod{2^{n}}$ implies
$x+3k+6k^{2}+4k^{3}\equiv r\pmod{2^{n-1}}$. The left hand part  of the latter congruence is
a polynomial $\phi (x,k)$ in $x,k$. It is equiprobable in view of $3.8, 2^{\circ }$,
since $d\phi \equiv dx+dk\pmod{2}$ (and hence this differential  vanishes modulo
2 nowhere) and $\phi \equiv x+k\pmod{2}$ is obviously an equiprobable modulo 2
function. This implies that the congruence $\phi (x,k)\equiv r\pmod{2^{n-1}}$ in
unknowns $x,k$ has exactly $2^{n-1}$ solutions in $\{0,1,\ldots  ,2^{n-1}-1\}^{(2)}$.
\par
If $z=2r$ for some $r\in \{0,1,\ldots  ,2^{n-1}-1\}$, then it follows that $y=2k$ for some
$k\in \{0,1,\ldots  ,2^{n-1}-1\}$; consequently, the congruence $f(x,y)\equiv z\pmod{2^{n}}$
implies the congruence $x+4k^{3}\equiv r\pmod{2^{n-1}}$. Again the function $\psi (x,k)$ in
the left hand part of the latter congruence is equiprobable in view of
$3.8,2^{\circ }$, since $d\psi \equiv dx\pmod{2}$ vanishes modulo 2 at no point of $({\Bbb Z}/2)^{(2)}$ and
$\psi \equiv x\pmod{2}$ is equiprobable modulo 2. From here, using an argument similar
to one of the previous case, we conclude that the congruence
$f(x,y)\equiv 2r\pmod{2^{n}}$ in unknowns $x,y$ has exactly $2^{n}$ solutions in
$\{0,1,\ldots  ,2^{n}-1\}^{(2)}$. Thus, $f$ is equiprobable.
\par
Now we are to begin a study of asymptotically ergodic functions in the class
of all uniformly differentiable modulo $p$ functions, which have integer-valued 
derivatives modulo $p$. It turnes out that these functions could be in
one variable only. To be more exact, the following theorem is true.
\proclaim
{ 3.11 Theorem} { Let a function $F=(f_{1},\ldots  ,f_{n})\colon{\Bbb Z}^{(n)}_{p}\rightarrow {\Bbb Z}^{(n)}_{p}$ 
be uniformly differentiable
modulo $p$ and asymptotically ergodic, and let it have integer-valued
derivatives modulo $p$. Then $n=1$.}
\endproclaim
\par We will need two lemmata. 
\proclaim 
{ 3.12 Lemma} { Let a function $f\colon{\Bbb Z}^{(n)}_{p}\rightarrow {\Bbb Z}_{p}$ be
uniformly differentiable modulo $p$, let it
have integer-valued derivatives modulo $p$, and let it vanish modulo $p^{k}$ 
{\rm (i.e., let it be congruent 0 modulo $p^{k}$)} for
some $k>N_{1}(f)$ at all points of ${\Bbb Z}^{(n)}_{p}$. Then each partial derivative modulo $p$
of the function $f$ vanishes modulo $p$ at all points of ${\Bbb Z}^{(n)}_{p}$.}
\endproclaim
\demo 
{Proof of the lemma 3.12} Each function $g_{i}(x_{0},x_{1},\ldots  ,x_{n})=x_{i}+x_{0}f(x_{1},\ldots  ,x_{n})$ for arbitrary
values of $x_{0},x_{1},\ldots  ,x_{i-1},x_{i+1},\ldots  ,x_{n}$ is a bijective modulo $p^{k}$ function in
variable $x_{i}, (i=1,2,\ldots  ,n)$. As $k>N_{1}(g_{i})=N_{1}(f)$, then according to 3.9, $g_{i}$
asymptotically preserves measure, and thus its derivative modulo $p$ vanishes
at no point of ${\Bbb Z}_{p}$. Moreover, the following is true:
$$
\frac{\partial_1}{\partial_1 x_{i}}g_{i}(u_{0},\ldots  ,u_{n})=1+u_{0}\cdot \frac{\partial_1}{\partial_1 x_{i}}f(u_{1},\ldots  ,u_{n})\not\equiv 0\pmod{p}
\eqnum{1}$$
\noindent for all $u_{0},\ldots  ,u_{n}\in {\Bbb Z}_{p}$. If
 $$
\frac{\partial_1}{\partial_1 x_{i}}f(u_{1},\ldots  ,u_{n})\equiv d\not\equiv 0\pmod{p}$$
 for some
$u_{1},\ldots  ,u_{n}\in {\Bbb Z}_{p}$, then choosing $u_{0}$ such that $u_{0}d\equiv -1\pmod{p}$ we get a
contradiction to (1). This proves the lemma.\qed 
\enddemo
\proclaim
{ 3.13 Lemma} { Let a function $H\colon{\Bbb Z}^{(n)}_{p}\rightarrow {\Bbb Z}^{(n)}_{p}$ 
be uniformly differentiable modulo $p$,
and let it have integer-valued derivatives modulo $p$. If $H$ is
bijective modulo $p^{k}$ and if $H$ induces a trivial  permutation
modulo $p^{k-1}$ {\rm (i.e., an identity transformation of $({\Bbb Z}/p^{k-1})^{(n)}$)}  
for some $k>N_{1}(H)+1$, then $H$ induces modulo $p^k$ {\rm (i.e., on
$({\Bbb Z}/p^{k})^{(n)}$)} either a
trivial permutation, or a permutation of order $p$.}
\endproclaim
\demo 
{ Proof of the lemma 3.13} Let $G$ be an arbitrary function which satisfies 
assumptions of the 
lemma, and let $N_{1}(G)=N_{1}(H)$. Represent both $H$ and $G$ in the following
form:
$$
H(x_{1},\ldots  ,x_{n})=(x_{1},\ldots  ,x_{n})+U(x_{1},\ldots  ,x_{n});
$$
$$
G(x_{1},\ldots  ,x_{n})=(x_{1},\ldots  ,x_{n})+V(x_{1},\ldots  ,x_{n}).
$$
\noindent Then both $U$ and $V$ are uniformly differentiable modulo $p$, have
integer-valued derivatives modulo $p$, and $N_{1}(U)=N_{1}(V)=N_{1}(H)$. Moreover,
both $U$ and $V$ vanish modulo $p^{k-1}$   on  ${\Bbb Z}^{(n)}_{p}$,
for  $k-1>N_{1}(H)$. Then
lemma 3.12 implies that $U_1^{\prime} =V_1^{\prime} =0$ at all  points  of  ${\Bbb Z}^{(n)}_{p}$.  
As $\| U\| _{p}\le p^{-k+1}$
and  $\| V\| _{p}\le p^{-k+1}$  everywhere  on  ${\Bbb Z}^{(n)}_{p}$,  then,  
implying 2.4, for all
$h_{1},\ldots  ,h_{n}\in {\Bbb Z}_{p}$ we  obtain, consequently, that
$$
\displaylines{H(G(h_{1},\ldots  ,h_{n}))=H((h_{1},\ldots  ,h_{n})+V(h_{1},\ldots  ,h_{n}))\hfill\cr
\equiv H(h_{1},\ldots  ,h_{n})+V(h_{1},\ldots  ,h_{n})H_1^{\prime} (h_{1},\ldots  ,h_{n})\cr
\equiv H(h_{1},\ldots  ,h_{n})+V(h_{1},\ldots  ,h_{n})+V(h_{1},\ldots  ,h_{n})U_1^{\prime} (h_{1},\ldots  ,h_{n})\cr
\hfill\equiv (h_{1},\ldots  ,h_{n})+U(h_{1},\ldots  ,h_{n})+V(h_{1},\ldots  ,h_{n})\pmod{p^k}.\cr}
$$ 
\noindent This implies, in particular,
that for all $s\in {\Bbb N}$ the following congruence holds:
$$
\displaylines{H^{s}(h_{1},\ldots  ,h_{n})=\underbrace{H(\ldots H}_{s \;\text{times}}(h_{1},\ldots  ,h_{n})\ldots  )\hfill\cr
\hfill\equiv (h_{1},\ldots  ,h_{n})+sU(h_{1},\ldots  ,h_{n})\pmod{ p^{k}}.\cr}
$$
\par 
As $U$ vanishes modulo $p^{k-1}$ everywhere, then the latter congruence
implies that $H^{p}(h_{1},\ldots  ,h_{n})\equiv (h_{1},\ldots  ,h_{n})\pmod{ p^{k}}$ for all $h_{1},\ldots  ,h_{n}\in {\Bbb Z}_{p}$. This
proves the lemma.\qed
\enddemo
\demo
{Proof of the theorem 3.11} Choose $k>N_{1}(F)+1$ such that $F$ is
transitive modulo $p^{n}$ for all $n\ge k-1$. The function $F$ induces a permutation
on $({\Bbb Z}/p^{k})^{(n)}$ which is denoted as $\sigma _{k}(F)$. Consider a permutation
$\sigma =\sigma _{k}(F)^{p^{(k-1)n}}$. As $F$ is transitive modulo $p^{k}$, 
the order of $\sigma $ is $p^{n}$ 
(and hence $\sigma$ is not  trivial).
\par
On the other hand, $\sigma =\sigma _{k}(F^{p^{(k-1)n}})$. But $F^{p^{(k-1)n}}$ is bijective modulo
$p^{k}$ and induces a trivial permutation modulo $p^{k-1}$ (the latter assertion
follows from transitivity of $F$ modulo $p^{k-1}$). 
Since $\sigma$ is not trivial, in view of 3.13 the
order of $\sigma $ must be $p$. Yet, according to the previous argument, the order of 
$\sigma$ is $p^n$, so necessarily $n=1$.\qed
\enddemo
\par
It is still an open problem to characterize asymptotically ergodic
functions in the class of all uniformly differentiable modulo $p$ functions
which have integer-valued derivatives modulo $p$, but if we additionally
assume that the function is uniformly differentiable modulo $p^{2}$ and has
integer-valued derivative modulo $p^{2}$, the following description can
be obtained. The method we prove the next theorem is in fact a generalization
to $p$-adic case of the idea originally applied by M. V. Larin to description
of transitive modulo $n$ polynomials over $\Bbb Z$, \cite{15}.
\proclaim
{ 3.14 Theorem} {Let a function $f\colon{\Bbb Z}_{p}\rightarrow {\Bbb Z}_{p}$ 
be uniformly differentiable modulo $p^{2}$ and let it
have integer-valued derivative modulo $p^{2}$. Then $f$ is
asymptotically ergodic if and only if it is transitive modulo $p^{N_{2}(f)+1}$ for
odd prime $p$ or, respectively, modulo $2^{N_{2}(f)+2}$ for $p=2$.}
\endproclaim
\par We need the following 
\proclaim 
{3.15 Lemma} { Let a function $f\colon{\Bbb Z}_{p}\rightarrow {\Bbb Z}_{p}$ 
be uniformly differentiable modulo $p$, and let it have
integer-valued derivative modulo $p$. If $f$ is transitive
modulo $p^{k}$ for some $k>N_{1}(f)$, then $f$ induces on ${\Bbb Z}/p^{k+1}$ 
a permutation, which
is either a single cycle of length $p^{k+1}$, or a product of $p$ pairwise disjoint
cycles of length $p^{k}$ each.}
\endproclaim
\demo 
{Proof of the lemma 3.15} For $i=0,1,2,\ldots$ we denote  via
$x_{i}=\delta _{i}(x)\in \{0,1,\ldots  ,p-1\}$ a value of $i$th digit
in canonic representation of 
 $p$-adic integer $x\in {\Bbb Z}_{p}$. Now the definition of uniform differentiability
 modulo $p$ implies that for an arbitrary $x\in {\Bbb Z}_{p}$ and $s\ge N_{1}(f)=N$ 
 there holds a congruence $f(x_0+x_1p+\cdots+x_{s-1}p^{s-1}+x_sp^s)\equiv
f(x_0+x_1p+\cdots+x_{s-1}p^{s-1})+x_sp^sf_1^{\prime}(x_0+x_1p+\cdots+x_{s-1}p^{s-1})\pmod {p^{s+1}}$.
The latter implies that
$$
\delta _{s}(f(x))\equiv \Phi _{s}(x_{0},\ldots  ,x_{s-1})+x_{s}f_1^{\prime}(x)\pmod{p},
\eqnum{1}$$
\noindent where $x_{i}=\delta _{i}(x)\in \{0,1,\ldots  ,p-1\}$ is the $i$-th $p$-adic digit of $x\in {\Bbb Z}_{p},
(i=0,1,2,\ldots   ); \Phi _{s}(x_{0},\ldots  ,x_{s-1})=\delta _{s}(f(x_{0}+x_{1}p+\cdots +x_{s-1}p^{s-1}))$.
\par
Since partial derivative  $f_1^{\prime}(x)$ modulo $p$ is periodic with
period $p^N$, it depends only on $x_{0},\ldots  ,x_{N-1}$, so (1) can be
represented in the form
$$
\delta _{s}(f(x))\equiv \Phi _{s}(x_{0},\ldots  ,x_{s-1})+x_{s}\Psi (x_{0},\ldots  ,x_{N-1})\pmod{p},
\eqnum{2}$$
\noindent where $\Psi (x_{0},\ldots  ,x_{N-1})=f_1^{\prime}(x)$. Applying
for the composition of functions
`rules of differentiation modulo $p^k$'  which were mentioned at 
the beginning of the section, we conclude that
for all $r=1,2,\ldots  $ the following congruence holds:
$$
(f^{r}(x))_1^{\prime}\equiv \prod \limits^{r-1}_{j=0} f_1^{\prime}(f^{j}(x))\pmod{p}.
\eqnum{3}$$
\noindent We recall that $f^{r}(x)=\underbrace{f(\ldots f}_{r\;\text{times}}(x)\ldots ), f^{0}(x)=x$. As $f$ is asymptotically
compatible, then transitivity of $f$ modulo $p^{k}$ for some $k\ge N$ implies
transitivity of $f$ modulo $p^{n}$ for all $k\ge n\ge N$ (see \cite{11},
theorems 2.10 and 1.4). Yet $f_1^{\prime}$ 
depends only on $x_{0},\ldots  ,x_{N-1}$, and $f$ is transitive modulo
$p^{N}$, so (3) implies that
$$
(f^{p^{n}}(x))_1^{\prime}\equiv \Biggl ( \prod \limits ^{p-1}_{u_{0},\ldots ,u_{N-1}=0} \Psi (u_{0},\ldots  ,u_{N-1})\Biggr )^{p^{n-N}}\pmod{p}.
\eqnum{4}$$
\noindent We denote the product in the brackets in the right hand part of (4)
as $\Pi $. Now, since  $f^{p^{n}}(x)$ is uniformly differentiable modulo $p$
and has
integer-valued derivative modulo $p$, in view of (2) and (4) we
conclude that
$$
\delta _{n}(f^{p^{n}}(x))\equiv \phi _{n}(x_{0},\ldots  ,x_{n-1})+x_{n}\Pi ^{p^{n-N}}\pmod{p},
\eqnum{5}$$
\noindent where $\phi _{n}(x_{0},\ldots  ,x_{n-1})=\delta _{n}(f^{p^{n}}(x_{0}+x_{1}p+\cdots  +x_{n-1}p^{n-1}))$. 
Since $f$ is a transitive
modulo $p^{n+1}$ function for $k\ge n\ge N$, the function $f^{p^{n}}$,
on the one hand, induces a trivial permutation modulo $p^n$, and on the
other hand, induces on each coset $a+p^{n}(\Bbb Z/p^{n+1})$ of the ring
$\Bbb Z/p^{n+1}$ a permutation, which is a cycle of length $p$.
This, in particular, means that the function in the right hand part of (5),
being considered as a function in variable $x_{n}$, must be a permutation,
moreover -- a cycle of length $p$ on $\{0,1,\ldots  ,p-1\}$. 
It is well known, however, that a polynomial $c+dy\in \Bbb Z[y]$ 
is transitive modulo  $p$ iff  $d\equiv 1\pmod p$ and 
$c\not\equiv 0 \pmod p$ (see e.g. \cite{2, Ch. 3, Theorem A}). 
This
implies, in particular, that $\Pi ^{p^{n-N}}\equiv 1\pmod{p}$, and hence 
$\Pi \equiv 1\pmod{p}$. Finally we obtain that
$$
\displaylines{f^{p^{k}}(x)\equiv f^{p^{k}}(x_{0}+x_{1}p+\cdots  +x_{k}p^{k})\hfill\cr
\hfill\equiv x_{0}+x_{1}p+\cdots  +x_{k-1}p^{k-1}+p^{k}(\phi _{k}(x_{0},\ldots  ,x_{k-1})+x_{k})\pmod{p^{k+1}}.\quad (6)\cr}
$$
\par The latter congruence implies that $f$ induces a permutation $\sigma $ modulo $p^{k+1}$.
Moreover, we assert that if $$\phi _{k}(x_{0},\ldots  ,x_{k-1})\not\equiv 0\pmod{p}$$ for some (equivalently, all) 
$x_{0},\ldots  ,x_{k-1}\in \{0,1,\ldots  ,p-1\}$, then $f$ is transitive modulo $p^{k+1}$; 
otherwise the permutation $\sigma $
is a product of exactly $p$ disjoint cycles of length $p^{k}$ each.
\par
To prove this assertion, consider some $u_{0},\ldots  ,u_{k}\in \{0,1,\ldots  ,p-1\}$ and
denote $C$ a cycle of the permutation $\sigma $ which contains the point
$u_{0}+u_{1}p+\cdots  +u_{k-1}p^{k-1}+x_{k}p^{k}\in {\Bbb Z}/p^{k+1}$. As $f$ is transitive modulo $p^{k}$ then (see
(6)) $p^k$ is a factor of $\vert C\vert $, the length of the cycle $C$. If $\phi _{k}(u_{0},\ldots  ,u_{k-1})\not\equiv 0\pmod{p}$,
then (6) implies that
$$
\displaylines{f^{p^{k}}(u_{0}+u_{1}p+\cdots  +u_{k-1}p^{k-1}+x_{k}p^{k})\hfill\cr
\hfill\not\equiv u_{0}+u_{1}p+\cdots  +u_{k-1}p^{k-1}+x_{k}p^{k}\pmod{p^{k+1}},
\quad(7)\cr}$$
\noindent i.e., that $\vert C\vert >p^{k}$. On the other hand, (6) implies that $\vert C\vert $ 
is a factor of $p^{k+1}$.
Finally we conclude that in this case $\vert C\vert =p^{k+1}$, i.e., $f$ is transitive
modulo $p^{k+1}$.
\par
If $\phi _{k}(u_{0},\ldots  ,u_{k-1})\equiv 0\pmod{p}$  holds for some $u_{0},\ldots  ,u_{k}\in \{0,1,\ldots  ,p-1\}$, then
this congruence holds for all $u_{0},\ldots  ,u_{k}\in \{0,1,\ldots  ,p-1\}$ (otherwise in
view of the previous case $f$ is transitive modulo $p^{k+1}$ and (7) holds for
all $u_{0},\ldots  ,u_{k}\in \{0,1,\ldots  ,p-1\}$ and the latter in view of (6) means that
$\phi _{k}(u_{0},\ldots  ,u_{k-1})\not\equiv 0\pmod{p}$, a contradiction). Then (6) implies that $\sigma ^{p^{k}}$ is
an identity permutation, i.e. $\vert C\vert =p^{k}$, as $p^{k}$ is a factor of $\vert C\vert$. 
This proves the lemma.\qed
\enddemo
\demo
{Proof of the theorem 3.14} 
During the proof
of the previous lemma we have established that if $f$ is transitive modulo
$p^{k}$ for some $k\ge N_{1}(f)$, then $f$ is transitive modulo $p^{n}$ for all $k\ge n\ge N_{1}(f)$. So
the `only if' part of the theorem is proved, as $N_{2}(f)+1>N_{1}(f)$.
\par
Now we have to prove that if $n\ge N_{2}(f)+1$ (resp., if $n\ge N_{2}(f)+2$ for
$p=2)$ and if $f$ is transitive modulo $p^{n}$, then it is transitive modulo $p^{n+1}$.
In view of lemma 3.15 it is sufficient to prove that for some $x\in {\Bbb Z}_{p}$ the
following condition holds:
$$
f^{p^{n}}(x)\not\equiv x\pmod{p^{n+1}}.
\eqnum{1}$$
\par As transitivity modulo $p^{n}$ implies transitivity modulo $p^{n-1}$, in
view of lemma 3.15 we have
$$
f^{p^{n-1}}(x)=x+p^{n-1}\xi (x),
\eqnum{2}$$
\noindent where $\xi \colon{\Bbb Z}_{p}\rightarrow {\Bbb Z}_{p}$ and $\xi (x)\not\equiv 0\pmod{p}$ for all $x\in {\Bbb Z}_{p}$ (otherwise 3.15 implies that
$f$ is not transitive modulo $p^{n}$, a contradiction to the assumption).
\par
Further, since $f$ is uniformly differentiable modulo $p^2$ and has integer-valued
derivative modulo $p^2$, then for all $r=1,2,\ldots$ a composition $f^{r}$ 
is uniformly differentiable modulo $p^{2}$
and has integer-valued derivative modulo $p^{2}$, and
$(f^{r}(x))_2^{\prime}\equiv \prod \limits^{r-1}_{j=0} f_2^{\prime}(f^{j}(x))\pmod{p^2}$ 
(see (3) of 3.15). Now, as $n-1\ge N_{2}(f)$, then taking into account
these considerations 
and an obvious (following from
(2)) equality
$f^{sp^{n-1}}(x)=f^{(s-1)p^{n-1}}(x+p^{n-1}\xi (x))$, where $s=1,2,\ldots  $,  
we successively calculate 
$$
\displaylines{f^{p^{n}}(x)\equiv f^{(p-1)p^{n-1}}(x)+p^{n-1}\xi (x)\prod \limits ^{(p-1)p^{n-1}-1}_{j=0}f_2^{\prime}(f^{j}(x))\hfill\cr
\hfill\equiv f^{(p-2)p^{n-1}}(x)+p^{n-1}\xi (x)\Biggl ( \prod \limits ^{(p-2)p^{n-1}-1}_{j=0} f_2^{\prime}(f^{j}(x))+ \prod \limits ^{(p-1)p^{n-1}-1}_{j=0}f_2^{\prime}(f^{j}(x))\Biggr )\cr
\hfill\equiv\ldots\equiv x+p^{n-1}\xi (x)\Biggl (1+\sum\limits^{p-1}_{i=1} \prod \limits ^{(p-i)p^{n-1}-1}_{j=0}f_2^{\prime}(f^{j}(x))\Biggr )\pmod{p^{n+1}}.\quad (3)\cr}$$
\par
Yet $f_2^{\prime}$  is a periodic function with period $p^{N_{2}(f)}$  
and $f$ is transitive modulo $p^{n-1}$, so  we conclude that for arbitrary
$i,j\in {\Bbb N}$ the following congruence holds:
 $$f_2^{\prime}(f^{j}(x))\equiv f_2^{\prime}(f^{j+ip^{n-1}}(x))\pmod{p^{2}}.$$
\par In view of the transitivity of $f$ modulo $p^{n-1}$ the latter congruence
implies that
$$
\prod \limits ^{(p-i)p^{n-1}-1}_{j=0}f_2^{\prime}(f^{j}(x))\equiv \alpha (x)^{p-i}\pmod{p^{2}},
$$
\noindent where
$$
\alpha (x)= \prod \limits ^{p^{n-1}-1}_{j=0}f_2^{\prime}(f^{j}(x)).
$$
\leftline{In view of (3) we now conclude that }
$$
f^{p^{n}}(x)\equiv x+p^{n-1}\xi (x)\Biggl (1+\sum\limits_{i=1}^{p-1}\alpha (x)^{i}\Biggr )\pmod{p^{n+1}}.
\eqnum{4}$$
Again, as $f_2^{\prime}$  modulo $p^{2}$ is periodic with period 
$p^{N_{2}(f)}$ and $f$ is transitive modulo $p^{n-1}$ for $n-1\ge N_{2}(f)$, then $\alpha (x)$ modulo $p^{2}$
does not depend on $x$. Moreover, we assert that $\alpha (x)\equiv 1\pmod{p}$.
\par
In fact, during the proof of 3.15 we have already established that if
$k\ge N_{1}(f)$ and if $f$ is a transitive modulo $p^{k}$ and uniformly differentiable 
modulo $p$ function with
integer-valued derivative modulo $p$, then
$$
\prod\limits^{p^{N_1 (f)}-1}_{j=0} f_1^{\prime}(f^{j}(x))\equiv 1\pmod{p} \eqnum{5}
$$
\noindent for all $x\in {\Bbb Z}_{p}$ (see the proof of (6) in 3.15). 
The definition of a
derivative modulo $p^{2}$ implies that $f_2^{\prime}(x)\equiv f_1^{\prime}(x)\pmod{p}$; consequently,
$$
\alpha (x)\equiv 1+p\beta \pmod{p^{2}}
\eqnum{6}$$
\noindent for some $\beta \in {\Bbb N}_{0}$. In view of (5) and (6), now (4) implies that 
$$
f^{p^{n}}(x)\equiv x+p^{n-1}\xi (x)\Biggl (p+p\beta \sum\limits_{i=1}^{p-1} i\Biggr )\pmod{p^{n+1}},
\eqnum{7}$$
\noindent and for $p\neq 2$ we conclude that
$$
f^{p^{n}}(x)\equiv x+p^{n}\xi (x)\pmod{p^{n+1}}.
$$
\noindent In view of 3.15 the latter proves the theorem for $p\neq 2$, since
$\xi (x)\not\equiv 0\pmod{p}$ (see the text which follows (2)).
\par
For the case $p=2$, the congruence (7) implies that
$$
f^{2^{n}}(x)\equiv x+2^{n}(1+\beta )\pmod{2^{n+1}}
\eqnum{8}$$
\noindent and to finish the proof it is sufficient to show that $\beta $ is even. 
\par For $n\ge N_{2}(f)+2$ the transitivity of $f$ modulo $2^{n}$ implies that $f$ is
transitive modulo $2^{N_{2}(f)+2}$, so in view of the definition of a derivative
modulo $p^2$  we have that
$$
f^{2^{N}}(x+2^{N}\xi )\equiv f^{2^{N}}(x)+2^{N}\xi \prod \limits_{j=0}^{2^N -1} f_2^{\prime}(f^{j}(x))\pmod{2^{N+2}}
\eqnum{9}$$
\noindent for $N=N_{2}(f)$, $\xi \in {\Bbb Z}_{2}$. As $f$ is transitive modulo $2^{N+2}$, then for arbitrary
$x\in \{0,1,\ldots  ,2^{N}-1\}$ and with $\xi$ running over $\{0,1,2,3\}$ the mapping
$$
\phi _{x}\colon\xi \mapsto\delta _{N}(f^{2^{N}}(x+2^{N}\xi ))+2\delta _{N+1}(f^{2^{N}}(x+2^{N}\xi ))
$$
\noindent is a cycle of length 4 on ${\Bbb Z}/4$. In view of (6), 
$$
\prod \limits_{j=0}^{2^N -1} f_2^{\prime}(f^{j}(x))\equiv 1+2\beta \pmod 4;
$$
\noindent so (9) implies that 
$$
\phi _{x}(\xi )\equiv c(x)+\xi (1+2\beta )\pmod 4,
\eqnum{10}$$
\noindent where $c(x)=\delta _{N}(f^{2^{N}}(x))+2\delta _{N+1}(f^{2^{N}}(x))$. 
But for each $x$ the mapping $\phi _{x}$ is
transitive modulo 4, so (10) in view of the  above mentioned transitivity criterion
for polynomials of degree 1 (see \cite {2 , Ch. 3, Theorem A})
implies that $\beta \equiv 0\pmod 2$.\qed
\enddemo
\remark 
{Note} {The analog of the theorem 3.14 generally does not hold for a function
which is uniformly differentiable modulo $p$.
Namely, for each  $n\in \Bbb N$ there exists a uniformly differentiable modulo
2 and compatible function 
$f\colon\Bbb Z_2\rightarrow \Bbb Z_2$ with $f^{\prime}_1=1$ everywhere
on $\Bbb Z_2$, 
$N_1(f)=1$,  
which is transitive modulo $2^k$ for $k=1,2,\ldots,n$, but which is not
transitive modulo 
$2^{k}$ for all  $k>n$. (By argument similar to applied below
one can construct a counterexample for  $p\ne 2$ as well.)}  
\par
Represent $x\in\Bbb Z_2$ in its canonic form 
$x=x_0+x_1\cdot 2+x_2\cdot 2^2+\ldots$, где $x_0,x_1,x_2\ldots \in \{0,1\}$.
Consider a function 
$$f(x)=\sum_{i=0}^{\infty}\phi_i(x_0,\ldots,x_i)\cdot 2^i,$$
where each $\phi_i(x_0,\ldots,x_i)$ is a Boolean polynomial, which is linear with respect to variable
$x_i$. In other words, $\phi_i(x_0,\ldots,x_i)=\psi_i(x_0,\ldots,x_{i-1})+x_i$
in the factor-ring
$\Bbb Z/2[x_0,\ldots,x_i]\big/(x_0^2-x_0,\ldots,x_i^2-x_i)$ of the
ring $\Bbb Z/2[x_0,\ldots,x_i]$ of all polynomials in variables $x_0,\ldots,x_i$
over $\Bbb Z/2$ with respect to the ideal, generated by $x_0^2-x_0,\ldots,x_i^2-x_i$ 
(we assume $\psi_0=1$). It is not difficult to see
that this function
$f$ is compatible (see 3.9 of \cite{11}). Direct calculations show
that for arbitrary
$s\in\Bbb N$ and $h\in\Bbb Z_2$ there holds a congruence 
$f(x+2^sh)\equiv f(x)+2^sh\pmod{2^{s+1}}$, i.e., that the function $f$ 
is uniformly differentiable 
modulo  2, and $f^{\prime}_1=1$ everywhere on $\Bbb Z_2$, with $N_1(f)=1$.
\par
Further, in the theory of Boolean functions there are well known
sufficient and necessary conditions for transitivity modulo
$2^n$ of the function $f$ of the considerd kind:
namely, it is transitive modulo  $2^n$ iff 
$\phi_i(x_0,\ldots,x_i)=\psi_i(x_0,\ldots,x_{i-1})+x_i$ for $i=1,2,\ldots,n-1$,
where each Boolean polynomial $\psi_i(x_0,\ldots,x_{i-1})$ for $i=1,2,\ldots,n-1$ 
is of odd weight
(that is, the number of all Boolean  vectors, satisfying it, is odd) 
and $\psi_0=1$. (This result, which is known as transitivity modulo $2^n$
criterion
for triangle transformations, belongs to mathematical folklore, so it
is difficult to refer the originating paper, yet a proof can be found in,
e.g., \cite{11}, see 4.8 there).
\par
Now choosing for a given  $n\in\Bbb N$ a function $f$ so that
$\psi_0=1$, 
with Boolean polynomials $\psi_i(x_0,\ldots,x_{i-1})$  
of odd weight for $i=1,2,\ldots,n-1$, and with Boolean polynomial 
$\psi_n(x_0,\ldots,x_{n-1})$ of even weight,
we obtain a function, which is transitive modulo 
$2^k$ for $k=1,2,\ldots,n$, but which is not transitive modulo
$2^{n+1}$. Then it is not transitive each modulo $2^k$ with $k>n$,
since, in view of compatibility of $f$, transitivity of $f$ modulo $2^{k+1}$
implies its transitivity modulo  $2^k$.
\endremark
\proclaim
{ 3.16 Corollary} { Let $A=\langle{\Bbb Z}_{p};\Omega \rangle$ be 
a universal algebra of finite signature
$\Omega $, and let all operations of $\Omega $ be uniformly differentiable
modulo $p^{2}$ functions with integer-valued derivatives modulo $p^2$. 
Then there exists
a positive rational integer $k(A)$ such that a polynomial $f(x)\in A[x]$ is
asymptotically ergodic if and only if it is transitive modulo $p^{k(A)}$.}
\endproclaim
\demo 
{ Proof} The proof of this corollary is similar to one of 3.10, $3^{\circ }$ 
and so is omitted.
We can take $k(A)=\max\{N_{2}(\omega ):\omega \in \Omega \}+\epsilon $, where $\epsilon =1$ 
if $p$ is odd, otherwise $\epsilon =2$.\qed
\enddemo
\head{4.} Hensel lift starting points.\endhead

The results of previous section show that for a class 
$\Cal D_1$ 
(respectively, $\Cal D_2$) of all uniformly differentiable modulo $p$
(respectively, modulo $p^2$) functions, which have integer-valued derivatives
modulo $p$
(respectively, modulo $p^2$), 
there exists a function $\zeta\colon \Cal D_1\rightarrow \Bbb N$ 
(respectively, $\eta\colon \Cal D_2\rightarrow \Bbb N$), such that a function
$f\in\Cal D_1$
(respectively, $f\in\Cal D_2$) is asymptotically measure-preserving (or
is ergodic)
iff it is bijective (respectively, transitive) modulo  $p^{\zeta(f)}$ 
(respectively, modulo $p^{\eta(f)}$). Theorems 3.9 and 3.14 give corresponding
estimates for  $\zeta(f)$ and
$\eta(f)$. 
\par
These estimates are sharp, i.e., there exist a compatible function 
$f\in\Cal D_1$ (respectively, $f\in\Cal D_2$) such that $f$ is bijective
(respectively, transitive) modulo 
$p^{N_1(f)}$ (respspectively, modulo $p^{N_2(f)}$ for $p\ne 2$, or modulo $2^{N_2(f)+1}$
for $p=2$), but $f$ is not measure-preserving (respectively, is not ergodic). 
For instance, a polynomial $f(x)=1+x^p$
is bijective modulo $p$, $N_1(f)=1$, but in force of $3.10, 1^{\circ}$ 
the polynomial $f$ is not bijective modulo  $p^2$, since
$f^{\prime}(z)\equiv 0\pmod p$ for all $z\in\Bbb Z_p$. 
\par
A corresponding example for theorem 3.14 in case $p\ne 2$ gives a function
$f(x)=(x+1)\odot_p 1$,
where $\odot_p$ is digitwise multiplication modulo $p$ of $p$-adic integers:
$\delta_i(x\odot_p y)\equiv\delta_i(x)\delta_i(y)\pmod p$ 
for all $i\in\Bbb N_0$. The function $f$ is uniformly differentiable,
its derivative is $0$ everywhere on $\Bbb Z_p$, and $N_2(f)=1$; at the
same time $f$ is transitive 
modulo $p$, but it is not even bijective (hence, is not transitive) modulo $p^2$.
\par
Nevertheless, boundaries for 
$\zeta(f)$ and
$\eta(f)$, 
which give, respectively, theorems 
3.9 and 
3.14, might differ significantly from the ones for various proper subclasses
of 
$\Cal D_1$ and of
$\Cal D_2$. For instance, for a function $f(x)=(ax+b)\XOR c$, 
with $a,b,c\in\Bbb N$, 
theorem 3.14 states that  $f(x)$ is asymptotically ergodic iff 
it is transitive modulo  $2^{\lfloor\log_2 c\rfloor+2}$, since this function
is uniformly differentiable and has a derivative which is  $a$ everywhere
on $\Bbb Z_2$, and
$N_2(f)=\lfloor\log_2 c\rfloor$.
Yet direct application of the above mentioned criteria of transitivity
modulo $2^n$
for triangle transformations and for polynomials of degree 1 over 
 $\Bbb Z$ immediately implies that $f$ is ergodic iff it is transitive
 modulo 4. So the problem of sharpening estimates of 
 $\zeta(f)$ and
$\eta(f)$ for various important from a certain view classes, which are narrower then $\Cal D_1$ and 
$\Cal D_2$, could be of interest.
\par 
In this section we study a class 
${\Cal A}$ of all compatible functions $f\colon{\Bbb Z}_{p}\rightarrow {\Bbb Z}_{p}$ 
such that, loosely speaking, coefficients of their interpolation series tends to 0 as fast
as
$i!$, or faster (recall that $\lim\limits_{i\to\infty}^p {i!}=0$).
More accurate, a function $f$, represented by interpolation series  $(\diamondsuit)$ 
(see section 2)
with $p$-adic integer coefficients $a_{i}$, belongs to ${\Cal A}$ iff it
is compatible, and a sequence
$\{\| {\frac {a{ } _{i}} {i!}}\| _{p}:i=0,1,2,\ldots\}$ is bounded, i.e.,
$\| {\frac {a{ } _{i}} {i!}}\| _{p}\leqslant     p^{\rho (f)}$ for some
${\rho (f)}\in {\Bbb N}_0$. Recall that according to the theorem 2.1, a
function $f$ 
represented by ($\diamondsuit$) is compatible iff
$\|a_i\|_p\le   p^{-\lfloor\log_p i\rfloor}$ for all 
$i\in\Bbb N$. 
\par
Class ${\Cal A}$ is rather wide: it contains all integer-valued compatible
analytic on $\Bbb Z_p$ functions, in particular, compatible functions which could
be defined
by integer-valued polynomials over  $\Bbb Q_p$. 
It is known (see \cite{3, Ch. 4, Theorem 4, p. 224}), that a function 
$f$ of the form ($\diamondsuit$) is analytic on $\Bbb Z_p$ iff 
$ \lim\limits_{i\to \infty}^p {\frac {a_{i}} {i!}}=0$.
\par 
So for the rest of this section we assume that $f\in {\Cal A}$. Put
$$
\lambda (f)=\min\biggl\{k\in \Bbb N: 2 {\frac {p^{k}-1} {p-1}} - k >\rho (f)\biggr\}.
$$
The following theorem is true.
\proclaim
{ 4.1 Theorem} { Let $f\in {\Cal A}$ and $p$ is an odd prime. The function $f$ is ergodic if
and only if it is transitive modulo $p^{\lambda (f)+1}$ {\rom(}if $p\neq 3${\rom)} or modulo $3^{\lambda (f)+2}$
{\rom(}if $p=3${\rom)}.}
\endproclaim
\par Since $f$ is compatible, then in view of 2.1 it can be represented  in the following
form:
$$
f(x)=b_{0}+\sum^{\infty }_{i=1}b_{i}p^{\left\lfloor{\log_{p}i}\right\rfloor}{ {x}\choose {i}},
$$
\noindent where $b_{j}\in {\Bbb Z}_{p}$ for $j=0,1,2,\ldots .$ Everywhere during the proof we assume that $f$
is represented in this form. Further $\lambda (f)$ is denoted as $\lambda $ and $p$ is assumed
to be an odd
prime. We will need some additional technical results.
\proclaim
{ 4.2 Lemma} { Under the assumptions of  theorem 4.1  the  following is true: 
$$
\displaylines{\hfill b_{i}\equiv 0\pmod p\hbox{, for }i\ge 2p^{\lambda };\hfill\cr
\hfill b_{i}\equiv 0\pmod{p^{2}}\hbox{, for }i\ge 3p^{\lambda }.\hfill\cr}
$$}
\endproclaim 
\demo
{ Proof of the lemma 4.2} If $b_{i}=0$, then the assertion of the lemma is trivial. Suppose that
$b_{i}\neq 0$. Represent $f$ as
$$
f(x)=b_{0}+\sum^{\infty }_{i=1}{\frac {1} {i!}}b_{i}p^{\left\lfloor{\log_{p}i}\right\rfloor} (x)_{i},
$$
\noindent where, we recall, $(x)_{i}=x(x-1)\cdots (x-i+1)$ (with $(x)_0=1$)
is $i$th descending factorial power of $x$. As $f\in {\Cal A}$, i.e., 
$$
\left\Vert b_{i}p^{\left\lfloor{\log_{p}i}\right\rfloor} \right\Vert_{p} \le p^{\rho (f)} \Vert i! \Vert_p
$$
\noindent then
\par
$$
\ord_{p}\,b_{i}\ge \ord_{p}\,i!-\left\lfloor{\log_{p}i}\right\rfloor-\rho (f),
\eqnum{1}$$
\noindent for all $i=1,2,\ldots  $. We recall that $\log_{p}\| a\| _{p}=-\ord_{p}\,a$, for $a\in {\Bbb Z}_{p}$.
Thus, the maximal $p$-prime factor of $a$ is exactly $p^{\ord_{p}\,a}$.
\par
In fact, the function $\kappa (i)=\ord_{p}\,i!-\left\lfloor{\log_{p}i}\right\rfloor$ is nondecreasing.
To prove this, note that, obviously, $\ord_{p}\,i!\ge \ord_{p}\,(i-1)!$. If
$\left\lfloor{\log_{p}i}\right\rfloor= \left\lfloor{\log_{p}(i-1)}\right\rfloor$ then $\kappa (i-1)\le \kappa (i)$.
\par
Assume $\left\lfloor{\log_{p}j}\right\rfloor>\left\lfloor{\log_{p}(j-1)}\right\rfloor$ for some positive rational integer $j$.
Evidently, $\left\lfloor{\log_{p}j}\right\rfloor+1$ is the number of significant
digits in the $p$-base expansion of $j$. 
Hence the case under consideration takes place exactly if and only
if $j-1=(p-1)+(p-1)p+\cdots  +(p-1)p^{n}=p^{n+1}-1$ for some $n\in {\Bbb N}_0$. But then
$\ord_{p}\,j!=\ord_{p}\,(j-1)!+n, \left\lfloor{\log_{p}(j-1)}\right\rfloor=n,\left\lfloor{\log_{p}j}\right\rfloor=n+1$, and so
$\kappa (j)>\kappa (j-1)$.
\par
Now it is sufficient to prove only that  $\kappa (2p^{\lambda })-\rho (f)\ge 1$  and $\kappa (3p^{\lambda })-
\rho (f)\ge 2$. We recall that $\ord_{p}\,i!={\frac {1} {p-1}}(i-\wt_{p}\,i)$,  where  $\wt_{p}\,i$  is the sum of
all digits in a $p$-base  expansion  of  $i$ (i.e.,  if
$i=i_{0}+i_{1}p+\cdots  +i_{s}p^{s}$, where $i_{0},\ldots  ,i_{s}\in \{0,1,\ldots  ,p-1\}$, then $\wt_{p}\,i=i_{0}+\cdots  +i_{s}$, see
e.g., \cite{6 , ch.1, section 2, exercise 13}).
\par
As $p\neq 2$, then $\kappa (2p^{\lambda })-\rho (f)={\frac {1} {p-1}}(2p^{\lambda }-2)-\lambda -\rho (f)\ge 1$ according to the
definition of $\lambda =\lambda (f)$. Hence, if $p\neq 3$, then
$$
\kappa (3p^{\lambda })-\rho (f)={\frac {1} {p-1}}(3p^{\lambda }-3)-\lambda -\rho (f)=\kappa (2p^{\lambda })+{\frac {1} {p-1}}(p^{\lambda }-1)-\rho (f)\ge 2.
$$
\leftline{So if $p\neq 3$ the lemma is proved. }
\par 
Finally, let $p=3$. Then
$$
\kappa (3p^{\lambda })-\rho (f)=\kappa (3^{\lambda +1})-\rho (f)={\frac {1} {2}}(3^{\lambda +1}-1)-\lambda -1-\rho (f)\ge 2,
$$
\leftline{otherwise in view of the inequality }
$$
3^{\lambda }-1-\lambda > \rho (f),
$$
\leftline{(which follows directly from the definition of $\lambda =\lambda (f))$ we get }
$$
{\frac {1} {2}}(3^{\lambda +1}-1)-\lambda -1 - 3^{\lambda }+1+\lambda <1,
$$
\noindent i.e., $3^{\lambda }-1<2$, and so $\lambda <1$, a contradiction. The lemma 
4.2 is proved.\qed
\enddemo
\proclaim
{ 4.3 Corollary} { Under assumptions of theorem 4.2, for $i\in {\Bbb N}$ the following is
true:
$$
{\frac {\Delta ^{i}f(x)} {i}}\equiv \cases 0\pmod{p^2},&\text{\rom{if} } i\ge 2p^{\lambda }+1;\cr 0\pmod p,&\text{\rom{if} }i\ge p^{\lambda }+1.\cr\endcases 
$$}
\endproclaim
\demo 
{Proof of the corollary  4.3}  As  $\Delta ^{j}{ {x}\choose {i}}={ {x}\choose {i-j}}$ if $i\ge j$  and  $\Delta ^{j}{ {x}\choose {i}}=0$ if $i<j$,
then 
\par
$$
{\frac {\Delta ^{i}f(x)} {i}}={\frac {1} {\hat\imath}} \sum^{\infty }_{j=i}b_{j}p^{\left\lfloor{\log_{p}j}\right\rfloor {-\ord_{p}\,j}}{ {x}\choose {j-i}},
$$
\noindent where $\hat\imath =ip^{-\ord_{p}\,i}\in {\Bbb Z}_{p}, \ord_{p}\,\hat\imath =0$. Now the result is obvious in view of lemma
4.2.\qed
\enddemo
\proclaim 
{ 4.4 Proposition} { Under assumptions of theorem 4.1 the function $f$
is uniformly differentiable modulo $p^{2}$, has integer-valued derivative
modulo $p^{2}$,
$N_{2}(f)\le \lambda (f)+1$. Moreover,
$$
f^\prime_2 (x)\equiv \sum^{2p{ } ^{\lambda }}_{i=1}(-1)^{i-1}{\frac {\Delta ^{i}f(x)} {i}}\pmod{p^{2}}.
$$}
\endproclaim
\demo 
{Proof of the proposition 4.4} To prove the first assertion of the
proposition we will demonstrate that there exists a function $f^\prime_2 \colon{\Bbb Z}_{p}\rightarrow {\Bbb Z}_{p}$ such
that for all $x,h\in {\Bbb Z}_{p}$ and $m\ge \lambda (f)+1$ the following congruence holds:
$$
f(x+p^{m}h)\equiv f(x)+p^{m}hf^\prime_2 (x)\pmod{p^{m+2}}.
\eqnum{1}$$
\noindent In view of the compatibility of $f$, it is sufficient to prove the
congruence (1) only for $h\in \{1,2,\ldots  ,p^{2}-1\}$ (for $h=0$ the congruence is
trivial). Applying Newton formula
$$
f(x+n)=\sum^{n}_{i=0}{ {n}\choose {i}}\Delta ^{i}f(x)
$$
\leftline{for $n=p^mh$, we have }
$$
f(x+p^{m}h)=f(x)+p^{m}h\phi _{m}(x,h),
\eqnum{2}$$
\leftline{where }
$$
\phi _{m}(x,h)=\sum^{p^{m}h}_{i=1}{ {p^{m}h-1}\choose {i-1}}{\frac {\Delta ^{i}f(x)} {i}} .
\eqnum{3}$$
\leftline{Hence in view of 4.3 for $m\ge \lambda +1$ we obtain: }
$$
\phi _{m}(x,h)\equiv \sum^{2p{ } ^{\lambda }}_{i=1}{ {p^{m}h-1}\choose {i-1}}{\frac {\Delta ^{i}f(x)} {i}}\pmod{p^{2}}.
\eqnum{4}$$
\noindent Further, for $i=1,2,\ldots  ,2p^{\lambda }$ the following obviuos
equalities hold:
$$
{ {p^{m}h-1}\choose {i-1}}=\prod^{i-2}_{k=0}{\frac {p^{m}h-(k+1)} {k+1}}=\prod^{i-1}_{j=1}\Biggl ( {\frac {h} {\hat\jmath}}p^{m-\ord_{p}\,j}-1\Biggr ).
\eqnum{5}$$
\noindent Here $\hat\jmath =jp^{-\ord_{p}\,j}$ is the unit of ${\Bbb Z}_{p}$, 
i.e., $\hat\jmath$ has multiplicative inverse
${\frac {1} {\hat\jmath}}$ in ${\Bbb Z}_{p}$; hence, each factor of the product in the right hand part of (5)
is $p$-adic integer.
\par
If $i\le p^{\lambda }$ then $m-\ord_{p}\,j\ge 2$ for all $j=1,2,\ldots  ,i-1$; so (5) implies that
$$
{ {p^{m}h-1}\choose {i-1}}\equiv (-1)^{i-1}\pmod{p^{2}}.
\eqnum{6}$$
If $p^{\lambda }+1\le i\le 2p^{\lambda }$ and $j\in \{1,2,\ldots  ,i-1\}$ then $m-\ord_{p}\,j=1$ only in the case
when simultaneously $j=p^{\lambda }$ and $m=\lambda +1$ hold; otherwise $m-\ord_{p}\,j\ge 2$.
Yet if $m-\ord_{p}\,j=1$ then
$$
{\frac {\Delta ^{i}f(x)} {i}}\equiv 0\pmod{p}
$$
\noindent (see 4.3); hence in both cases we have that
$$
\Biggl ( {\frac {h} {\hat\jmath}}p^{m-\ord_{p}\,j}-1\Biggr ){\frac {\Delta ^{i}f(x)} {i}}\equiv - {\frac {\Delta ^{i}f(x)} {i}}\pmod{p^{2}}.
$$
\leftline{So in view of (5) we conclude that }
$$
{ {p^{m}h-1}\choose {i-1}}{\frac {\Delta ^{i}f(x)} {i}}\equiv (-1)^{i-1} {\frac {\Delta ^{i}f(x)} {i}}\pmod{p^{2}}.
\eqnum{7}$$
\noindent for all $i=1,2,\ldots, 2p^\lambda$. Now (4), (6), (7) together imply that
$$
\phi _{m}(x,h)\equiv \sum^{2p{ } ^{\lambda }}_{i=1}(-1)^{i-1}{\frac {\Delta ^{i}f(x)} {i}}\pmod{p^{2}}
$$
\noindent and in view of (2), (3), (4) this completes the proof 
of proposition 4.4.\qed
\enddemo
\proclaim
{4.5. Lemma} { Under assumptions of theorem 4.1, there exists a function 
$\theta \colon{\Bbb Z}_{p}\rightarrow {\Bbb Z}_{p}$ such that
for arbitrary $x,h\in {\Bbb Z}_{p}$ the 
following congruence holds:
$$
f(x+p^{\lambda }h)\equiv f(x)+p^{\lambda }hf_2^{\prime}(x)+p^{\lambda +1}h^{2}\theta (x)\pmod{p^{\lambda +2}}.
$$
\noindent The function $\theta $ satisfies the following condition: for arbitrary $a,b\in {\Bbb Z}_{p}$ the congruence
$a\equiv b\pmod{p^{\lambda}}$ implies $\theta (a)\equiv \theta (b)\pmod{p}$. 
Moreover, one may put
$$
\theta (x)=\sum^{p-1}_{j=2}(-1)^{j}\sum^{j-1}_{i=1}{\frac {1} {i}} {\frac {\Delta ^{jp^{\lambda -1}}f(x)} {jp^{\lambda -1}}}+\sum^{p-1}_{k=1}(-1)^{k-1}{\frac {\Delta ^{kp^{\lambda -1}+p^{\lambda }}f(x)} {kp^{\lambda }}} + {\frac {\Delta ^{2p^{\lambda }}f(x)} {2p^{\lambda +1}}}.
$$}
\endproclaim
\demo 
{Proof of the lemma 4.5} Firstly we prove that the function $\theta $ defined by
the latter equality is integer-valued 
on ${\Bbb Z}_{p}$.  
Since $f$ is compatible, each fraction $\frac{\Delta^s f(x)}{s}$ for $s=1,2,3,\ldots$,
is $p$-adic integer (see 3.1 of \cite {11}). So it is sufficient to prove only that for all
 $k\in \{1,2,\ldots  ,p-1\}$ both functions $\alpha (x)$ and $\beta _{k}(x)$ (defined below)
are integer-valued 
on ${\Bbb Z}_{p}$. By definition,
$$
\alpha (x)={\frac {\Delta ^{2p^{\lambda }}f(x)} {2p^{\lambda +1}}};\quad \beta _{k}(x)={\frac {\Delta ^{kp^{\lambda -1}+p^{\lambda }}f(x)} {kp^{\lambda }}}.
$$
\par
 Since 
$$
\Delta ^{i}f(x)=\sum^{\infty }_{j=i}b_{j}p^{\left\lfloor{\log_{p}j}\right\rfloor}{ {x}\choose {j-i}}
\eqnum{1}$$
\leftline{for $i=1,2,3,\ldots  $ and }
$$
b_{j}p^{\left\lfloor{\log_{p}j}\right\rfloor}\equiv 0\pmod{ p^{\lambda +1}}
$$
\noindent for all integer rationals $j\ge 2p^{\lambda }$ (see 4.2), 
then $\alpha (x)\in {\Bbb Z}_{p}$. If $j\ge kp^{\lambda -1}+p^{\lambda }$
then $\left\lfloor{\log_{p}j}\right\rfloor\ge \lambda $; hence (1) implies 
that $\beta _{k}(x)\in {\Bbb Z}_{p}$.
\par
Now we prove that for all
$a,b\in {\Bbb Z}_{p}$ the congruence $a\equiv b\pmod{p^{\lambda}}$ 
implies $\theta (a)\equiv \theta (b)\pmod{p}$. 
In view of (1) and 4.2 the following congruence holds:
$$
\alpha (x)\equiv {\frac {1} {2}}\sum^{3p^{\lambda }-1}_{j=2p{ } ^{\lambda }}{\frac {1} {p}} b_{j}{ {x}\choose {j-2p{ } ^{\lambda }}}\pmod{p}.
\eqnum{2}$$
\par We recall a statement of the well known Lucas theorem
 (for a proof see e.g 
\cite 4): if $a=\sum^{\infty }_{i=0}a_{i}p^{i}$ and
$b=\sum^{N}_{i=0}b_{i}p^{i}$ are, respectively, canonic representations
of $p$-adic
integer $a$ and of nonnegative integer rational $b$
(i.e., $a_{i},b_{i}\in \{0,1,\ldots  ,p-1\}$ for $i=0,1,2,\ldots$), then 
$$
{a\choose b}\equiv {{a_0}\choose {b_0}}{{a_1}\choose {b_1}}\cdots {{a_N}\choose {b_N}}\pmod{p}.$$ 

\par So, if $a\equiv b\pmod{p^{\lambda}}$, then Lucas theorem implies that for all $j=2p^{\lambda },2p^{\lambda }+1,\ldots  ,3p^{\lambda }-1$ the
following congruence holds:
$$
{ {a}\choose {j-2p{ } ^{\lambda }}}\equiv { {b}\choose {j-2p{ } ^{\lambda }}}\pmod{p}.
$$
\leftline{Thus, (2) implies that }
$$
\alpha (a)\equiv \alpha (b)\pmod{p}.
\eqnum{3}$$
\leftline{Further, combining (1) and 4.2 we obtain that}
$$
\beta _{k}(x)\equiv {\frac {1} {k}}\sum^{2p^{\lambda }-1}_{j=kp^{\lambda -1}+p{ } ^{\lambda }} b_{j}{ {x}\choose {j-kp^{\lambda -1}-p{ } ^{\lambda }}}\pmod{p}
$$
\noindent for all $k=1,2,\ldots  ,p-1$. Now applying Lucas theorem once
again, we conclude that
$$\beta _{k}(a)\equiv \beta _{k}(b)\pmod{p}\eqnum {4}$$
\noindent for $a\equiv b\pmod{p^{\lambda}}$.
\par
Lastly, assuming
$$
\gamma _{k}(x)={\frac {\Delta ^{kp^{\lambda -1}_{}}f(x)} {kp^{\lambda -1}}},
$$
\noindent in view of (1) we conclude that for $k=1,2,\ldots  ,p-1$ the following congruence holds:
$$
\gamma _{k}(x)\equiv {\frac {1} {k}}\sum^{p^{\lambda }-1}_{j=kp^{\lambda -1}} b_{j}{ {x}\choose {j-kp^{\lambda -1}}}\pmod{p}.
$$
\leftline{Again, applying Lucas theorem, we conclude that }
$$
\gamma _{k}(a)\equiv \gamma _{k}(b)\pmod{p}
\eqnum{5}$$
\noindent for $a\equiv b\pmod{p^{\lambda}}$. Hence in view of (3) -- (5) the congruence $a\equiv b\pmod{p^{\lambda}}$
implies the congruence $\theta (a)\equiv \theta (b)\pmod{p}$.
\par
Now we prove the rest of the lemma. As $f$ is compatible, during the
proof we may assume that $h\in {\Bbb N}$ (case $h=0$ is trivial). According to 4.4
(see (2)--(5) there) the following is true:
$$
f(x+p^{\lambda }h)\equiv f(x)+p^{\lambda }h\phi (x,h)\pmod{p^{\lambda +2}},
\eqnum{6}$$
\leftline{where }
$$
\phi (x,h)\equiv \sum^{2p^{\lambda }}_{i=1}{ {p^{\lambda }h-1}\choose {i-1}}{\frac {\Delta ^{i}f(x)} {i}}\pmod{ p^2}.
\eqnum{7}$$
\leftline{and, besides, }
$$
{ {p^{\lambda }h-1}\choose {i-1}}=\prod^{i-1}_{j=1}({\frac {h} {\hat\jmath}}p^{\lambda -\ord_{p}\,j}-1)
\eqnum{8}$$
\leftline{for $i=1,2,\ldots  ,2p^{\lambda }$.}
As $f$ is compatible, then, according to 3.40 of \cite{11}, 
$$\frac {\Delta ^{i}f(x)} {i}\equiv 0\pmod p$$
in all cases with the exception of, possibly, a case when
$i$ is of the form $i=tp^s$ 
for suitable $t\in\{1,2,\ldots p-1\}$ and $s\in\Bbb N_0$.
Thus, if $i\le p^{\lambda -1}$, 
as well as if simultaneously $p^{\lambda -1}<i<p^{\lambda }$ and $p^{\lambda -1}$
is not a factor of $i$, the equality (8) implies:
$$
{ {p^{\lambda }h-1}\choose {i-1}}{\frac {\Delta ^{i}f(x)} {i}}\equiv (-1)^{i-1}{\frac {\Delta ^{i}f(x)} {i}}\pmod{ p^2}.
\eqnum{9}$$
\leftline{Let $i=kp^{\lambda -1}$ for $k\in \{2,3,\ldots  ,p-1\}$. Then (8) implies: }
$$
{ {p^{\lambda }h-1}\choose {i-1}}\equiv (-1)^{kp^{\lambda -1}-1}+(-1)^{k}ph\sum^{k-1}_{j=1}{\frac {1} {j}}\pmod{ p^2}.
\eqnum{10}$$
\noindent Further, if $p^{\lambda }\le i\le 2p^{\lambda }$ and $\ord_{p}\,i\neq \lambda ,\lambda -1$ then (1)
(together with congruence following it) imply that 
$$
{\frac {\Delta ^{i}f(x)} {i}}\equiv 0\pmod{ p^2}.
\eqnum{11}$$
\par
Now we have to study the only two remaining cases: $i=\nu p^{\lambda }$ for $\nu \in \{1,2\}$ and $i=kp^{\lambda -1}+p^{\lambda }$ for
$k\in \{1,2,\ldots  ,p-1\}$. The latter one in view of 4.3 and (8) implies that
$$
{ {p^{\lambda }h-1}\choose {i-1}}{\frac {\Delta ^{i}f(x)} {i}}\equiv (-1)^{i-1}{\frac {\Delta ^{i}f(x)} {i}}+(-1)^{k-1}h{\frac {\Delta ^{i}f(x)} {i}}\pmod{ p^2}.
\eqnum{12}$$
\noindent Further, for $k=1,2,\ldots  ,p-1$ the following trivial equality holds
in ${\Bbb Q}_{p}$:
$$
\biggl (1+{\frac {p} {k}}\biggr ) {\frac {\Delta ^{kp^{\lambda -1}+p^{\lambda }}f(x)} {kp^{\lambda -1}+p^{\lambda }}} = {\frac {\Delta ^{kp^{\lambda -1}+p^{\lambda }}f(x)} {kp^{\lambda -1}}}
\eqnum{13}$$
\noindent From here  in view of 4.3 we conclude that
$$
{\frac {\Delta ^{kp^{\lambda -1}+p^{\lambda }}f(x)} {kp^{\lambda -1}+p^{\lambda }}}\equiv 0\pmod{p}
$$
\noindent and since ${\frac {p} {k}}\in {\Bbb Z}_{p}$ and $\ord_p{\frac{p}{k}}=1$, the equality (13) implies that 
$$
{\frac {\Delta ^{kp^{\lambda -1}+p^{\lambda }}f(x)} {kp^{\lambda -1}+p^{\lambda }}}\equiv {\frac {\Delta ^{kp^{\lambda -1}+p^{\lambda }}f(x)} {kp^{\lambda -1}_{}}} \pmod{ p^2}.
$$
\leftline{Hence, applying (12) for $i=kp^{\lambda -1}+p^{\lambda }$ , we have that}
$$
\displaylines {{ {p^{\lambda }h-1}\choose {kp^{\lambda -1}+p^{\lambda }-1}}{\frac {\Delta ^{kp^{\lambda -1}+p^{\lambda }}f(x)} {kp^{\lambda -1}+p^{\lambda }}}\hfill\cr
\hfill{}=(-1)^{kp^{\lambda -1}+p^{\lambda }-1} {\frac {\Delta ^{kp^{\lambda -1}+p^{\lambda }}f(x)} {kp^{\lambda -1}+p^{\lambda }}}+(-1)^{k-1}ph\beta _{k}(x)\pmod{ p^2}. \qquad (14)\cr}$$
\par In case $i=p^{\lambda }$, the equality (8) implies that
$$
{ {p^{\lambda }h-1}\choose {p^{\lambda }-1}}\equiv (-1)^{p^{\lambda}-1}-ph\sum^{p-1}_{j=1}{\frac {1} {j}}\equiv (-1)^{p^{\lambda }-1}\pmod{ p^2},
\eqnum{15}$$
\noindent since for $p\neq 2$ the following congruences hold in ${\Bbb Q}_{p}$: $\sum^{p-1}_{j=1}{\frac {1} {j}}\equiv \sum^{p-1}_{j=1}j\equiv 0\pmod{p}$.
\par
Finally, for $i=2p^{\lambda }$, applying (8) and 4.3, we conclude that
$$\displaylines{{ {p^{\lambda }h-1}\choose {2p^{\lambda }-1}}{\frac {\Delta ^{2p^{\lambda }}f(x)} {2p^{\lambda }}}\equiv (-1)^{2p^{\lambda }-1} {\frac {\Delta ^{2p^{\lambda }}f(x)} {2p^{\lambda }}}+h{\frac {\Delta ^{2p^{\lambda }}f(x)} {2p^{\lambda }}}
\hfill\cr
\hfill{}\equiv (-1)^{2p^{\lambda }-1} {\frac {\Delta ^{2p^{\lambda }}f(x)} {2p^{\lambda }}}+hp\alpha (x)\pmod{ p^2},\quad (16)\cr}$$ 
\noindent where $\alpha (x)\in {\Bbb Z}_{p}$, as it was shown above.
\par 
Now by the combination of (6), (7), (9), (11), (14), (15), (16) with
4.4 we finish the proof of the lemma 4.5.\qed
\enddemo
\proclaim
{4.6 Lemma} { Under assumptions of theorem 4.1, for all $x,h\in {\Bbb Z}_{p}$ the
following congruence holds:
$$
f_2^\prime (x+p^{\lambda }h)\equiv f_2^\prime (x)+2ph\theta (x)\pmod{p^2}.
$$
\noindent Here $\theta $ is the function defined in 4.5. }
\endproclaim
\demo 
{Proof of the lemma 4.6} In view of 4.4 the following is true: 
$$
f^{\prime }_2(x+p^{\lambda }h)\equiv \sum^{2p{ } ^{\lambda }}_{i=1}(-1)^{i-1}{\frac {\Delta ^{i}f(x+p^{\lambda}h)} {i}}\pmod{p^2}.
\eqnum{1}$$
\noindent For $i=1,2,\ldots  ,2p^{\lambda }$ the previous lemma implies that 
$$
\displaylines{{\frac {\Delta ^{i}f(x+p^{\lambda }h)} {i}}\equiv {\frac {\Delta ^{i}f(x)} {i}} +hp^{\lambda -\ord_{p}\,i}{\frac {\Delta ^{i}f^{\prime }_2(x)} {\hat\imath}}\hfill\cr
\hfill {}+h^{2}p^{\lambda +1-\ord_{p}\,i}{\frac {\Delta ^{i}\theta _{}(x)} {\hat\imath}}\pmod{p^2},
\quad (2)\cr}$$
\noindent where $\hat\imath =ip^{-\ord_{p}\,i}$ is a unit in ${\Bbb Z}_{p}$, i.e., it has a multiplicative inverse
${1\over \hat\imath}\in{\Bbb Z}_{p}$.
\par
The term of order 2 (with respect to $h$) in (2) may not vanish modulo $p^{2}$ 
only if $i\in \{p^{\lambda },2p^{\lambda }\}$. Yet, as 
$\Delta ^{j}{ {x}\choose {\nu}}={ {x}\choose {\nu-j}}$ 
for $\nu\ge j$  and  $\Delta ^{j}{ {x}\choose {\nu}}=0$ for $\nu<j$,
then for all $j\in {\Bbb N}$ we have
$$
\Delta ^{j}f(x)=\sum^{\infty }_{\nu =j}b_{\nu }p^{\left\lfloor{\log_{p}\nu }\right\rfloor}{ {x}\choose {\nu -j}}.
\eqnum{3}$$
\noindent Consequently, if $j\in \{p^{\lambda },2p^{\lambda }\}$, then 
$$
{\frac {\Delta ^{j+kp^{\lambda -1}_{}}f(x)} {kp^{\lambda -1}_{}}}\equiv 0\pmod{p}.
\eqnum{4}$$
\noindent for $k\in\{1,2,\ldots,p-1\}$. Further, for 
$j\in \{p^{\lambda },2p^{\lambda }\}$ the equality (3) in view of  4.3
implies that
$$\displaylines{\hfill{\frac {\Delta ^{j+kp^{\lambda -1}+p^{\lambda }}f(x)} {kp^{\lambda }_{}}}\equiv 0\pmod{p},\hfill (5)\cr
\hfill{\frac {\Delta ^{j+2p^{\lambda }}f(x)} {2p^{\lambda }}}\equiv 0\pmod{p}.\hfill (6)\cr}$$
\noindent Now, by the definition of $\theta $, combining together (4), (5), (6) we conclude
that ${\frac {\Delta ^{i}\theta _{}(x)} {\hat\imath}}\equiv 0\pmod{p}$ for $i\in \{p^{\lambda }, 2p^{\lambda }\}$, and thus
$$
h^{2}p^{\lambda +1-\ord_{p}\,i}{\frac {\Delta ^{i}\theta _{}(x)} {\hat\imath}}\equiv 0\pmod{p^2}
\eqnum{7}$$
\leftline{for all $i=1,2,\ldots  ,2p^{\lambda }$. }
\par
The term of order 1 in (2) may not vanish modulo $p^{2}$
only for $i\in \{1,2,\ldots  ,2p^{\lambda }\}$ such that $\ord_{p}\,i\ge \lambda -1$, i.e., for
$$
i\in \{p^{\lambda },2p^{\lambda },kp^{\lambda-1},kp^{\lambda -1}+p^{\lambda }:k=1,2,\ldots  ,p-1\}.
$$
\noindent Combining together 4.3, 4.4 and 3.4 of \cite{11} we already referred
(see argument which follows (8) in the proof of 4.5), we have 
$$
f^{\prime }_2(x)\equiv {\frac {\Delta ^{p^{\lambda }}f(x)} {p^{\lambda }}}+\sum^{\lambda -1}_{t=0}\sum^{p-1}_{\tau =1}(-1)^{\tau -1}{\frac {\Delta ^{\tau p^{t}_{}}f(x)} {\tau p^{t}}}\pmod{p},
\eqnum{8}$$
\leftline{and hence }
$$
\Delta ^{i}f^{\prime }_2 (x)\equiv {\frac {\Delta ^{i+p^{\lambda }}f(x)} {p^{\lambda }}}+\sum^{\lambda -1}_{t=0}\sum^{p-1}_{\tau =1}(-1)^{\tau -1}{\frac {\Delta ^{i+\tau p^{t}_{}}f(x)} {\tau p^{t}}}\pmod{p}.
\eqnum{9}$$
\noindent This for $i\in \{kp^{\lambda -1}+p^{\lambda }:k=1,2,\ldots  ,p-1\}$ 
in force of (3) and 4.2 implies that
$\Delta ^{i}f^{\prime }_2(x)\equiv 0\pmod{p}$, and consequently
$$
hp{\frac {\Delta ^{kp^{\lambda -1}+p^{\lambda }}f^{\prime }_2(x)} {k+p}}\equiv 0\pmod{p^2}
\eqnum{10}$$
\noindent for $k=1,2,\ldots  ,p-1$ (since multiplicative inverse ${1\over k+p}$
of $k+p$ is in $\Bbb Z_p$). 
\par 
If $i\in \{kp^{\lambda -1}:k=1,2,\ldots  ,p-1\}$ then in view of 4.2, (3) and (9) we have:
$$
\Delta ^{kp^{\lambda -1}}f^{\prime }_2(x)\equiv {\frac {\Delta ^{kp^{\lambda -1}+p^{\lambda }}f_{}(x)} {p^{\lambda }}} +\sum^{p-k-1}_{\tau =1}(-1)^{\tau -1}{\frac {\Delta ^{(\tau +k)p^{\lambda -1}_{}}f(x)} {\tau p^{\lambda -1}}} \pmod{p}.
\eqnum{11}$$
\noindent If $i=2p^{\lambda }$ then 4.4 implies that
$$
\Delta ^{2p^{\lambda }}f^{\prime }_2(x)\equiv \sum^{2p^{\lambda }}_{j=1}(-1)^{j-1}{\frac {\Delta ^{j+2p^{\lambda }_{}}f(x)} {j_{}}}\pmod{p^2}.
$$
\leftline{This in view of (3) and 4.2 implies that }
$$
\Delta ^{2p^{\lambda }}f^{\prime }_2(x)\equiv 0\pmod{p^2}.
\eqnum{12}$$
\par Now we consider a case $i=p^{\lambda }$. Proposition 4.4 implies that
$$
\Delta ^{p^{\lambda }}f^{\prime }_2(x)\equiv \sum^{1+p^{\lambda }}_{j=1}(-1)^{j-1}{\frac {\Delta ^{j+p^{\lambda }_{}}f(x)} {j_{}}}\pmod{p^2},
\eqnum{13}$$
\noindent since, combining together (3) and 4.2, 
for $j=p^{\lambda }+1,\ldots  ,2p^{\lambda }$ we  conclude
that
$$
{\frac {\Delta ^{j+p^{\lambda }_{}}f(x)} {j_{}}}\equiv 0\pmod{p^2}.
$$
Moreover, (3) implies that the latter congruence holds also for all $j\le p^{\lambda }-1$,
such that $j\neq kp^{\lambda -1}$, where $k=1,2,\ldots  ,p-1$. Thus, (13) implies that
$$
\Delta ^{p^{\lambda }}f^{\prime }_2(x)\equiv {\frac {\Delta ^{_{}2p^{\lambda }}f_{}(x)} {p^{\lambda }}} +\sum^{p-1}_{k=1}(-1)^{k-1} {\frac {\Delta ^{kp^{\lambda -1}+p^{\lambda }}f(x)} {kp^{\lambda -1}}}\pmod{p^2}.
\eqnum{14}$$
Now, substituting (7), (10), (11), (12), (14) to (2) and summarizing
up all the obtained congruences for $i$ ranging from 1 to $2p^{\lambda }$, in view of
(1) and 4.4 we conclude that
$$
\displaylines{f^{\prime }_2 (x+p^{\lambda }h)\equiv f^{\prime }_2(x)+hp\Biggl (\sum^{p-1}_{k=1}{\frac {(-1)^{k-1}} {k}}\sum^{p-k-1}_{\tau =1}(-1)^{\tau -1}{\frac {\Delta ^{(\tau +k)p^{\lambda -1}_{}}f(x)} {\tau p^{\lambda -1}}}\hfill\cr
\hfill {}+ \sum^{p-1}_{k=1}(-1)^{k-1}{\frac {\Delta ^{kp^{\lambda -1}+p^{\lambda }}f(x)} {kp^{\lambda }}}\Biggr ) \hfill\cr
\hfill {}+ h\sum^{p-1}_{k=1}(-1)^{k-1}{\frac {\Delta ^{kp^{\lambda -1}+p^{\lambda }}f(x)} {kp^{\lambda -1}}}
+h {\frac {\Delta ^{2p^{\lambda }_{}}f(x)} {p^{\lambda }}}(\bmod p^{2}).
\quad (15)\cr}$$
\noindent We recall that here and after all calculations are performed in the
field ${\Bbb Q}_{p}$, and by the above agreement the congruence $\xi \equiv 0\pmod{p^{k}}$ for $\xi \in {\Bbb Q}_{p}$
and positive integer rational $k$ means that $\| \xi \| _{p}=p^{-k}$ (hence, 
$\xi $ is a
$p$-adic integer). Proceeding with this note, we conclude that for
$k,\tau \in \{1,2,\ldots  ,p-1\}$ the following equalities hold in ${\Bbb Q}_{p}$:
$$
\displaylines{\sum^{p-1}_{k=1}{\frac {(-1)^{k-1}} {k}}\sum^{p-k-1}_{\tau =1}(-1)^{\tau -1}{\frac {\Delta ^{(\tau +k)p^{\lambda -1}_{}}f(x)} {\tau p^{\lambda -1}}}
\hfill\cr
\hfill
=\sum^{p-1}_{m=1}(-1)^{m}\sum^{}_{k+\tau =m} {\frac {1} {k\tau }}\cdot {\frac {\Delta ^{mp^{\lambda -1}_{}}f(x)} {p^{\lambda -1}}}
=2\sum^{p-1}_{m=1}(-1)^{m}\sum^{m-1}_{\tau =1} {\frac {1} {\tau }}\cdot {\frac {\Delta ^{mp^{\lambda -1}_{}}f(x)} {mp^{\lambda -1}}},
\quad (16)\cr}$$
\noindent since for $k,\tau \in \{1,2,\ldots  ,p-1\}$ it is obvious that
$$
\sum^{}_{k+\tau =m}{\frac {1} {k\tau }} = \sum^{}_{k+\tau =m}{\frac {1} {(m-\tau )\tau }} = {\frac {1} {m}}\sum^{}_{k+\tau =m}({\frac {1} {\tau }} + {\frac {1} {m-\tau }}) = {\frac {2} {m}} \sum^{m-1}_{\tau =1}{\frac {1} {\tau }} .
$$
\noindent Besides, as it was shown during the proof of 4.5, both $\alpha (x)$ and $\beta _{k}(x)$ are $p$-adic integers for $k=1,2,\ldots  ,p-1$ and $x\in {\Bbb Z}_{p}$; thus
$$
2hp\alpha (x)=h {\frac {\Delta ^{2p^{\lambda }}f(x)} {p^{\lambda }}};\quad  hp\beta _{k}(x)=h {\frac {\Delta ^{kp^{\lambda -1}+p^{\lambda }}f(x)} {kp^{\lambda -1}}},
\eqnum{17}$$
\noindent where all the factors are $p$-adic integers. Now the assertion of the lemma
follows from (15), (16), (17) and definition of the function $\theta$.\qed
\enddemo
\demo
{Proof of the theorem 4.1}
Finishing the proof of the theorem 4.1, note that according
to 4.4 there holds an inequality $N_{2}(f)\le \lambda (f)+1$. Thus, by 3.14 it is sufficient only to show that
if $p\neq 3$ and $f$ is transitive modulo $p^{\lambda (f)+1}$ then it is transitive modulo
$p^{\lambda (f)+2}$. In turn, for this purpose in view of 3.15 it is sufficient only
to prove that
$$
f^{p^{\lambda +1}}(x)\not\equiv x\pmod{ p^{\lambda +2}}
\eqnum{1}$$
\noindent at least for one $x\in {\Bbb Z}_{p}$. Further we merely calculate 
$f^{p^{\lambda +1}}(x)\bmod{ p^{\lambda +2}}$.
\par Under the assumptions we have made above, $f$ 
is transitive modulo $p^{\lambda }$, since $f$ is
compatible. Then by 3.15 we conclude that for all $x\in {\Bbb Z}_{p}$
$$
f^{p^{\lambda }}(x)=x+p^{\lambda }\xi (x),\quad \xi (x)\not\equiv 0\pmod p,
\eqnum{2}$$
\noindent where $\xi \colon{\Bbb Z}_{p}\rightarrow {\Bbb Z}_{p}$ is a function defined everywhere on ${\Bbb Z}_{p}$.
\par
We assert that for each $i=0,1,2,\ldots  $ the following congruence holds:
$$
\displaylines{f^{p^{\lambda }+i}(x)\equiv f^{i}(x)+p^{\lambda }\xi (x)\prod^{i-1}_{j=0}f^{\prime }_2(f^{j}(x))\hfill\cr                                  
\hfill {}+p^{\lambda +1}\xi (x)^{2}\prod^{i-1}_{j=0}f^{\prime }_2(f^{j}(x))\sum^{i-1}_{k=0} {\frac {\theta (f^{k}(x))} {f^{\prime }_2(f^{k}(x))}}\prod^{k-1}_{\tau =0}f^{\prime }_2(f^{\tau }(x))\pmod{ p^{\lambda +2}}\qquad (3)\cr }$$
\noindent Recall that the sum (resp., product) over the empty set of indexes is
assumed to be 0 (resp., 1). Note also that since $f$ is transitive modulo
$p^{\lambda +1}$ it is bijective modulo $p^{\lambda +1}$. Consequently, $f$ is bijective modulo
$p^{\lambda },\ldots  ,p^{2},p$ since $f$ is compatible. Hence $f^{\prime }_1(x)\not\equiv 0\pmod p$ for all $x\in {\Bbb Z}_{p}$ (see
the proof of 3.9) and thus $f^{\prime }_2(x)\not\equiv 0\pmod p$ either (as $f^{\prime }_2(x)\equiv f^{\prime }_1(x)\pmod p$).
So all the denominators in (3) have multiplicative inverses in ${\Bbb Z}_{p}$; thus,
during the proof of (3) and further we assume that all the calculations
are performed in ${\Bbb Z}_{p}$.
\par
We can easily prove (3) by the induction on $i$. If $i=0$, then (3)
trivially follows from (2). Assume that (3) is true for $i=m-1$. As
$$
f^{p^{\lambda }+m}(x)=f(f^{p^{\lambda }+m-1}(x))
\eqnum{4}$$
\noindent then,  assuming in (3) that $i=m-1$, substituting (3) to (4),
applying 
4.5 and a congruence 
$(f^k(x))^{\prime}_2\equiv\prod^{k-1}_{j=0}f^{\prime }_2(f^{j}(x))\pmod {p^2}$,
we prove the congruence (3) for $i=m$, in view of
compatibility of $f$,  by obvious direct
calculations. We omit details.
\par
Now we apply (3) to calculate $f^{p^{\lambda +1}}(x)\bmod{ p^{\lambda +2}}$. Put
$$
\displaylines{A _{i}(x)=(f^i(x))^{\prime}_2=\prod^{i-1}_{j=0}f^{\prime }_2(f^{j}(x));\hfill\cr
B _{i}(x)=(f^i(x))^{\prime}_2\sum^{i-1}_{k=0} {\frac {(f^{k }(x))^{\prime }_2} {f^{\prime }_2(f^{k}(x))}}\theta (f^{k}(x))=
\hfill\cr\hfill=\Biggl(\prod^{i-1}_{j=0}f^{\prime }_2(f^{j}(x))\Biggl)\cdot\Biggr(\sum^{i-1}_{k=0} {\frac {\theta (f^{k}(x))} {f^{\prime }_2(f^{k}(x))^{2}}} \prod^{k}_{\tau =0}f^{\prime }_2(f^{\tau }(x))\Biggr).\cr}
$$
\noindent Lemma 4.6 implies that 
$$
f^{\prime }_2(a+p^{\lambda }h)\equiv \left\{\matrix f^{\prime }_2(a)\pmod {p^{2}},&\hbox{if }h=0;\\ f^{\prime }_2(a)\pmod p,&\hbox{if }h\neq 0.\endmatrix\right.
\eqnum{5}$$
\noindent As $f$ is transitive modulo $p^{\lambda }$, then (5) implies that
$f^{\prime }_2(f^{k}(x))\equiv f^{\prime }_2(f^{r}(x))\pmod p$ 
as soon as $k\equiv r\pmod {p^\lambda}$. Besides, by 4.5 the
latter condition implies that $\theta (f^{k}(x))\equiv \theta (f^{r}(x))\pmod p$.
\par
 Further,
$$
\prod^{p^{\lambda }-1}_{\tau =0}f^{\prime }_2(f^{\tau }(x))\equiv 1\pmod p.
\eqnum{6}$$
\noindent This has been already proven in 3.15 (see proof of (6) there), since 4.5 implies that
$N_{1}(f)\le \lambda $. Consequently,
$$
\prod^{_{}k}_{\tau =0}f^{\prime }_2(f^{\tau }(x))\equiv \prod^{_{}r}_{\tau =0}f^{\prime }_2(f^{\tau }(x))\pmod p
$$
\noindent as soon as $k\equiv r\pmod {p^\lambda}$.
\par
Finally we conclude that for every $t\in {\Bbb N}$
$$
B _{t p^\lambda}(x)\equiv t\sum^{p^{\lambda }-1}_{\tau =0} {\frac {\theta (f^{\tau }(x))} {f^{\prime }_2(f^{\tau }(x))^{2}}}\prod^{\tau }_{\nu =0}f^{\prime }_2(f^{\nu }(x))\equiv tB _{p^\lambda}(x)\pmod p.
\eqnum{7}$$
\noindent Now we calculate $A_{t p^\lambda}(x)\pmod{p^{2}}$ for $t\in {\Bbb N}$. The congruence (3) in view of
(6) implies that
$$
f^{kp^{\lambda }+\tau }(x)\equiv f^{\tau }(x)+kp^{\lambda }\xi (x)\prod^{\tau -1}_{j=0}f^{\prime }_2(f^{j}(x))\pmod{p^{\lambda +1}}
\eqnum{8}$$
\noindent for all $k\in {\Bbb N}$ and all $\tau \in \{0,1,\ldots  ,p^{\lambda }-1\}$. As 
$$
A_{t p^\lambda}(x)=\prod^{t-1}_{k=0}\prod^{p^{\lambda }-1}_{\tau =0}f^{\prime }_2(f^{kp^{\lambda }+\tau }(x)),
$$
\noindent then in view of (5) and 4.6 the congruence (8) implies that 
$$
A_{t p^\lambda}(x)=\prod^{t-1}_{k=0}\prod^{p^{\lambda }-1}_{\tau =0}f^{\prime }_2\biggl (f^{\tau }(x)+kp^{\lambda }\xi (x)\prod^{\tau -1}_{j=0}f^{\prime }_2(f^{j}(x))\biggr )\pmod{p^{2}},
$$
\leftline{or, applying 4,6, }
$$
\displaylines{A_{t p^\lambda}(x)=\prod^{t-1}_{k=0}\prod^{p^{\lambda }-1}_{\tau =0}\Biggl (f^{\prime }_2 (f^{\tau }(x))+2kp\xi (x)\theta (f^{\tau }(x))\prod^{\tau -1}_{j=0}f^{\prime }_2(f^{j}(x))\Biggr )
\hfill\cr 
\quad\equiv \prod^{t-1}_{k=0}\Biggl (\prod^{p^{\lambda }-1}_{\tau =0}f^{\prime }_2(f^{\tau }(x))
\hfill\cr\hfill
{}+2kp\xi (x)\sum^{p^{\lambda }-1}_{s=0}\theta (f^{s}(x)){\frac {\prod^{p^{\lambda }-1}_{j=0}f^{\prime }_2(f^{j}(x))} {f^{\prime }_2(f^{s}(x))}}\prod^{s-1}_{j=0}f^{\prime }_2(f^{j}(x))\Biggr )\pmod {p^{2}}.\qquad (9)\cr}$$
\par
According to (6),
$$
\prod^{p^{\lambda }-1}_{j=0}f^{\prime }_2(f^{j}(x))=1+p\epsilon 
$$
\noindent for suitable $\epsilon \in {\Bbb Z}_{p}$; consequently, (9) implies that
$$
 \displaylines{A_{t p^\lambda}(x)\equiv \prod^{t-1}_{k=0}\Biggl (1+p\epsilon +2kp\xi (x)\sum^{p^{\lambda }-1}_{s=0}\theta (f^{s}(x)){\frac {\prod^{s-1}_{j=0}f^{\prime }_2(f^{j}(x))} {f^{\prime }_2(f^{s}(x))}}\Biggl )\equiv
\hfill\cr
\quad
\equiv 1+tp\epsilon +2p\xi (x)\Biggl(\sum^{t-1}_{k=0}k\Biggr)\cdot\Biggl(\sum^{p^{\lambda }-1}_{s=0}\theta (f^{s}(x)){\frac {\prod^{s-1}_{j=0}f^{\prime }_2(f^{j}(x))} {f^{\prime }_2(f^{s}(x))^{2}}}\Biggr)\equiv
\hfill\cr
\hfill
{}\equiv 1+tp\epsilon +pt(t-1)\xi (x)\sum^{p^{\lambda }-1}_{s=0}\theta (f^{s}(x)){\frac {\prod^{s-1}_{j=0}f^{\prime }_2(f^{j}(x))} {f^{\prime }_2(f^{s}(x))^{2}}}\pmod {p^{2}}.
\qquad (10)\cr}$$
\noindent Now combining together (2), (3), (7) and (10) we conclude that
$$
\displaylines{f^{(t+1)p^{\lambda }}(x)\equiv f^{tp^{\lambda }+p^{\lambda }}(x)
\hfill\cr
\hfill
{}\equiv f^{tp^{\lambda }_{}}(x)+p^{\lambda }\xi (x)+\epsilon tp^{\lambda +1}\xi (x)+p^{\lambda +1}t^{2}\xi (x)^{2}B _{p^\lambda}(x)\pmod{ p^{\lambda +2}}.
\qquad (11)\cr}$$
\noindent Finally, combining (11), (2) with obvious induction on $n$ we obtain
that
$$
\displaylines{f^{np^{\lambda }_{}}(x)\equiv x+np^{\lambda }\xi (x)+\epsilon p^{\lambda +1}\xi (x){\frac {n(n-1)} {2}}
\hfill\cr\hfill
{}+p^{\lambda +1}\xi (x)^{2}B _{p^\lambda}(x){\frac {n(n-1)(2n-1)} {6}}\pmod{ p^{\lambda +2}}\cr}
$$
\noindent or, in particular, 
$$f^{p^{\lambda +1}_{}}(x)\equiv x+p^{\lambda +1}\xi (x)\pmod{ p^{\lambda +2}},$$
\noindent since $p\neq 2,3$. But the latter congruence in view of (2) implies that 
$$
f^{p^{\lambda +1}_{}}(x)\not\equiv x\pmod{ p^{\lambda +2}}.
$$
\noindent This finally proves the theorem 4.1\qed
\enddemo
\remark
{Note}With the use of theorem 4.1 we can determine whether a given
integer-valued and compatible polynomial $f(x)\in {\Bbb Q}_{p}[x]$ is ergodic.
Represent $f(x)$ in the form $f(x)={g(x)\over r}$, where $r\in {\Bbb Z}_{p}$ and $g(x)\in {\Bbb Z}_{p}[x]$ and at least
one coefficient of $g(x)$ is coprime with $p$. In fact, we can take $r$ to be a common denominator 
of all coefficients of $f(x)$ represented as irreducible fractions. Here we assume that $f(x)$ is represented 
in the basis $(x)_0=1, (x)_1=x, (x)_2=x(x-1), \ldots$ of descending factorial
powers,  or in a standard basis $1, x, x^2, \ldots.$
Then $\rho (f)=\ord_{p}\,r$, and $\rho (f)$ does not depend on the choice of the basis. We
recall that $p^{\ord_{p}\,r}$ is the greatest power of $p$ which is a factor of $r$. Now we easily
find $\lambda (f)$ and determine whether $f$ is transitive on 
${\Bbb Z}/p^{\lambda (f)+1}$ 
(e.g., by direct
calculations). In view of 4.1 for $p\ne 2,3$ this is equivalent to the ergodicity
of $f(x)$ (for $p=3$ one should study transitivity of $f$ on ${\Bbb Z}/p^{\lambda (f)+2}$).

Moreover, it is possible for each prime  $p$ to determine, whether a polynomial
 $f(x)\in {\Bbb Q}_{p}[x]$ is integer-valued, compatible and ergodic, by
 calculating its values at $O(\deg f)$ points.
Namely, the following is true.
\endremark 
\proclaim
{4.7 Proposition}{A polynomial $f(x)\in {\Bbb Q}_{p}[x]$ is integer-valued,
compatible and ergodic iff the mapping 
$$z\mapsto f(z)\bmod p^{\lfloor
\log_p (\deg f)\rfloor +3},$$ with $z$ ranging over $
\{0,1,\ldots,p^{\lfloor
\log_p (\deg f)\rfloor +3}-1\}$, defines a compatible and transitive function
on the residue class ring  $\Bbb Z/p^{\lfloor
\log_p (\deg f)\rfloor +3}$.}
\endproclaim
\demo
{Proof} Coefficients $a_i\in\Bbb Q_p$ ($i=0,1,\ldots,d$) of the
polynomial $f(x)$ of degree $d$, which is represented in the form
$
f(x)=\sum^{d}_{i=0}a_{i}{x\choose i}
$
(see $(\diamondsuit)$ of section 2), are defined by the values  this polynomial
$f(x)$ takes at the points  $0,1,\ldots,d$.
In other words, all values $f(0), f(1),\ldots,f(d)$ are $p$-adic integers
iff all coefficients  $a_i\in\Bbb Q_p$ ($i=0,1,\ldots,d$)
are $p$-adic integers, i.e., iff a polynomial $f(x)$ is integer-valued
(see the beginning of section  2). By the analogy, in view of the theorem
  2.1, a polynomial $f(x)$
preserves all congruences of the ring  $\Bbb Z/p^{\lfloor \log_p d \rfloor+1}$ iff
$\|a_i\|\le    p^{-\lfloor\log_p i\rfloor}$ for all
$i=1,2,\ldots,d$, i.e., iff  $f(x)$ is compatible on $\Bbb Z_p$.
In other words,  to determine whether a polynomial $f(x)$ is integer-valued
and compatible it is sufficient (and necessary) 
to determine whether it induces a compatible function on the ring
ограничиться проверк$\Bbb Z/p^k$ for some (arbitrarily fixed)
произвольным образом$k\ge \lfloor\log_p d\rfloor +1$.

In force of theorem  4.1, for $p\ne 2$, an integer-valued and compatible
polynomial $f(x)$ is ergodic iff 
it is transitive modulo $p^k$ for any arbitrarily fixed  $k\ge\lambda (f)+2$. 
Representing $f(x)$ as 
$
f(x)=b_{0}+\sum^{d}_{i=1}b_{i}p^{\left\lfloor{\log_{p}i}\right\rfloor}{ {x}\choose {i}},
$ где $b_{j}\in {\Bbb Z}_{p}$ for $j=0,1,2,\ldots $, we conclude that 
$\rho(f)$ is the least nonnegative integer rational, which is not less
that each of  
$
\ord_{p}\,i!-\left\lfloor{\log_{p}i}\right\rfloor-\ord_{p}\,b_{i}
$ ($i=1,2,\ldots,d$). Thus, since a function $\ord_{p}\,i!-\left\lfloor{\log_{p}i}\right\rfloor$ 
is nondecreasing (see proof of the lemma 4.2), then each 
$k\in \Bbb N$, which satisfies inequality $ 2 {\frac {p^{k}-1} {p-1}} - k >\ord_{p}\,d!-\left\lfloor{\log_{p}d}\right\rfloor$,
will satisfy inequality $k\geqslant\lambda(f)$. Yet since $\ord_{p}\,d!={\frac {1} {p-1}}(d-\wt_{p}\,d)$,  
where $\wt_{p}\,d$  is a sum of all digits in $p$-base expansion of $d$, then,
choosing any
$k\in \Bbb N$, which satisfy inequality 
$$ 2 {\frac {p^{k}-1} {p-1}} - k >\frac{d}{p-1},\tag{1}$$
we obtain that $k\ge\lambda(f)$. Elementary considerations, however, show that 
$k=\lfloor\log_p d\rfloor+1$ satisfies inequality (1), thus proving 
the proposition for $p\ne 2$.

In case $p=2$ a polynomial $f(x)\in \Bbb Q_2[x]$ of degree $d$ is integer-valued,
compatible and ergodic iff 
it is of a form
$$
f(x)=1+x+\sum^{d}_{i=0}b_{i}\,2^{\left\lfloor{\log_{2}(i+1)}\right\rfloor+1}{ {x}\choose {i}},\tag{2}
$$
where $b_{i}\in {\Bbb Z}_{2}$, $i=0,1,2,\ldots,d $ (see theorem 2.3). Since
coefficients 
of the polynomial $f(x)$ in its representation in  a basis $x\choose i$, $i=0,1,2,\ldots$, are
uniquelly defined by the values of $f(z)$ at the points $z=0,1,\ldots,d$, then
to verify conditions 
(2) for the polynomial $f(x)$ it is sufficient to calculate its values at
the points  $z=0,1,\ldots,2^r-1$,
where $r\in\Bbb N$ is an arbitrarily fixed number satisfying the inequality
$d\le     2^r-1$.
So one can take, for instance, $r=\lfloor\log_2 (d+1)\rfloor+1$, or $r=\lfloor\log_2 d\rfloor+3$.
This finishes the proof of 4.7. \qed        
\enddemo
\remark
{Note} Proposition 4.4 shows that for $p\ne 2$ a function $f\in\Cal A$ satisfies
assumptions of the proposition  3.9; hence, since $N_1(f)\le     N_2(f)$, 
a function $f$ 
preserves measure iff it is bijective modulo $p^{\lambda(f)+2}$.
By the argument similar to those of the proof of 
proposition  4.7,
one could prove the following
\endremark
\proclaim
{4.8 Proposition}{A polynomial $f(x)\in {\Bbb Q}_{p}[x]$ is integer-valued,
compatible and measure-preserving iff
the mapping
$$z\mapsto f(z)\bmod
p^{K_f},$$ 
with $K_f={\lfloor
\log_p (\deg f)\rfloor +3}$ and $z$ ranging over
$0,1,\ldots,p^{K_f}-1$, induces a compatible and bijective function 
on the ring 
$\Bbb Z/p^{K_f}$.\qed}
\endproclaim

Again, estimates of $\zeta(f)$ and  $\eta(f)$ we mentioned at the beginning
of the section, might be sharpened for various important proper subclasses of 
$\Cal A$ in comparison with given by the theorem  4.1 and propositions
4.7 and 4.8.  A case of analytic on  $\Bbb Z_p$ functions (i.e., functions
which can be represented by convergent everywhere on $\Bbb Z_p$ power series)
seems to be of importance.

It is well known (see e.g. \cite{3, Ch. 14. Section 4}) that power series
$\sum^{\infty }_{i=0}c_{i}x^{i}$ ($c_i\in\Bbb Q_p,\ i=0,1,2\ldots$) converges
everywhere on $\Bbb Z_p$ iff $\lim\limits_{i\to\infty}^p c_{i}=0$;
under the latter condition the series defines a continues function on $\Bbb Z_p$.
Of course, in general this function may not be integer-valued,
not speaking about compatibility. Consider, however, a particular case,
when all coefficients $c_i$ are $p$-adic integers. Namely, in the ring
$\Bbb Z_p[[x]]$ of all formal power series in variable $x$ over a ring
$\Bbb Z_p$ consider a set ${\Cal C}(x)$ of all series 
$$
s(x)=\sum^{\infty }_{i=0}c_{i}x^{i}\quad (c_i\in\Bbb Z_p,\ i=0,1,2\ldots),
$$
which converge everywhere on  $\Bbb Z_p$. In other words, $s(x)\in {\Cal C}(x)$
iff $\lim\limits_{i\to\infty}^p c_{i}=0$. Under these assumptions the series
$s(x)\in{\Cal C}(x)$ defines on  $\Bbb Z_p$ an integer-valued function
$s:\Bbb Z_p\rightarrow \Bbb Z_p$. It turnes out that this function $s$
is uniformly differentiable and has integer-valued derivative everywhere
on  $\Bbb Z_p$.

Consider a formal derivative $s^\prime(x)\in \Bbb Z_p[[x]]$  of the series
$s(x)$:
$$s^{\prime} (x)=\sum_{i=1}^\infty ic_ix^{i-1}.$$
Since $0\le   \|ic_i\|_p=\|i\|_p\|c_i\|_p\le    \|c_i\|_p$, 
and $\lim\limits_{i\to\infty}^p c_{i}=0$, then $\lim\limits_{i\to\infty}^p ic_{i}=0$,
and hence $s^\prime (x)\in{\Cal C}(x)$. We assert that the function $s^\prime: \Bbb Z_p\rightarrow \Bbb Z_p$
is a derivative of a function $s: \Bbb Z_p\rightarrow \Bbb Z_p$ with respect
to $p$-adic distance.

Indeed, it is known that in the ring ${\Bbb Z}_{p}[[x,y]]$
of all formal power series in variables $x, y$ over ${\Bbb Z}_{p}$ the
following equality holds:
$$
s(x+y)=\sum^{\infty }_{i=0}{\frac {s^{(i)}(x)} {i!}}y^{i},
$$
where $s^{(i)}(x)\in \Bbb Z_p[[x]]\ (i=1,2,\ldots)$ is $i$th formal derivative
of the series $s(x)$, and $s^{(0)}(x)=s(x)$. By the assertion proven above, 
$s^{(i)}(x)\in{\Cal C}(x)$ for all $i=0,1,2,\ldots$. Thus,
$$
{\frac {s^{(i)}(u)} {i!}}=\sum^{\infty }_{j=i}c_{j}{ {j}\choose {i}}u^{j-i}\in
\Bbb Z_p
$$
for each  $u\in {\Bbb Z}_{p}$. But
$$
\left\Vert{{\frac {s^{(i)}(u)} {i!}}}\right\Vert_{p}=\left\Vert{\sum^{\infty }_{j=i}c_{j}{ {j}\choose {i}}u^{j-i}}\right\Vert_{p}\leqslant     \max\{\| c_{j}\| _{p}:j=i,i+1,\ldots \},
$$
and consequently,
$$
\lim\limits^p_{i\to\infty}{\frac {s^{(i)}(u)} {i!}}=0,$$
since $ \lim\limits^p_{i\to\infty} c_{i}=0$. Thus, for each  $u\in\Bbb Z_p$
we have that
$$
s(u+y)=\sum^{\infty }_{i=0}{\frac {s^{(i)}(u)} {i!}}y^{i}\in{\Cal C}(y).\eqnum{\spadesuit} 
$$
Finally, if $s(x)\in\Cal
C(x)$, then Taylor series $(\spadesuit )$ at the point $u\in\Bbb Z_p$ converges
to $s$
everywhere on $\Bbb Z_p$. In particular, for $h\in\Bbb Z_p$ we obtain 
$$
s(u+h)=s(u)+s^\prime(u)h+\alpha(u,h),
$$
with $\lim\limits^p_{h\to 0}\dfrac{\alpha(u,h)}{h}
=\lim\limits^p_{h\to 0}h\sum^{\infty }_{i=2}{\frac {s^{(i)}(u)} {i!}}h^{i-2}=0$,
since $\sum^{\infty }_{i=2}{\frac {s^{(i)}(u)} {i!}}h^{i-2}\in\Bbb
Z_p$ in view of the equality $\lim\limits^p_{i\to\infty}{\frac {s^{(i)}(u)} {i!}}=0$,
which just has been proven above. So, $s^\prime(u)$ is a derivative of the function
$s$ at the point $u$. Thus, the set ${\Cal C}(x)$ is closed with respect
to differentiations, and all functions, defined by series of ${\Cal C}(x)$,
are infinitely many times differentiable. 

Further, let
$$
s(x)=\sum^{\infty }_{i=0}s_{i}{ {x}\choose {i}}
$$
be an interpolation series for the function  $s(x)\in {\Cal C}(x)$. We
assert that $\dfrac{s_i}{i!}$ is $p$-adic integer for all $i=0,1,2,\ldots  $. 
Actually,
$$
s(x)=\sum^{\infty }_{k=0}c_{k}x^{k}=\sum^{\infty }_{k=0}c_{k}\sum^{k}_{i=0}S_{2}(k,i)i!{ {x}\choose {i}}=\sum^{\infty }_{i=0}i!{ {x}\choose {i}}\sum^{\infty }_{k=i}S_{2}(k,i)c_{k},
$$
where $S_{2}(k,i)$ is  Stirling number. Since $\lim\limits_{i\to\infty}^p c_{i}=0$,
then $\lim\limits_{k\to\infty}^p S_{2}(k,i)c_{k}=0$, because all Stirling
numbers $S_{2}(k,i)$ are integer rationals, i.e., $\|{S_{2}(k,i)}\|_{p}\le   1$.
Consequently, the series $\sum^{\infty }_{k=i}S_{2}(k,i)c_{k}$ converges
to some $A_{i}\in {\Bbb Z}_{p}$ for all  $i=0,1,2,\ldots  $. This proves
our assertion, since
$$s_{i}=i!A_{i}\quad (i=0,1,2,\ldots  ).\eqnum\bigstar $$

Put
$$\Cal B(x)=\bigg\{f(x)=\sum^{\infty }_{i=0}a_{i}{ {x}\choose {i}}: \dfrac{a_i}{i!}\in\Bbb
Z_p, \quad i=0,1,2,\ldots\bigg\}.$$
In other words, ${\Cal B}(x)$ is a ring of all formal descending factorial
power series over $\Bbb Z_p$. Each series $f(x)\in\Cal B(x)$ correctly
defines on $\Bbb Z_p$ an integer-valued and uniformly continuos function 
$f:\Bbb Z_p\rightarrow \Bbb Z_p$ (see the beginning of the section 2).
This function $f$ is compatible in view of 2.1, since we have shown during
the proof of the lemma 4.2 that $\ord_p\,(i!)-\lfloor\log_pi\rfloor$
is nonnegative and nondescending function on $\Bbb N_0$. Denote via ${\Cal B}$ (respectively,
via ${\Cal C}$) a class of all functions defined by all series of ${\Cal B(x)}$
(respectively, of ${\Cal C(x)}$). Obviously, ${\Cal B}(x)$, ${\Cal B}$,
 ${\Cal C}(x)$, ${\Cal C}$ are rings.

Further, any two distinct series of ${\Cal B}(x) $ (respectively, of
${\Cal C}(x) $ ) define two distinct functions on $\Bbb Z_p$: for the
series
of $\Cal B(x)$ see the beginning of the section 2.  As for the series of
$\Cal C(x)$, note, that the above mentioned interpolation series for  
$s(x)\in\Cal C(x)$ defines a function, which is identically 0 on $\Bbb
Z_p$ iff all its coefficients $s_i$ are 0 (hence, $A_i=0$, $i=0,1,2,\ldots$),
see $(\bigstar)$. Yet $A_{i}=\sum^{\infty }_{k=i}S_{2}(k,i)c_{k}$, hence
$c_{i}=\sum^{\infty }_{k=i}S_{1}(k,i)A_{k}=0$, where $S_1(k,i), S_2(k,i)$
are Stirling numbers of respective kind, and the assertion follows. 
Thus, the rings ${\Cal B}(x)$ and ${\Cal B}$ (respectively,
${\Cal C}(x)$ and ${\Cal C}$) are isomorphic; so further we do not differ series
from the function it defines.

Note also that the incluion ${\Cal B}\supset {\Cal C}$  (see $(\bigstar)$)
is strict. Obviously, $f(x)=\sum_{i=0}^\infty (x)_i\in\Cal B$, since
$f(x)=\sum_{i=0}^\infty i!{x\choose i}$. Yet $f(x)\notin \Cal C$. Moreover,
this
function is not even analytic on $\Bbb Z_p$: according to \cite{3, Ch.4,
Theorem 4} a function represented by the interpolation series $(\diamondsuit )$
of the section 2 is analytic on $\Bbb Z_p$ iff  $\lim\limits^p_{i\to\infty}\frac{a_i}{i!}=0$.

So, a function of $\Cal B$ (in contrast to one of $\Cal C$), generally
speaking, can not be represented by Taylor series which is convergent everywhere
on $\Bbb Z_p$. Newertheless, all functions of $\Cal B$ are differentiable
at all points of $\Bbb Z_p$, and $\Cal B$ is closed with respect to differentiations:
if $f\in \Cal B$, then $f^\prime\in \Cal B$.

To prove the latter assertion, recall  that
a uniformly continuous on $\Bbb Z_p$ function $f$, which is represented
by the interpolation series $(\diamondsuit)$, is differentiable everywhere
on $\Bbb Z_p$ iff
$$\lim\limits^p_{i\to\infty}\frac{a_{i+n}}{i}=0\eqno{(\blacklozenge)}$$ 
for all $n\in\Bbb N_0$ (see \cite{3, Ch. 13, Theorem 2}). The lattter condition
obviously holds for  $f\in\Cal B$, since  
 $\ord_p\,a_i\geqslant\ord_p\,(i!)=\frac{1}{p-1}(i-\wt_p\,i)$, and 
 $\lfloor\log_pi\rfloor\geqslant\ord_p\,i$ for
 all $i=0,1,2,\ldots$. Thus, a derivative $f^\prime$ of the function $f$
 is defined everywhere on $\Bbb Z_p$, and
$$f^\prime(x)=\sum_{i=1}^\infty (-1)^{i+1}\frac{\Delta^if(x)}{i},$$
in case this series is convergent. Yet
$\frac{\Delta^if(x)}{i}={1\over i}\sum_{j=i}^\infty a_j{x\choose{j-i}}$,
consequently,
$$\sum_{i=1}^\infty (-1)^{i+1}\frac{\Delta^if(x)}{i}=\sum_{k=0}^\infty{x\choose
k} \sum^\infty_{i=1}(-1)^{i+1}\frac{a_{k+i}}{i}.$$
But the series $\sum^\infty_{i=1}(-1)^{i+1}\frac{a_{k+i}}{i}$, in view of
 $(\blacklozenge)$, for each $k\in\Bbb N_0$ converges to a certain $S_k\in\Bbb Q_p$, 
 and
 $\ord_p\,\frac{a_{k+i}}{i}=\ord_p\,a_{k+i}-\ord_p\,i\ge\ord_p\,((k+i)!)-\lfloor\log_pi\rfloor
={1\over{p-1}}(i+k-\wt_p\,(i+k)) -\lfloor\log_pi\rfloor=
{1\over{p-1}}(i-\wt_p\,i) -\lfloor\log_pi\rfloor+{1\over{p-1}}(k-\wt_p\,k)+
{1\over p-1}(\wt_p\,k-\wt_p\,(i+k)+\wt_p\,i)\ge{1\over{p-1}}(k-\wt_p\,k)=\ord_p\,(k!)$.
(The latter inequality holds since
${1\over{p-1}}(i-\wt_p\,i)\ge
\lfloor\log_pi\rfloor$ and ${1\over p-1}(\wt_p\,k-\wt_p\,(i+k)+\wt_p\,i)=\ord_p\,{{i+k}\choose
i}\ge 0$)). Thus, ${\dfrac{S_k}{k!}}\in\Bbb Z_p$ for all $k\in\Bbb
N_0$; hence $f^\prime\in\Cal B$. 

With the use of these results now we are able to prove the following
\proclaim
{4.9 Theorem} A function $f\in \Cal B$ preserves measure iff it is bijective
modulo $p^2$. The function $f$ is ergodic iff it is transitive modulo $p^2$
{\rom(}for $p\ne 2,3${\rom)}, or modulo $p^3$ {\rom(}for $p\in\{2,3\}${\rom)}.
\endproclaim
\demo
{Proof} The definition of $\Cal B$ immediately implies that $\rho(f)=0$
for each $f\in {\Cal B}$, hence, $\lambda (f)=1$. Thus, for $p\ne 2$ 
the second assertion of the theorem follows from 4.1.

To prove the first assertion, in view of 3.9 it is sufficient to demonstrate
that $f$ is uniformly differentiable modulo $p$, and $N_1(f)\le     1$;
that is
$$f(z+p^kr)\equiv f(z)+p^krf^\prime (z)\pmod {p^{k+1}} \tag1$$
for all $z,r\in\Bbb Z_p$ and $k=1,2,\ldots$. Since $f,f^\prime\in\Cal B$,
these both functions are compatible, so it is sufficient to prove (1) for 
$z,r\in\Bbb N_0$. Since for $r=0$ the congruence (1) is trivial, we may
additionally assume that $p^kr=n\in\Bbb N$.

Further, since $\frac{f(z+n)-f(z)}{n}=\sum_{i=1}^\infty{{n-1}\choose{i-1}}\frac{\Delta^i
f(x)}{i}$, а $f^\prime (z)=\sum_{i=1}^\infty (-1)^{i+1}\frac{\Delta^if(z)}{i}$,
then to prove (1) it is sufficient to prove that
$$\sum_{i=1}^\infty\biggl({{n-1}\choose{i-1}}-(-1)^{i+1}\biggr)\frac{\Delta^if(z)}{i}\equiv 0\pmod p.\tag2$$
Yet $\frac{\Delta^if(x)}{i}={1\over i}\sum_{j=i}^\infty a_j{x\choose{j-i}}$,
thus, in view of 4.2, for $p\ne 2$ there holds a congruence 
 $\frac{\Delta^if(x)}{i}\equiv
0\pmod p$ for
all $i\ge 2p$. So within this case (2) is equivalent  to the congruence
$$\sum_{i=1}^{2p-1}\biggl({{n-1}\choose{i-1}}-(-1)^{i+1}\biggr)\frac{\Delta^if(z)}{i}\equiv 0\pmod p.\tag3$$
Since $f$ is compatible, then
$\frac{\Delta^if(x)}{i}\not\equiv 0\pmod p$ only for, might be, $i=sp^m$,
($m\in\Bbb N_0$, $s\in\{1,2,\ldots,p-1\}$) --- see \cite{11, lemma 3.4}.
Now, since $n=p^kr$, (3) immediately follows from the already mentioned
Lucas theorem, thus proving the first assertion of 4.9 for $p\ne 2$.

Now, if $p=2$, then (2) is equivalent to
$$\sum_{i=0}^\infty\biggl({{2^kr-1}\choose{2^i-1}}+1\biggr)\frac{\Delta^{2^i}f(z)}{2^i}\equiv 0\pmod 2.\tag4$$
Yet since ${a_j\over j!}\in \Bbb Z_2$ for all $j=0,1,2,\dots$, then
$\ord_2\,a_{2^i+m}\ge\ord_2\,(2^i)!= 2^i-1$ for all $m=0,1,2,\ldots$; 
consequently, $\frac{\Delta^{2^i}f(z)}{2^i}\not\equiv
0\pmod2$ only for, might be, $i=0$, thus proving (4).

Finally, the rest part of the assertion of theorem 4.9 for $p=2$ follows
from 2.3: as
 $\ord_{2}\,i!\le  \ord_{2}\,a_{i}$ for all $i=0,1,2,\ldots  $, and
 $\ord_{2}\,i!=i-\wt_2(i)$, 
then by an elementary argument it is not difficult to show that
$\lfloor\log_2 (i+1)\rfloor+1\le i-\wt_2(i)\le \ord_2\,a_i$ for $i\ge
4$; and $\ord_2\,a_i\ge 3$. This implies that necessary and sufficient
conditions of ergodicity of a function expressed as interpolation series
$(\diamondsuit)$ of section 2
hold for all coefficients $a_i$ with $i\ge 4$. These conditions for the
rest of the coefficents are equivalent to the transitivity of $f$ modulo 8,
since $a_i\equiv 0\pmod 8$ for $i\geqslant 4$.\qed
\enddemo
\remark
{Note} Theorem 4.9 demonstrates that sufficient and necessary conditions
of transitivity modulo $p^n$ for the polynomials with integer rational
coefficients established by M. V. Larin in \cite{15} remain valid for a
wider class (namely, $\Cal B$) of functions. It turnes out, however, that
all these functions modulo each $p^n$ could be expressed as polynomials
with rational integer coefficients.
\endremark

Namely, from the  definition of a class $\Cal B$ it easily follows that
each function $f\in\Cal B$ is uniformly approximated by polynomials over
$\Bbb Z_p$: for each $n\in\Bbb N$ there exists a polynomial $f_n(x)\in\Bbb Z_p[x]$,
such that $f(z)\equiv f_n(z)\pmod {p^n}$ for all $z\in\Bbb Z_p$. Actually,
the series $\sum_{j=0}^\infty r_j{x\choose j}$ defines a function, which
is identically 0 modulo $p^n$ iff all $r_j\equiv 0\pmod{p^n}$ (see \cite{11,
proposition 4.2}). So we may put $f_n(x)=\sum_{i=0}^{\omega(n)}a_i{x\choose i}$, where
$\omega (n)=\max\{j\in\Bbb N_0: {1\over p-1}(j-\wt_p\,j)<n\}$.      

It turnes out that the inverse assertion is also true: if a function 
$f:\Bbb Z_p\rightarrow\Bbb
Z_p$ is uniformly approximated by polynomials over $\Bbb Z_p$ in the above
mentioned sence, then it lies in $\Cal B$. To prove this assertion, assume
that $f(z)\equiv f_i(z)\pmod{p^i}$
for all $z\in\Bbb Z_p$, where  $f_i(x)\in\Bbb Z_p[x]$, $i=1,2,\ldots \ $.
Each polynomial $f_i(x)$ of degree $d_i$ admits one and the only representation 
as interpolation series $(\diamondsuit)$ of section 2:
$f_i(x)=\sum_{j=0}^{d_i}a_{ij}{x\choose j}$, where $a_{ij}\in\Bbb Z_p$ and
$\ord_p\,a_{ij}\ge\ord_p\,(j!)$ in view of $(\bigstar)$, since, obviously,
$f_i\in\Cal C\subset\Cal B$. For a given function $f$ each polynomial $f_i(x)$
is unique up to the summand which induces an identically 0 modulo $p^i$
function. So we may assume that $d_i=\omega (i)$ (see above); then coefficients
of the polynomial $f_i(x)$ are defined uniquelly up to the summands with
$p$-adic 
norms not exceeding $p^{-i}$. This implies that
$a_{i+1,j}\equiv a_{ij}\pmod{p^i}$ (we assume $a_{ij}=0$ for $j>\omega
(i)$). Hence, $\lim\limits^p_{i\to\infty}a_{ij}=a_j\in\Bbb Z_p$, and
${a_j\over j!}\in\Bbb Z_p$. Consequently, the series $\sum_{i=0}^\infty a_i{x\choose i}$
defines a uniformly continuous on $\Bbb Z_p$ function ${\tilde f}\in\Cal B$,
which must be equal to $f$, since $f(z)\equiv f_i(z)\equiv {\tilde f} (z)\pmod{p^i}$ 
for all $z\in\Bbb Z_p$
and all $i=1,2,\ldots \ $.

Now we define a non-Archimedian pseudo-valuation on $\Cal B$ as 
$\max\{\|f(z)\|_p\colon z\in\Bbb Z_p\}$ for  $f\in \Cal B$. The just proven results
imply that with respect to the distance $D_p$, induced by this pseudo-valuation, the
ring $\Cal B$ is a complete metric space; actually, $\Cal B$ is a completion
with respect to $D_p$ of the space  ${\Cal P}\subset{\Cal
C}$ of all functions induced on $\Bbb Z_p$ by polynomials over $\Bbb Z$
(in particular, the space $\Cal B$ is separable).

This implies, in turn, that $\Cal B$ (contrasting to $\Cal C$) is closed with respect to composition
of functions: if $f,g\in\Cal B$ then $f(g)\in\Cal B$. In fact, let $g$
be uniformly approximated by the sequence
$\{g_n(x)\in\Bbb Z_p[x]: n=1,2,\ldots\}$, that is, $g_n(z)\equiv g(z)\pmod{p^n}$
for all $z\in\Bbb Z_p$. The compatibility of the function $f$ imples then
that $D_p(f(g),f(g_n))\le  p^{-n}$, i.e., for $n\to\infty$ the sequence $f(g_n)$ tends to
$f(g)$ with respect to distance $D_p$. But $f(g_n)\in \Cal
B$ for each $n=1,2,\ldots$: if $f$ is uniformly approximated by the sequence
$\{f_m(x)\in\Bbb Z_p[x]: m=1,2,\ldots\}$, then 
$f_m(g_n(z))\equiv f(g_n(z))\pmod{p^m}$ for all $z\in\Bbb Z_p$. Hence,
the sequence $\{f_m(g_n(x))\in\Bbb Z_p[x]: m=1,2,\ldots\}$ 
tends to the function $f(g_n)$ with respect to the distance $D_p$, and
$f_m(g_n)\in\Cal B$, since it is a polynomial over $\Bbb Z_p$.
Consequently, $f(g)\in\Cal B$ in view of completeness of $\Cal
B$.

Thus, we have proven the following
\proclaim
{4.10 Proposition} The ring $\Cal B$ is a separable and complete with respect
to the distance $D_p$ metric space of functions, which are differentiable everywhere
on $\Bbb Z_p$. $\Cal B$ is closed with respect to compositions of functions
and with respect to differentiations. A countable set $\Cal P$ of all polynomials
over $\Bbb Z$ is a dence subset of $\Cal B$.
\qed
\endproclaim

To make use of criterion 4.9 for the applications to pseudorandom number
generation it is important to have a huge stock of examples of functions
of $\Cal B$ which are to be implemented as computer programs. As we have mentioned
above, all polynomials over $\Bbb Z_p$ are in $\Cal B$.

Rational over $\Bbb Z_p$ functions, that is, functions of the form
$f(x)=\frac{u(x)}{v(x)}$, where $u(x), v(x)\in\Bbb Z_p[x]$, are also in
$\Cal B$, providing  the denominator vanishes modulo $p$ nowhere on $\Bbb
Z_p$ (in view of compatibility it is sufficient to verify the latter condition only
for the points of $\{0,1,\ldots,p-1\}$). Indeed, for each $z\in\Bbb
Z_p$ the element $v(z)$  is not 0 modulo $p$, and hence has a multiplicative
inverse in the ring $\Bbb Z/p^n$. Thus
$\frac{u(z)}{v(z)}\equiv u(z)v(z)^{\phi(p^n)-1}\pmod{p^n}$, where $\phi$
is Euler totient function. Hence, the function $f$ could be uniformly approximated
by polynomials $u(x)v(x)^{\phi(p^n)-1}\in \Bbb Z_p[x]$, $n=1,2,\ldots \ $; 
hence, it is in  $\Cal B$ in force of 4.10.

Another type of functions of $\Cal B$ are exponential ones.
For instance, consider a function $a^x$ with $a\equiv 1\pmod
p$ (hence,  $a=1+pr$ for suitable $r\in\Bbb Z_p$). Then 
$a^x=\sum_{i=0}^{\infty}p^ir^i\binom{x}{i}$, and it is well known
(see e.g. \cite {3, Ch. 14, Section
5}), that for $p\ne 2$ this function is analytic on $\Bbb Z_p$ (hence,
lies in $\Cal C$). If $p=2$ and $r$ is odd, then $a^x$ is not analytic
on $\Bbb Z_2$, thus not in $\Cal C$. Newertheles, within the latter case
$a^x$ is in $\Cal B$, since $\ord_2\,(i!)=i-wt_2\,i$ and hence
$(1+2r)^x=\sum_{i=0}^\infty 2^ir^i{x\choose i}\in\Cal B$. It is not difficult
to see that the function $(1+4r)^x$ is in $\Cal C$. So, summarizing all
these considerations, if $a\in\Bbb Z_p$, $a\equiv
1\pmod p$, then the function $a^x$ is in $\Cal B$.

Exponential functions of the considered type are particular cases 
of functions
of more general form $u^{v}$, where $u(z)\equiv 1\pmod p$ for all $z\in\Bbb Z_p$.
\proclaim
{4.11 Lemma} Let  $u,{v}\colon\Bbb Z_p\rightarrow \Bbb Z_p$ be compatible functions
and let $u(z)\equiv 1\pmod p$ for all $z\in\Bbb Z_p$ {\rom (}{\rm so it is
sufficient
to verify the latter condition only for} $z=0,1,\ldots,p-1${\rom )}. Then
the function $f(z)=u(z)^{v(z)}$ is correctly defined for all $z\in\Bbb Z_p$,
integer-valued and compatible. Moreover, if $w,{v}\in\Cal B$, $u(z)=1+pw(z)$, 
then  $f\in\Cal B$.  
\endproclaim 
\demo 
{Proof} The above considerations of functions of type $a^x$ with 
$a\equiv 1\pmod p$
immediately imply that 
the function $f$ is correctly defined on $\Bbb Z_p$ and that it is integer-valued.
To prove the compatibility of $f$, note, that for arbitrary
$b,c,d\in\Bbb Z_p$ and $n=1,2,\ldots$ one has
$(a+p^nb)^{c+p^nd}=(a+p^nb)^{c}((a+p^nb)^{p^n})^{d}$, since elmentary properties
of powers are of the same form both in real and $p$-adic cases, see
\cite {3, Ch. 14, Section 5}. As both $u$ and $v$ are compatible functions,
then for arbitrary $z,r\in\Bbb Z_p$ there exist $s,t\in\Bbb Z_p$, such
that $(u(z+p^nr))^{v(z+p^nr)}=(u(z)+p^nt)^{v(z)+p^ns}$; hence
$(u(z+p^nr))^{v(z+p^nr)}= (u(z)+p^nt)^{v(z)}((u(z)+p^nt)^{p^n})^{s}\equiv
(u(z)+p^nt)^{v(z)}\pmod{p^n}$, in view of the congruence $(u(z)+p^nt)^{p^n}\equiv
1\pmod{p^n}$. The latter congruence is to be proven.

As $u(z)\equiv 1\pmod p$, then for a suitable $k\in\Bbb Z_p$ we have
$u(z)+p^nt=1+pk$. Yet $(1+pk)^{p^n}=\sum_{i=0}^{p^n}k^ip^i{p^n\choose i}=\sum_{i=0}^{p^n}k^i{p^i\over
i!}(p^n)_i\equiv 1\pmod{p^n}$, since ${p^i\over i!}\in\Bbb Z_p$. Finally,
denoting by $\overline {v(z)}=v(z)\bmod{p^n}$ the least nonnegative residue
of $v(z)$ modulo $p^n$, for a suitable $h\in\Bbb Z_p$ we obtain 
$f(z+p^nr)\equiv (u(z)+p^nt)^{v(z)}=(u(z)+p^nt)^{\overline{v(z)}}(u(z)+p^nt)^{p^nh}\equiv
(u(z)+p^nt)^{\overline{v(z)}}=\sum_{i=0}^{\overline{v(z)}} u(z)^{{\overline{v(z)}}-i}p^{ni}t^i{{\overline{v(z)}}\choose
i}=(u(z))^{\overline{v(z)}}\equiv (u(z))^{\overline{v(z)}}(u(z))^{p^nh}=(u(z))^{v(z)}$,
where $\equiv$ stands for congruence modulo $p^n$. Thus,
$f$ is compatible.

To prove the rest of the lemma, note, that for each $z\in\Bbb Z_p$  and
each $n=1,2,\ldots$ the congruence 
$(u(z))^{v(z)}\equiv \sum_{i=0}^n (u(z)-1)^i{v(z)\choose i}\pmod{p^n}$
holds, since $\|u(z)-1\|_p\le    {1\over p}$. This implies that
\roster 
\item 
all functions
$f_n=\sum_{i=0}^n {p^i\over i!}(v)_iw^i$ are in $\Cal B$, since all ${p^i\over i!}$
are $p$-adic integers (see above);
\item
the sequence $\{f_n: n=1,2,\ldots\}$ tends to
$f$ with respect to the distance $D_p$.
\endroster
Now (1)--(2) imply that $f\in\Cal B$ in force of 4.10.\qed
\enddemo

With the use of these results one may construct explicit forms of various ergodic functions
to be performed by a computer. For instance, the following is true.
\proclaim
{4.12 Proposition} For $g\in\Cal B$ the function $f(x)=1+x+p^2g(x)$ is
ergodic.
\endproclaim
\demo
{Proof} For $p\notin\{2,3\}$ the assertion trivially follows from 4.9.
For $p\in\{2,3\}$ in view of 4.9 it is sufficient to show that $f$ is transitive
modulo $p^3$. In turn, to demonstrate the latter it is sufficient to prove
only that $f^{kp^2}(0)\not\equiv 0\pmod{p^3}$ for $k=1,2,\ldots,p-1$, since
in $f$ is transitive modulo $p^2$ and hence, in view of its compatibility,
induces on $\Bbb Z/p^3$ a permutation with each cycle length  being a
multiple of $p^2$. Yet since for all $i=0,1,2\ldots $ the compatibility
of $g$ implies that $f^i(0)\equiv i+p^2\sum_{j=0}^{i-1}g(j)\pmod{p^3}$, 
then $f^{kp^2}(0)\equiv kp^2+p^2\sum_{j=0}^{kp^2-1}g(j)
\equiv kp^2+p^2\sum_{z=0}^{p-1}g(z)pk\equiv kp^2\pmod{p^3}$,
since (again in view of the compatibility of $g$) a congruence
$s\equiv r\pmod p$ implies the congruence
$p^2g(r)\equiv p^2g(s)\pmod{p^3}$.\qed
\enddemo

\head{5.} {Applications: a discussion.}\endhead

The results obtained in previous sections might have applications to design
pseudorandom number generators  which have relatively simple program implementation,
generate purely periodic sequences of numbers of $\{0,1,\ldots, m-1\}$
and provide certain guarantee for the statistical quality of these sequences,
their uniform distribution at the first turn. Speaking about relatively
simple program implementation, we mean that the considered generators have
certain parameters which are critical to the performance, and which one may vary
to achieve the desired performance without affecting the quality.

In case $m=p^k$ is a power of a prime $p$, these sequences might be generated
as the first order recurrence sequences satisfying the relation 
$x_{n+1}\equiv f(x_n)\pmod m$, where 
$f\colon\Bbb Z_p\rightarrow\Bbb Z_p$ is any compatible and ergodic
function of the considered in previous sections.
In this case for each $k=1,2,\ldots$ we obtain a purely periodic sequence
with period length $p^k$, with each element of $\{0,1,\ldots, p^k-1\}$
occuring at the period exactly once (in particular, the generated sequence
is uniformly distributed).

An important indicator of statistical quality of the sequence is
the distribution of $(r+1)$-tuples $\{(x_{n},x_{n+1},\ldots,x_{n+r}):
n=0,1,2,\ldots\}.$  Ideally, the sequence
$\{\bold u_n=({x_{n}\over p^k},{x_{n+1}\over p^k},\ldots,{x_{n+r}\over p^k}):n=0,1,2,\ldots\}$
of points of $(r+1)$-dimensional Euclidean space should be uniformly distributed
in the unit hypercube for all $r$. By no means this can be achieved
for periodic sequences. For such sequences there are some popular tests of 
quality, based on certain characteristics of families of hyperplanes,
which are parallel one to another, and
which union contain all the points corresponding to the sequences 
of $(r+1)$-tuples
(see e.g. \cite {2, section 3.3.4}).

Note, that  if for some $c,c_0,\ldots, c_{r}\in\Bbb Z$ the congruences
$$c+\sum_{i=0}^{r} c_i x_{n+i}\equiv
0\pmod{p^k}, \qquad (n=0,1,2,\ldots) \eqnum{\blacktriangle}$$  
hold, then  all the points $\bold u_n$ fall into the hyperplanes
$h+\sum_{i=0}^{r} c_i X_{i}=0$, which are parallel one to another.
For linear congruential generators such families of parallel hyperplanes
exist even for $r=2$, not depending on $k$ (see the introduction).

Note, that if $(\blacktriangle )$ holds for some $k$, then for all
$j=1,2,\ldots$  for the members of
the sequence $\{x_n\}$ hold relations $p^jc+\sum_{i=0}^{r} p^jc_i x_{n+i}\equiv
0\pmod{p^{k+j}}$. The relations of the latter kind will be temporarily 
and loosely
defined as
trivial. Trivial relations always exist: for instance, choosing certain
$K\in\Bbb N$ in view of the ergodicity of $f$ we obtain for all $k\ge K$
the trivial relations $p^{k-K}x_{n+p^K}\equiv p^{k-K}x_n\pmod {p^k}$.
Speaking informally, the triviality of relations just means that their
coefficients tend to 0 whereas $k$ tends to infinity, i.e. trivial  relations
are those which
degenerate to $0=0$ in $\Bbb Z_p$.

For an important wide class of nonlinear congruential generators  we
prove that if the dimension of hyperplanes, which are parallel one to
another, and which union contains all points $\bold u_n$, $(n=0,1,2,\ldots)$,
does not tend to infinity together with $k$, then this family of hyperplanes
is defined by trivial relations.

Now we give exact statements.
\proclaim
{5.1 Proposition} Let $f\in\Bbb Q_p[x]$ be an integer-valued, compatible
and ergodic polynomial of degree $d$ over a field $\Bbb Q_p$ of $p$-adic
numbers
{\rom(}{\rm all these polynomials for $p=2$ are completely characterized
by theorem 2.3; for odd $p$ see 2.4, 4.7 and a note preceding 4.7}{\rom)}. Let, further, 
$r$ be a positive integer rational such that for each $k\in\Bbb N$
there exist $c,c_0,\ldots, c_{r}\in\Bbb Z_p$, 
which satisfy $(\blacktriangle )$ and not all of which are $0$ modulo $p$.
Then $d=1$.
\endproclaim
We will need the following
\proclaim
{5.2 Lemma} Under the assumptions of proposition 5.1 let
$c,c_0,\ldots, c_{r}\in\Bbb Z_p$ be  not depending on $k$, that is, let there
exist $c,c_0,\ldots, c_{r}\in\Bbb Z_p$ satisfying $(\blacktriangle )$ for
all $k\in\Bbb N$ simultaneously. Then $d=1$.
\endproclaim
\demo
{Proof of the lemma 5.2} As $f$ is ergodic, then $d\ne 0$. Assume that $d>1$. 
Consider $w(x)=c+\sum_{i=0}^r c_i f^i(x)$. As $w(x)$ is a composition of
integer-valued and compatible polynomials over $\Bbb Q_p$, 
then $w(x)\in\Bbb Q_p[x]$ is integer-valued and compatible. Yet each
$f^i(x)$ has degree $d^i$; hence, since $d>1$, then $w(x)$, being a sum
of polynomials of pairwise distinct degrees, must be a polynomial
of nonzero degree. 

Yet, since $x_{n+i}\equiv f^i(f^n(x_0))\pmod {p^k}$,
the assumptions of the lemma imply that $w(x_n)\equiv
0\pmod{p^k}$ for all $n=0,1,2,\ldots$. In other words,
$w(z)\equiv 0\pmod {p^k}$ for all $z\in\Bbb Z_p$, since $x_n$ takes all
values in $\{0,1,\ldots,p^k-1\}$ in view of the ergodicity of $f$, and $w(x)$
is compatible. The assumptions of the lemma now imply that $w(z)\equiv 0\pmod {p^k}$
for all $z\in\Bbb Z_p$ and all $k=1,2,\ldots$. Consequently, $w(z)=0$ for
all $z\in\Bbb Z_p$ and hence polynomial $w(x)$ must be 0 in the ring $\Bbb Q_p[x]$. 
A contradiction
proving the lemma.
\qed\enddemo
\demo
{Proof of the proposition 5.1} By the assumption, for each $k\in\Bbb N$
the set $\Cal L_k$ of all $\bold c=(c,c_0,\ldots,c_r)\in\Bbb Z_p^{r+2}$, $\|\bold
c\|_p=1$ with $c,c_0,\ldots,c_r$ satisfying $(\blacktriangle )$, is not
empty. Obviously, $\Cal L_1\supset\Cal L_2\supset\ldots$, since $f$ is
compatible.

Further, we assert that each set $\Cal L_k$ is closed in the topology of metric space
$\Bbb Z_p^{r+2}$. Actually, if $\bold c\in\Cal L_k$,
$\bold c^{\prime}\in \Bbb Z_p^{r+2}$, $\|\bold c - \bold c^{\prime}\|\le
p^{-s}$, $s\ge k$, then $\bold c^{\prime}=\bold c+p^s \bold z$ for a suitable
$\bold z\in\Bbb Z_p^{r+2}$. Hence, $\|\bold c^{\prime}\|_p=1$ and
$\bold c^{\prime}$ satisfies $(\blacktriangle )$; consequently,
$\bold
c^{\prime}\in\Cal L_k$.

Now we apply to the sequence $\Cal L_1\supset\Cal L_2\supset\ldots$ the
$p$-adic analog of the lemma on the imbedded closed intervals of real analysis.
The analog of this lemma holds for the topological spaces
of much more general type ---
see e.g. the theorem in \cite{16, Ch. 3, section 34, I}, from which the $p$-adic 
case could be
easily deduced. Thus, we conclude
that the intersection of this sequence is not empty. That is, there
exists $\bold c^{\prime\prime}\in\Bbb
Z_p^{r+2}$ which satisfies the assumptions of lemma 5.2. Yet then $d=1$.\qed
\enddemo
From here we deduce the following
\proclaim
{5.3 Theorem} Let $f\in\Bbb Q_p[x]$ be an integer-valued compatible and
ergodic polynomial with $\deg f>1$, and let there exists $r\in\Bbb N$ such
that for each $k\in\Bbb N$ the linear complexity over the ring
$\Bbb Z/p^k$ of the recurrence sequence $\{x_n\}$, defined by
the recurrence relation $x_{n+1}\equiv f(x_n)\pmod{p^k}$, does not
exceed $r$. In other words, let there exist
$c^{(k)},c_0^{(k)},\ldots,c_r^{(k)}\in\Bbb Z_p $ such that
$$c^{(k)}+\sum_{i=0}^{r} c_i^{(k)} x_{n+i}\equiv
0\pmod{p^k}\qquad (n=0,1,2,\ldots). \eqnum{\blacktriangleleft}$$
Then $\lim\limits^p_{k\to\infty} c^{(k)}=\lim\limits^p_{k\to\infty} c^{(k)}_1=
\ldots =\lim\limits^p_{k\to\infty} c^{(k)}_r=0$.
\endproclaim
\demo
{Proof} To start with, we note, that from the proofs of both lemma 5.2
and proposition 5.1 it follows that they remain true if we let $k$
within their statements range over arbitrary infinite subset
of $\Bbb N$. 

Now for each $k\in\Bbb N$ choose (and fix) $c^{(k)},c^{(k)}_0,
c^{(k)}_1,\ldots,c^{(k)}_r\in\Bbb Z_p^{(r+2)}$ satisfying 
$(\blacktriangleleft)$. Put $\bold c_k=(c^{(k)},c^{(k)}_0,
c^{(k)}_1,\ldots,c^{(k)}_r)\in\Bbb Z_p^{(r+2)}$. In view of 5.1 then
$\|\bold c_k\|_p<1$ for all $k\in\Bbb N$. Denote $\Cal N =\{k\in\Bbb N:\|\bold
c_k\|_p>p^{-k}\}$. In other words, $k\notin \Cal N$ iff $(\blacktriangleleft)$
is equivalent to a congruence $0\equiv 0\pmod
{p^k}$. 

It is obvious that if $\Cal N$ is finite, then the conclusion
of the theorem is true. Let $\Cal N$ be infinite.

For $k\in\Cal N$ put $ \hat\bold c_k=
\|\bold c_k\|_p\bold c_k$ and denote  $\hat\Cal N$ a set of all $m\in\Bbb N$
such that $p^k\|\bold c_k\|_p=p^m$ for a suitable  $k\in\Cal
N$. In other words, we replace each  $(\blacktriangleleft)$ with the 
equivalent system of congruences
$$\hat c^{(k)}+\sum_{i=0}^{r} \hat c_i^{(k)} x_{n+i}\equiv
0\pmod{p^m}\qquad (n=0,1,2,\ldots),$$
where $(\hat c^{(k)},\hat c^{(k)}_0,
\hat c^{(k)}_1,\ldots,\hat c^{(k)}_r)=\hat\bold
c_k$, $p^m=p^k\|\bold c_k\|_p$.

If the set $\hat\Cal N$ is finite, the conclusion of the theorem is obviously
true. If $\hat\Cal N$ is infinite, then, since  $\|\hat\bold c_k\|_p=1$,
in view of 5.1 and the note at the beginning of the proof we conclude that
$\deg f=1$. A contradiction.\qed
\enddemo
In the statement of the theorem 5.3 we mention a notion of linear complexity
of a sequence over a ring. This is commonly used (especially in cryptography) 
characteristic of a quality of a sequence. Lemma 5.2 in these terms
asserts that the sequence $\{x_i=f(x_{i-1}):
i\in\Bbb N\}$ has infinite linear complexity over $\Bbb Z_p$,
providing $f\in\Bbb Q_p[x]$ is integer-valued compatible ergodic polynomial
of degree $d>1$. This assertion could be slightly strengthened.
\proclaim
{5.4 Corollary} If $f\in\Bbb Q_p[x]$ is an integer-valued compatible ergodic
polynomial of degree $d>1$, then a recurrence sequence $\{x_n\}$, which satisfy recurrence
relation $x_{n+1}=f(x_n)$, has infinite linear complexity over 
$\Bbb Q_p$.
\endproclaim
\demo
{Proof} If for suitable $c,c_0,\ldots,c_r\in\Bbb Q_p $, which are not 0
simultaneously, the equality $c+\sum_{j=0}^{r} c_jx_{n+j} =0$ holds for
all $n=0,1,2,\ldots$, then the equality $hc+\sum_{j=0}^{r} hc_jx_{n+j} =0$
with $h=1$, if $c,c_0,\ldots,c_r\in\Bbb Z_p$, and $h=\|(c,c_0,\ldots, c_r)\|_p$
otherwise, holds either. In view of compatibility of $f$ the conclusion
now follows from 5.2.\qed
\enddemo
\remark
{Note} The assumption $f\in\Bbb Q_p[x]$  within statements of 5.1--5.4
can not be omitted. For instance, let $p=2$ and let 
$$f(x)=1+x+4(-1)^{1+x}=
1+x+\sum_{j=0}^{\infty}(-1)^j2^{j+2}{x\choose j}.$$
By the theorem 2.3, the integer-valued function $f$ is compatible and ergodic.
However, it is easy to see that the recurrence sequence $\{x_n\in\Bbb Z_2\}$
with recurrence relation $x_{n+1}=f(x_n)$ satisfy the relation $x_{n+2}=x_n+2$,
i.e.,
has linear complexity 2 over $\Bbb Z_2$.
\endremark
 
We should notice that in this section we use the notion of linear complexity
of a sequence over a ring in a somewhat broader sence than it is commonly
used. More often the linear complexity of a sequence $\{x_n\}$ of elements
of a commutative ring $R$ is understood as the smallest $r>0$ such that there exist 
$c_0,\ldots,c_{r-1}\in R$ which satisfy simultaneously all equations
$x_{n+r}=\sum_{j=0}^{r-1} c_jx_{n+j}$ for $n=0,1,2,\ldots$. We, in distinction
from it,
admit nonzero constant term, as well as relations where all coefficients are zero
divisors (yet not all 0 simultaneously; in the assertion of 5.3 the latter,
however, is not important). If $R$ is a field, then both notions basically
do not differ one from another:
if a sequence satisfies a relation $c+\sum_{i=0}^{r} c_i x_{n+i}=
0$ with $c_r\ne 0$, then it satisfies the relation $x_{n+r+1}=c_r^{-1}c_0x_n-\sum_{j=0}^{r-1}c_r^{-1}(c_j-c_{j+1})x_{n+j+1}$.
Our definition seems to us some more convenient for geometric interpretations,
see above. 

In other words, we have shown that, loosely speaking, nonlinear ergodic
polynomial
generators are absolutely nonlinear --- the sequences they produce can
not be implemented as linear recurrences over $\Bbb Q_p$. 
We do not discuss here what is the impact of these results on testing
of the corresponding generators with the above mentioned statistical tests, leaving
this issue to the forthcoming paper. We only note that they give some evidence
that nonlinear congruential generators with integer-valued compatible ergodic
polynomials over $\Bbb Q$ as  state change functions in practice will pass 
the tests.

Properly restated analogs of 5.1--5.4 hold for composite $m=p_1^{k_1}\cdots p_s^{k_s}$,
which is a product of powers of distinct primes $p_1,\ldots, p_s$, providing
the transformation $f$ preserves all congruences of the ring $\Bbb Z_{p_1}\times \cdots\times \Bbb
Z_{p_s}$. In connection with congruential generators modulo a composite
$m$ we also note that one can take $f$ to be a function, defined on the
set $\Bbb N_0$ of all nonnegative integer rationals, which takes values
in $\Bbb Z$, preserves all congruences of the ring $\Bbb Z$ and which is
ergodic as a function of integer $p$-adic variable for all $p\in\{p_1,\ldots,p_s\}$.
These functions may also be constructed with the use of the results of
the paper.

For instance, such functions may be found in the class
$$\Cal B_0=\bigg\{\sum_{i=0}^{\infty}a_i\,(x)_i : a_i\in\Bbb Z; i=0,1,2,\ldots\bigg\},$$ 
where, we recall, $(x)_i$ is $i$th descending factorial power of $x$: 
$(x)_0=1$, $(x)_i=x(x-1)\cdots (x-i+1)$ for all $i=1,2,\ldots$. It is obvious
that $\Cal B_0$ is a proper subclass of the class $\Cal B$ (studied in section
4) for each prime $p$ (the definition of $\Cal B$, we recall, depends on
$p$). Since $\Cal B$ consists of functions, which preserve all congruences
of the ring $\Bbb Z_p$, then each function $g$ of $\Cal B_0$ preserves
all congruences of the ring $\Bbb Z$, that is, for each $a,b\in\Bbb N_0$
and each natural number $N>1$ a congruence $a\equiv b\pmod
N$ implies a congruence $g(a)\equiv g(b)\pmod N$. So as a state change function
of a pseudorandom generator we can take, for instance,
$$f(x)=1+x+p_1^2\cdots p_s^2g(x)\qquad (g\in\Cal B_0); $$ 
in view of 4.12 $f$ is ergodic as a function of integer $p_j$-adic variable
for all $j=1,2,\ldots,s$. That is, $f$ is transitive modulo $p_j^k$ for
all $k=1,2,\ldots$
and  for all $j=1,2,\ldots,s$. Thus, $f$ is transitive modulo
each $p_1^{t_1}\cdots p_s^{t_s}$, for arbitrary $t_1,\ldots,t_s\in\Bbb N$.
In particular, $f$ is transitive modulo $m$, and hence a pseudorandom number
generator with state change function $f$ and arbitrary initial state 
$x_0\in\{0,1,\ldots, m-1\}$ produces a purely periodic sequence
of period length $m$, and each number
of $\{0,1,\ldots, m-1\}$ occurs at the period of this sequence exactly once.

Obviously, $\Cal B_0$ contains all polynomials with rational integer coefficients,
so if $g(x)\in\Bbb Z[x]$ is a polynomial of degree $d\ge 1$, then the performance
of the correspondig pseudorandom generator is equivalent to $d$ additions
and
$d$ multiplications modulo $m$ of integer rationals. Obviosly, $\Cal B_0$
consists not only of polynomials over $\Bbb Z$. It is not clear, however,
whether it contains other `natural' functions which admit relatively simple
program implementation.

Moreover, if $m$ is arbitrary, it is not clear enough, which functions
should be considered as `natural', and which should not. If by `natural'
functions one understands  compositions of arithmetical operations (addition,
subtraction, multiplication, division, raising to a positive integer
power, exponentiation) then the functions of this kind could be constructed, for instance,
with the use of 2.3, 2.4, and 4.9 combined with
4.11 and 4.12. So, theorems 2.3--2.4 imply that a polynomial $f(x)\in\Bbb
Q[x]$ of a form
$$
f(x)=1+x+\sum^{d }_{i=0}c_{i}
\,p_1^{\left\lfloor \log_{p_1}(i+1)\right\rfloor+1}\cdots\,p_s^{\left\lfloor \log_{p_s}(i+1)\right\rfloor+1}{ {x}\choose {i}},
$$
\noindent for arbitrary $c_0, c_1, c_2 \ldots \in {\Bbb Z}$ is transitive
modulo arbitrary natural number  $M>1$, which is a product of powers of 
$\{p_1,\ldots,p_s\}$; in particular
$f(x)$ is transitive modulo $m$. Hence, the performance of the corresponding pseudorandom
generator is equivalent to    $d$ multiplications,
$d$ additions,  $d+1$ reductions some moduli and one division of integer 
rationals. 

From the above formula it follows that, for instance, a polynomial
$f(x)=1+x+{5\over 18}(x)_6$
is transitive modulo
$10^k$ for all $k=1,2,\ldots$. In a similar way, with the use of 2.5 and
4.11 (or 4.9 together with 4.11) one may construct generators which use
exponentiations. For instance, a function $f(x)=1+x+201^x$ (or, more generally, 
a function $f(x)=1+x+(1+200u(x))^{w(x)}$ with $u(x),v(x)\in\Bbb Z[x]$), 
as well as a function 
$f(x)=1+x+201^{201^x}$ are transitive modulo
$10^k$ for all  $k=1,2,\ldots$ (see 4.9 and 4.11); the same is true for
the function $f(x)=1+x+100\cdot
11^x$ (see 2.5 and 4.11). Judging by the number of publications on inversive
generators, taking a multiplicative inverse (or, generally, raising to
negative powers) modulo $m$ also should be considerd as
`natural' operations. Generators of this kind also could be constructed
with the use of results of the paper: for instance, taking $w(x)=-1, v(x)=x$
in the just mentioned example, one obtains a function
$f(x)=1+x+(1+200x)^{-1}$, which is transitive modulo $10^k$ for all $k=1,2,\ldots$.
 
 We note that during the past decade there were intensive studies of such
pseudorandom generators, as power
 generator ($f(x)=x^r, r\in\Bbb N$), exponential generator ($f(x)=a^x$),
twice exponential generator ($f(x)=a^{b^x}$) and inversive generator ($f(x)=a+bx^{-1}$ or
$f(x)=(a+bx)^{-1}$). The examples of generators, which are mentioned above
in the section, and which use compositions of arithmetical operations,
including exponentiation and raising to negative power, as we see, are somewhat
distinct from
the ones usually studied (by summand $1+x$, at the first turn).
These distinctions practically do not worsen the performance of the corresponding
programs. However, these distinctions do not allow
us to apply to the generators considered in this paper the results on already studied
generators. It would be very useful to  study the possibility of such transfer, 
since in this area there are a lot of important results
belonging to different authors (unfortunately, we could not present even a
short survey here because of hudge number of these).

At the same time, all the generators, introduced in this paper, are modulo given
$m$ equivalent to generators with recurrence relation 
$x_{n+1}\equiv f_m(x_n)\pmod m$, where $f_m(x)\in\Bbb
Q[x]$ (this is an immediate consequence of $p$-adic Weierstrass theorem,
for the latter
see e.g.
\cite{3, Ch. 10, Theorem 1}). Hence, all the results, obtained in literature
for so-called polynomial congruential generators, could be immediately applied
to generators, considered in this paper, at least, under extra restriction
$f_m(x)\in\Bbb Z[x]$.

We should note also, that a number of generators, studied in literature,
concern a specific case, when $m$ is a product of two distinct large primes.
The results of the current paper are of little interest for this particular case, since with
the use of these
results one can construct generators, which are either equivalent modulo
a prime divisor $p$ of $m$ to linear congruential generator, or involve some
given in advance 
transitive modulo $p$ polynomial of degree $>1$. The latter must be constructed beforehand
and then `adjusted' to make it transitive modulo some $p^s$, with $s$
satisfying assumptions of 3.14, 4.1 or 4.9. The methods of such `adjustment'
we hope to publish in one of forthcoming papers, and here we restrict ourselves
with an example. For instance, using these techniques, and choosing a transitive
modulo 5 polynomial 
$1+3x^3$, it is possible to construct a polynomial $1-127x-152x^3+152x^5$,
which is transitive modulo each $10^k$, with arbitrary $k=1,2,\ldots$.

So in view of these considerations, methods of construction of pseudorandom
generators, developed in the paper, could give the best effect if applied
to 
the case when $m$ is a product of relatively
small primes raised to relatively large powers. Thus the case $m=2^s$ is a natural
focuse point, being the easiest for program implementations, since the
reduction
of a positive integer rational
modulo $2^s$ is merely a truncation of all its 2-base expansion senior
bits, starting with the $s$th one (our numbering of digits starts with
0). It is this case, which leads to the most natural (judging by program
implementation) operations other than the above mentioned arithmetical
ones, namely, to bitwise logical operations like  $\XOR, \OR, \AND$ and other bitwise operations with nonnegative rational
integer operands, represented as 2-base expansions. And,
luckely, there is a complete description of measure-preserving (or ergodic)
functions in this case --- see section 2 of the paper.

The obtained results make it possible to construct pseudorandom
number generators, which satisfy some requirements to performance, statistics
or cryptographical security. This theme will be thoroughly studied in forthcoming
papers. Here we briefly note only that application of equiprobable functions,
which are also studied in the paper, as  output functions of congruential
generators with ergodic state change functions, allows us, preserving uniformity
of distribution, to eliminate one more well known disadvantage of congruential
generators, the so-called `low bit effect'. The latter demonstrates each sequence
$\{x_n\}$, satisfying recurrence relation $x_{n+1}\equiv f(x_n)\pmod{2^k}$
with compatible $f\colon\Bbb Z_2\rightarrow\Bbb Z_2$: a sequence,
composed of $j$th digits of each $x_n$, has a period length at most $2^{j+1}$ only.
Methods of remedy will be also studied in one of the future papers.

\enddocument
\end